\documentclass[11pt,reqno]{amsart}
\usepackage{graphicx}
\usepackage{amsmath,amssymb,mathrsfs,cite}
\usepackage[usenames]{color}
\usepackage[colorlinks=true,linkcolor=blue]{hyperref}
\usepackage{bm}
\usepackage{tikz}
\usepackage{graphicx}
\usepackage{exscale}
\usepackage{relsize}
\usepackage{verbatim}
\usepackage{enumerate}
\usepackage{booktabs,makecell,multirow}
\usepackage{booktabs}
\usepackage{extarrows}
\usepackage{enumitem}
\usepackage{subfigure}

\newtheorem{theorem}{Theorem}[section]

\newtheorem{lemma}[theorem]{Lemma}
\newtheorem{proposition}[theorem]{Proposition}

\theoremstyle{definition}

\theoremstyle{remark}
\newtheorem{remark}[theorem]{Remark}
\numberwithin{equation}{section}
\begin{document}
	\title[Factorization method for inverse elastic scattering]{Factorization method for inverse elastic cavity scattering}
	
	\author{Shuxin Li}
	\address{School of Science, Beijing University of Posts and Telecommunications,
		Beijing 100876, China.}
	\email{Lshuxin99@163.com}
	
	\author{Junliang Lv}
	\address{School of Mathematics, Jilin University, Changchun 130012, China.}
	\email{lvjl@jlu.edu.cn}
	
	\author{Yi Wang}
	\address{School of Mathematics, Jilin University, Changchun 130012, China.}
	\email{wangyi173@163.com}

	
	\begin{abstract}
		This paper is concerned with the inverse elastic scattering problem to determine the shape and location of an elastic cavity. By establishing a one-to-one correspondence between the Herglotz wave function and its kernel, we introduce the far-field operator which is crucial in the factorization method. We present a theoretical factorization of the far-field operator and rigorously prove the properties of its associated operators involved in the factorization.
		Unlike the Dirichlet problem where the boundary integral operator of the single-layer potential involved in the factorization of the far-field operator is weakly singular, the boundary integral operator of the conormal derivative of the double-layer potential involved in the factorization of the far-field operator with Neumann boundary conditions is hypersingular, which forces us to prove that this operator is isomorphic using Fredholm's theorem.
		Meanwhile, we present theoretical analyses of the factorization method for various illumination and measurement cases, including compression-wave illumination and compression-wave measurement, shear-wave illumination and shear-wave measurement, and full-wave illumination and full-wave measurement. In addition, we also consider the limited aperture problem and provide a rigorous theoretical analysis of the factorization method in this case. Numerous numerical experiments are carried out to demonstrate the effectiveness of the proposed method, and to analyze the influence of various factors, such as polarization direction, frequency, wavenumber, and multi-scale scatterers on the reconstructed results.
	\end{abstract}
	\maketitle
	
	\section{introduction}
	The inverse elastic scattering problems have attracted a lot of attention from scientists and engineers due to their diverse applications in many scientific areas \cite{ABGKLA15,K79}.
	For example, in geophysical exploration, one often detects the existence of faults, ore, or other geological features by analyzing the scattering and propagation of elastic waves through the geological landscape, which can be used for resource exploration and natural disaster hazard assessment. Another example is that in medical imaging, one can utilize the scattering of elastic waves to construct images of internal organs of the human body for diagnosis and monitoring of diseases.
	In contrast to acoustic and electromagnetic waves \cite{CK13,KRMK89,NLW23},
	elastic waves exhibit greater physical complexity due to the coupling of longitudinal and transverse wave. Therefore, it presents considerable mathematical and computational challenges for solving inverse elastic scattering problems.
	
	There have been extensive results in theoretical
	analyses and numerical methods for inverse elastic scattering problems of impenetrable elastic obstacles
	in the past decades.
	In theory, uniqueness results using infinitely many incident plane waves have been established in
	\cite{HKS12,HH93} with full-phase data and in
	\cite{CD24} by the reference ball technique with phaseless data. However, for the general shape of the elastic obstacles, the uniqueness result using one single incident wave is still a challenging open problem, although some researchers \cite{EH19,EY10,HKS12,LX10} gave uniqueness analysis under certain geometric assumptions. In numerics,
	the existing methods can be broadly classified into two categories: the quantitative method and the qualitative method. The formers includes the domain derivatives method \cite{LWWZ16}, the continuation method \cite{YLLY19}, and the iterative method \cite{DLL19,GM11}.
	The latter mainly refers to the sampling method and its variants, such as the linear sampling method \cite{A01}, the factorization method \cite{EH19,HKS12}, and the direct sampling method \cite{JLX18}.
	These methods have also been applied to a variety of inverse elastic scattering problems in waveguides \cite{BBC18,BLL11}, crack \cite{K96}, and periodic structures \cite{HLZ13}.

	Reconstructing algorithms is an important research topic in the field of inverse problems. To tackle the challenges caused by the nonlinearity and
	ill-posedness of the inverse scattering problem,
	the iterative method is a suitable approach to approximate the target. However, the iterative method exhibits local convergence and requires a priori information about the geometry and boundary condition of
	the problem; see \cite{LLW23,BLL13}. Most importantly,
	the iterative method generally requires the solution of forward scattering problems \cite{JLLWWZ22,JLLZ17,JLLW17,JLLW18} at each iteration step, which affects the computational efficiency. To avoid solving a series of forward problems, numerous non-iterative reconstruction
	techniques have been introduced, among which sampling methods are of great importance. Roughly speaking, sampling methods are based on choosing an appropriate indicator function to decide whether a point or line lies inside or outside the scatterer. This makes it easy to implement and relatively independent of the geometry and physical properties of the scatterer. Such advantages make it an efficient numerical method for determining the approximate location of an unknown scatterer. Subsequently, it can serve as a good initial guess in an iterative process for achieving a more accurate reconstruction of the scatterer; see \cite{ZS18}.
	
	Among various sampling methods, the factorization method, originally proposed by Kirsch \cite{K98} for inverse acoustic scattering problem, continues to be extensively studied. Compared with the linear sampling method \cite{CK96}, the factorization method better portrays the characteristics of points or lines inside or outside the scatterer in the numerical computation, and provides sufficient and necessary computation criterion for mathematical analysis \cite{KAG07}. Thus, the factorization method offers an approach to establish the uniqueness of inverse scattering problem. According to current literature, Alves and Kress \cite{AK02} were the first to establish the theoretical foundations of this method for the three-dimensional elastic scattering problem, and Arens \cite{A01} extended these theoretical results to two-dimensional rigid elastic scatterers and gave numerical calculations. After that, the factorization method has been systematically investigated for a wide range of
	possible elastic scatterer cases: obstacles \cite{EH19,HKS12}, penetrable bodies \cite{CKAGK06}, cracks \cite{GWY18}, periodic structures \cite{HLZ13},  or their combinations \cite{XY22,JY21}.

	In the present paper, we consider an inverse elastic
	scattering problem to determine the shape and location of a cavity, which satisfies the Neumann boundary condition. Such a problem has been investigated mathematically in \cite{GS12} for uniqueness  and numerically for shape reconstruction in \cite{GM11,YL23}
	using the iterative method and in \cite{JLX18}
	using the direct sampling method. In current paper, the factorization method is considered to solve this problem. We present a theoretical factorization of the far-field operator of an elastic cavity and rigorously prove the properties of the associated operator involved in the factorization. Unlike the Dirichlet problem, the boundary integral operator of the conormal derivative of the double-layer potential involved in the far-field factorization of the cavity is hypersingular, so this operator is not compact. The method used for the Dirichlet problem no longer holds. In order to deal with this issue, we use the Fredholm theorem to prove that this operator is isomorphic. Furthermore, we establish the theoretical foundations of the factorization method for a variety of illumination and measurement cases, in particular the limited aperture case. In addition, we give a possible reason why our analytical method does not extend directly to the impedance boundary condition case.
	
	The outline of this paper is as follows. In Section \ref{Section_2},
	we mathematically formulate
	the model of the inverse elastic scattering problem for a cavity and define our notations used throughout this paper. In Section \ref{Section_3}, the factorization form of the far field operator is derived and properties of the involved operators are presented.
	Section \ref{Section_4} is devoted to establishing the theoretical justification of the factorization method.
	Numerical examples are
	presented to illustrate the feasibility and effectiveness of our method in Section \ref{Section_5}.
	Eventually, the paper is concluded with some general
	remarks and directions for future research in Section \ref{Section_6}.

	\section{Model Formulation}\label{Section_2}
	Consider the elastic scattering by an impenetrable cavity $D$ with the Lipschitz boundary $\partial D$. The exterior domain $\mathbb{R}^2\setminus\overline{D}$ is filled with a homogeneous and isotropic elastic medium with Lam\'{e} constants $\lambda$ and $\mu$ satisfying $\mu>0$ and $\lambda+\mu>0$. Denote by $\bm{\nu}$ the unit outward normal vector to $\partial D$ and by $\bm{\tau}$  the unit tangential vector to $\partial D$. Then incident field $\bm{u}^{\mathrm{i}}$ is either generated by a compressional plane wave
	\begin{equation*}
		\bm{u}^{\mathrm{i}}(\bm{x},\bm{d};p)=\bm{d}
		e^{\mathrm{i}k_p\bm{x}\cdot\bm{d}},
	\end{equation*}
	or a shear plane wave
	\begin{equation*}
		\bm{u}^{\mathrm{i}}(\bm{x},\bm{d};s)=\bm{d}^\bot
		e^{\mathrm{i}k_s\bm{x}\cdot\bm{d}},
	\end{equation*}
	where $\bm{d}=(\cos\theta,\sin\theta)^\top$, $\bm{d}^\bot=(-\sin\theta,\cos\theta)^\top$, $\theta\in[0,2\pi]$ is the incident angle, and
	\begin{equation*}
		k_p=\frac{\omega}{\sqrt{\lambda+2 \mu}},\quad
		k_s=\frac{\omega}{\sqrt{\mu}},
	\end{equation*}
	are the compressional wavenumber and the shear wavenumber, respectively. We introduce notations used
	throughout the paper.  For a vector $\bm{x}\in\mathbb{R}^2$, we introduce two unit vector $\hat{\bm{x}}=\bm{x}/|\bm{x}|$ and $\bm{x}^\bot$
	obtained by rotating $\hat{\boldsymbol{x}}$ anticlockwise by $\pi/2$. We abbreviate $\mu \Delta \boldsymbol{u}+(\lambda+\mu) \operatorname{grad} \operatorname{div} \boldsymbol{u} $ by $\Delta^* \boldsymbol{u}$. The differential operator $\mathrm{grad}^\bot$ and $\mathrm{div}^\bot$ defined by
	$$
	\operatorname{grad}^{\perp} u:=\left[-\frac{\partial u}{\partial x_2}, \frac{\partial u}{\partial x_1}\right]^{\top},\quad \operatorname{div}^{\perp} \boldsymbol{u}:=\frac{\partial u_2}{\partial x_1}-\frac{\partial u_1}{\partial x_2}.
	$$

	Denote by $\bm{u}^{\mathrm{s}}$ the scattered field and by $\bm{u}=\bm{u}^{\mathrm{i}}+\bm{u}^{\mathrm{s}}$ the total field. Then the elastic scattering problem for cavity is to find the solution $\bm{u}^{\mathrm{s}}$ to
	the Navier equation
	\begin{equation}\label{navier equation}
		\mu \Delta \boldsymbol{u}^{\mathrm{s}}+(\lambda+\mu) \operatorname{grad} \operatorname{div} \boldsymbol{u}^{\mathrm{s}}+ \omega^2 \boldsymbol{u}^{\mathrm{s}}=\bm{0},\quad \text{in}~~\mathbb{R}^2\setminus\overline{D},
	\end{equation}
	which satisfies the Neumann boundary condition
	\begin{equation}\label{Neumann_BC}
		\bm{T}_{\bm{\nu}}\bm{u}^{\mathrm{s}}=-\bm{T}_{\bm{\nu}}
		\bm{u}^{\mathrm{i}},
		\quad \text{on}~~\partial D,
	\end{equation}
	where $\omega>0$ is the circular frequency, and the traction (conormal derivative) $\bm{T}_{\bm{\nu}}$ on $\partial D$ is given by
	\begin{equation*}
		\boldsymbol{T}_{\bm{\nu}} \boldsymbol{v}:=2 \mu \frac{\partial \boldsymbol{v}}{\partial \boldsymbol{\nu}}+\lambda \boldsymbol{\nu} \operatorname{div} \boldsymbol{v}-\mu \boldsymbol{\nu}^{\perp} \operatorname{div}^{\perp} \boldsymbol{v}.
	\end{equation*}
	In addition, to ensure uniqueness, the scattered field
	$\bm{u}^{\mathrm{s}}$ is required to satisfy the Kupradze radiation condition
	$$
	\lim _{r \rightarrow \infty} \sqrt{r}\left(\frac{\partial \bm{u}^{\mathrm{s}}_p}{\partial r}-\mathrm{i} k_p \bm{u}^{\mathrm{s}}_p\right)=0, \quad \lim _{r \rightarrow \infty} \sqrt{r}\left(\frac{\partial \bm{u}^{\mathrm{s}}_s}{\partial r}-\mathrm{i} k_s \bm{u}^{\mathrm{s}}_s\right)=0, \quad r=|\bm{x}|,
	$$
	uniformly in all directions. Here, the compressional wave $\bm{u}^{\mathrm{s}}_p$ and the shear wave $\bm{u}^{\mathrm{s}}_s$ of $\bm{u}^{\mathrm{s}}$ are defined by
	$$
	\bm{u}^{\mathrm{s}}_p:=-\frac{1}{k_p^2} \operatorname{grad} \operatorname{div} \bm{u}^{\mathrm{s}}, \quad \boldsymbol{u}^{\mathrm{s}}_s=-\frac{1} {k_{\mathrm{s}}^2}\operatorname{grad}^{\perp} \operatorname{div}^{\perp} \boldsymbol{u}^{\mathrm{s}}.
	$$
	The well-posedness of the forward scattering problem (\ref{navier equation})-(\ref{Neumann_BC}) can be found in \cite{ABGKLA15,BHSY18,KK99}.
	
	The fundamental solution to the Navier equation is given by
	\begin{equation}
		\mathbf{\Gamma}(\bm{x}, \bm{y}):=\frac{\mathrm{i}}{4 \mu} H_0^{(1)}\left(k_s|\bm{x}-\bm{y}|\right) \mathbf{I}+\frac{\mathrm{i}}{4 \omega^2} \operatorname{grad}_{\bm{x}} \operatorname{grad}_{\bm{x}}^{\top}\left[H_0^{(1)}\left(
		k_s|\bm{x}-\bm{y}|\right)
		-H_0^{(1)}\left(k_p|\bm{x}-\bm{y}|\right)\right],
	\end{equation}
	where $ H_0^{(1)}$ is the Hankel function of the first kind of order zero and $\mathbf{I}\in \mathbb{R}^{2 \times 2} $ is the identity matrix. Since the medium is homogeneous, one have $\boldsymbol{\Gamma}(\boldsymbol{x}, \boldsymbol{y})$ is symmetric and
	$$
	\bm{\Gamma}(\boldsymbol{y}, \boldsymbol{x})=\bm{\Gamma}(\boldsymbol{x}, \boldsymbol{y}), \quad \boldsymbol{x} \neq \boldsymbol{y}.
	$$
	The function $\mathbf{\Gamma}$ can also be decomposed into the compressional part $\mathbf{\Gamma}_p$ and shear part $\mathbf{\Gamma}_s$, i,e, $\mathbf{\Gamma}=\mathbf{\Gamma}_p+\mathbf{\Gamma}_s$.
	Moreover, $\bm{\Gamma}$ satisfies the Kupradze
	radiation condition. Thus, $\mathbf{\Gamma}$ has an  asymptotic behaviour of the from
	\begin{equation}\label{Gamma_asy}
		\boldsymbol{\Gamma}(\boldsymbol{x}, \boldsymbol{y})=\frac{e^{\mathrm{i} k_{p}|\boldsymbol{x}|}}{\sqrt{|\boldsymbol{x}|}} \boldsymbol{\Gamma}_{p}^{\infty}(\hat{\boldsymbol{x}}, \boldsymbol{y})+\frac{e^{\mathrm{i} k_{s}|\boldsymbol{x}|}}{\sqrt{|\boldsymbol{x}|}} \boldsymbol{\Gamma}_{s}^{\infty}(\hat{\boldsymbol{x}}, \boldsymbol{y})+\mathrm{O}\left(|\boldsymbol{x}|^{-3 / 2}\right), \quad|\boldsymbol{x}| \rightarrow \infty,
	\end{equation}
	uniformly in all directions $\hat{\boldsymbol{x}} \in \mathbb{S}$ where the functions $\boldsymbol{\Gamma}_{p}^{\infty}$ and $\boldsymbol{\Gamma}_{s}^{\infty}$ are known as the compressional part and the shear part of the far field patterns of $\boldsymbol{\Gamma}$, respectively.
	Here,
	\begin{equation}\label{Gamma}
		\begin{aligned}
			\boldsymbol{\Gamma}_{p}^{\infty}(\hat{\boldsymbol{x}}, \boldsymbol{y})=\frac{1}{\lambda+2 \mu} \frac{e^{\mathrm{i} \pi / 4}}{\sqrt{8 \pi k_{p}}} e^{-\mathrm{i} k_{p} \hat{\boldsymbol{x}} \cdot \boldsymbol{y}} \hat{\boldsymbol{x}} \otimes \hat{\boldsymbol{x}}, \quad
			\boldsymbol{\Gamma}_{s}^{\infty}(\hat{\boldsymbol{x}}, \boldsymbol{y})=\frac{1}{\mu} \frac{e^{\mathrm{i} \pi / 4}}{\sqrt{8 \pi k_{s}}} \mathrm{e}^{-\mathrm{i} k_{s} \hat{\boldsymbol{x}} \cdot \boldsymbol{y}} \hat{\boldsymbol{x}}^{\perp} \otimes \hat{\boldsymbol{x}}^{\perp}.
		\end{aligned}
	\end{equation}
	where '$\otimes$' denotes the tensor of two vectors.
	
	Using the Betti's formula, the integral representation of the radiating solution $\bm{u}^{\mathrm{s}}$ to Navier equation can be derived as
	\begin{equation}\label{integral re}
		\bm{u}^{\mathrm{s}}(\bm{x})=\int_{\partial D}\left\{\left[\bm{T}_{\bm{\nu}(\bm{y})} \mathbf{\Gamma}(\bm{x},\bm{y})\right]^{\top} \bm{u}^{\mathrm{s}}(\bm{y})-\mathbf{\Gamma}(\bm{x},\bm{y}) \bm{T}_{\bm{\nu}(\bm{y})} \bm{u}^{\mathrm{s}}(\bm{y})\right\} \mathrm{d} s(\bm{y}),
		\quad \bm{x} \in \mathbb{R}^2 \backslash \overline{D},
	\end{equation}
	where $\bm{T}_{\bm{\nu}}\mathbf{\Gamma}=(\bm{T}_{\bm{\nu}}
	\mathbf{\Gamma}_1
	,\bm{T}_{\bm{\nu}}\mathbf{\Gamma}_2)$ with $\mathbf{\Gamma}_j$ being the $j$th column of $\mathbf{\Gamma}$.
	It follows from (\ref{Gamma}) and (\ref{integral re}) that $\bm{u}^{\mathrm{s}}$ has an asymptotic behaviour of the form
	\begin{equation}\label{scattered_asy}
		\bm{u}^{\mathrm{s}}(\bm{x})=\frac{\mathrm{e}^{\mathrm{i} k_p|\bm{x}|}}{\sqrt{|\bm{x}|}}\bm{u}_p^{\infty} (\hat{\bm{x}})+\frac{\mathrm{e}^{\mathrm{i} k_s|\bm{x}|}}{\sqrt{|\bm{x}|}} \bm{u}_s^{\infty} (\hat{\bm{x}})+O\left(|\bm{x}|^{-3/2}\right), \quad|\bm{x}| \rightarrow \infty,
	\end{equation}
	uniformly in all directions $\hat{\bm{x}}$.
	The fields $\bm{u}_p^{\infty} (\hat{\bm{x}})$ and $\bm{u}_s^{\infty} (\hat{\bm{x}})$ are known as the compressional part and the shear part of far field pattern $\bm{u}^{\infty}$, respectively, and are given by
	\begin{equation}\label{ufar_p}
		\begin{aligned}
			\bm{u}_{p}^\infty(\hat{\bm{x}})&=\frac{1}{\lambda+2 \mu} \frac{\mathrm{e}^{\mathrm{i} \pi / 4}}{\sqrt{8 \pi k_{p}}} \int_{\partial D}\left\{\left[\bm{T}_{\bm{\nu}(\bm{y})}\left (\hat{\bm{x}} \hat{\bm{x}}^{\top} \mathrm{e}^{-\mathrm{i} k_p \hat{\bm{x}} \cdot \bm{y}}\right)\right]^{\top} \bm{u}^{\mathrm{s}}(\bm{y})-\hat{\bm{x}} \hat{\bm{x}}^{\top} \mathrm{e}^{-\mathrm{i} k_p \hat{\bm{x}} \cdot \bm{y}} \bm{T}_{\bm{\nu}(\bm{y})} \bm{u}^{\mathrm{s}}(\bm{y})\right\} \mathrm{d} s(\bm{y}), \\
			&=\frac{1}{\lambda+2 \mu} \frac{\mathrm{e}^{\mathrm{i} \pi / 4}}{\sqrt{8 \pi k_p}} \int_{\partial D}\left\{\boldsymbol{T}_{\bm{\nu}(\bm{y})} \left(e^{-\mathrm{i} k_p \hat{\bm{x}} \cdot \bm{y}} \hat{\bm{x}}\right) \cdot \bm{u}^{\mathrm{s}}(\bm{y})-\hat{\bm{x}} \cdot \bm{T}_{\nu(\bm{y})} \bm{u}^{\mathrm{s}}(\bm{y}) e^{-\mathrm{i} k_p \hat{\bm{x}} \cdot \bm{y}}\right\} \hat{\bm{x}} \mathrm{~d} s(\bm{y})
		\end{aligned}
	\end{equation}
	and
	\begin{equation}\label{ufar_s}
		\begin{aligned}
			\bm{u}_{s}^{\infty}(\hat{\bm{x}})&=\frac{1}{\mu} \frac{\mathrm{e}^{\mathrm{i} \pi / 4}}{\sqrt{8 \pi k_{\mathrm{s}}}} \int_{\partial D}\left\{\left[\bm{T}_{\bm{\nu}(\bm{y})}\left(\left(\mathbf{I}-\hat{\bm{x}} \hat{\bm{x}}^{\top}\right) \mathrm{e}^{-\mathrm{i} k_s \hat{\bm{x}} \cdot \bm{y}}\right)\right]^{\top} \bm{u}^{\mathrm{s}}(\bm{y})-\left(\mathbf{I}-\hat{\bm{x}} \hat{\bm{x}}^{\top}\right) \mathrm{e}^{-\mathrm{i} k_s \hat{\bm{x}} \cdot \bm{y}} \bm{T}_{\bm{\nu}(\bm{y})} \bm{u}^{\mathrm{s}}(\bm{y})\right\} \mathrm{d} s(\bm{y}),\\
			&=\frac{1}{ \mu} \frac{\mathrm{e}^{\mathrm{i} \pi / 4}}{\sqrt{8 \pi k_s}} \int_{\partial D}\left\{\boldsymbol{T}_{\bm{\nu}(\bm{y})} \left(e^{-\mathrm{i} k_s \hat{\bm{x}} \cdot \bm{y}} \hat{\bm{x}}^\bot\right) \cdot \bm{u}^{\mathrm{s}}(\bm{y})-\hat{\bm{x}}^\bot \cdot \bm{T}_{\bm{\nu}(\bm{y})} \bm{u}^{\mathrm{s}}(\bm{y}) e^{-\mathrm{i} k_s \hat{\bm{x}} \cdot \bm{y}}\right\} \hat{\bm{x}}^\bot \mathrm{~d} s(\bm{y})
		\end{aligned}
	\end{equation}
	for all $\hat{\bm{x}}\in\mathbb{S}$. The reduction of the two equations above uses $\hat{\bm{x}}\hat{\bm{x}}^\top=\hat{\bm{x}}\otimes\hat{\bm{x}}$ and $\mathbf{I}-\hat{\bm{x}} \hat{\bm{x}}^{\top}=\hat{\bm{x}}^\bot\otimes\hat{\bm{x}}^\bot$.
	From (\ref{integral re}), it can be seen that $\bm{u}_{p}^{\infty}$ and $\bm{u}_{s}^{\infty}$ are the
	far field patterns of the compressional wave $\bm{u}_{p}^{\mathrm{s}}$ and the shear wave $\bm{u}_{s}^{\mathrm{s}}$, respectively. We can further
	deduce the asymptotic behaviour of the from
	\begin{equation}\label{Tscattered_asy}
		\bm{T}_{\hat{\bm{x}}}\boldsymbol{u}^{\mathrm{s}}
		(\boldsymbol{x})=\frac{\mathrm{i}\omega^2}{k_p}
		\frac{e^{\mathrm{i} k_p|\boldsymbol{x}|}}{\sqrt{|\boldsymbol{x}|}} \boldsymbol{u}_p^{\infty}(\hat{\boldsymbol{x}})
		+\frac{\mathrm{i}\omega^2}{k_s}\frac{e^{\mathrm{i} k_s|\boldsymbol{x}|}}{\sqrt{|\boldsymbol{x}|}} \boldsymbol{u}_s^{\infty}(\hat{\boldsymbol{x}})
		+O\left(|\boldsymbol{x}|^{-3 / 2}\right), \quad|\boldsymbol{x}| \rightarrow \infty.
	\end{equation}
	
	From (\ref{ufar_p}) and (\ref{ufar_s}), it is easy to see that  $\bm{u}_{p}^{\infty}(\hat{\bm{x}})=
	u_{p}^{\infty}(\hat{\bm{x}})\hat{\bm{x}}$ and $\bm{u}_{s}^{\infty}(\hat{\bm{x}})=
	u_{s}^{\infty}(\hat{\bm{x}})\hat{\bm{x}}^\bot$. The far field patterns $\bm{u}^\infty$ can be rewritten as
	\begin{equation*}
		\bm{u}^\infty(\hat{\bm{x}})=u_{p}^{\infty}(\hat{\bm{x}})\hat{\bm{x}}
		+u_{s}^{\infty}(\hat{\bm{x}})\hat{\bm{x}}^\bot.
	\end{equation*}
	Thus, we define $(u_{p}^{\infty},u_{s}^{\infty})$ to be the far field pair of $\bm{u}^\infty$ corresponding to the scattered field $\bm{u}^{\mathrm{s}}$.
	The fields $\bm{u}^{\mathrm{s}}(\bm{x},\bm{d};p)$, $\bm{u}^\infty(\hat{\bm{x}},\bm{d};p)$,
	and $(u_p^\infty(\hat{\bm{x}},\bm{d};p),
	u_s^\infty(\hat{\bm{x}},\bm{d};p))$
	represent the scattered field, the far field patterns and the far field pair corresponding to an
	incident compressional plane wave $\bm{u}^{\mathrm{i}}(\bm{x},\bm{d};p)$, respectively; the fields $\bm{u}^{\mathrm{s}}(\bm{x},\bm{d};s)$,
	$\bm{u}^\infty(\hat{\bm{x}},\bm{d};s)$,
	and $(u_p^\infty(\hat{\bm{x}},\bm{d};s),
	u_s^\infty(\hat{\bm{x}},\bm{d};s))$
	represent the scattered field, the far field patterns and the far field pair corresponding to an
	incident shear plane wave $\bm{u}^{\mathrm{i}}(\bm{x},\bm{d};s)$, respectively.
	
	One can use a superposition of two plane waves $\bm{u}^{\mathrm{i}}(\bm{x},\bm{d})=
	\bm{u}^{\mathrm{i}}(\bm{x},\bm{d};p)+
	\bm{u}^{\mathrm{i}}(\bm{x},\bm{d};s)$ as the incident field. From the linearity of direct scattering problem, the corresponding scattered field is given by
	\begin{equation*}
		\bm{u}^{\mathrm{s}}(\bm{x},\bm{d})=\bm{u}^{\mathrm{s}}
		(\bm{x},\bm{d};p)+
		\bm{u}^{\mathrm{s}}(\bm{x},\bm{d};s)
	\end{equation*}
	and has the far field pattern
	\begin{equation*}
		\begin{aligned}
			\bm{u}^\infty(\hat{\bm{x}},\bm{d})&=
			\bm{u}^\infty(\hat{\bm{x}},\bm{d};p)
			+\bm{u}^\infty(\hat{\bm{x}},\bm{d};s)\\
			&=\left(u_p^\infty(\hat{\bm{x}},\bm{d};p)+
			u_p^\infty(\hat{\bm{x}},\bm{d};s)\right)
			\hat{\bm{x}}
			+\left(u_s^\infty(\hat{\bm{x}},\bm{d};p)+
			u_s^\infty(\hat{\bm{x}},\bm{d};s)\right)
			\hat{\bm{x}}^\bot.
		\end{aligned}
	\end{equation*}
	
	The inverse problems involved in the present paper can be stated as follows:
	\vskip0.2cm
	\textbf{FF case}: Determine $\partial D$ from the knowledge
	of the far field pattern $\bm{u}^\infty(\hat{\bm{x}},\bm{d})$ for all incidence directions $\bm{d}\in \mathbb{S}$ and all observation directions $\hat{\bm{x}}\in\mathbb{S}$ corresponding to the incident plane wave $\bm{u}^{\mathrm{i}}(\bm{x},\bm{d})$.
	\vskip0.2cm
	\textbf{PP case}: Determine $\partial D$ from the knowledge
	of the compressional part $\bm{u}_p^\infty(\hat{\bm{x}},\bm{d};p)$ for all incidence directions $\bm{d}\in \mathbb{S}$ and all observation directions $\hat{\bm{x}}\in\mathbb{S}$ corresponding to the incident plane compressional wave $\bm{u}^{\mathrm{i}}(\hat{\bm{x}},\bm{d};p)$.
	\vskip0.2cm
	\textbf{SS case}: Determine $\partial D$ from the knowledge
	of the shear part $\bm{u}_s^\infty(\hat{\bm{x}},\bm{d};s)$ for all incidence directions $\bm{d}\in \mathbb{S}$ and all observation directions $\hat{\bm{x}}\in\mathbb{S}$ corresponding to the incident plane shear wave $\bm{u}^{\mathrm{i}}(\hat{\bm{x}},\bm{d};s)$.
	\vskip0.2cm
	\textbf{Limited aperture case}: Determine $\partial D$ from the knowledge
	of the far field pattern $\bm{u}^\infty(\hat{\bm{x}},\bm{d})$ for all incidence directions $\bm{d}\in\mathbb{U}$ and partial observation directions $\hat{\bm{x}}\in\mathbb{U}$ corresponding to the incident plane wave $\bm{u}^{\mathrm{i}}(\bm{x},\bm{d})$.
	Here, $\mathbb{U}$ is a non-empty subset of $\mathbb{S}$.
	\vskip0.2cm
	We conclude this section with an introduction to the Hilbert space used in the following. Denote by $[L^2(\partial D)]^2$ and $[H^{1/2}(\partial D)]^2$ the usual Sobolev space of vector field on $\partial D$. $[H^{-1/2}(\partial D)]^2$ is the dual space $[H^{1/2}(\partial D)]^2$ with respect to the inner product in $[L^2(\partial D)]^2$. $[H^1_{\mathrm{loc}}(\mathbb{R}^2\backslash\overline{D})]^2$
	is the space of all functions $\bm{u}:\mathbb{R}^2\backslash\overline{D}
	\rightarrow\mathbb{C}^2$ such that $\bm{u}\in[H^1((\mathbb{R}^2\backslash\overline{D})\cap B)]^2$ for all open balls $B$ containing $\overline{D}$.
	Let $\mathbb{S}:=\left\{\boldsymbol{x} \in \mathbb{R}^2:|\boldsymbol{x}|=1\right\}$ denote the unit disk in $\mathbb{R}^2$. We define
	\begin{equation*}
		\mathcal{L}^2(\mathbb{S}):=\left\{(g_p,g_s): \bm{g}(\bm{d})=g_p(\bm{d})\bm{d}
		+g_s(\bm{d})\bm{d}^\bot, g_p\in L^2(\mathbb{S}), g_s\in L^2(\mathbb{S}),\bm{d}\in\mathbb{S}
		\right\},
	\end{equation*}
	equipped with the inner product
	$$
	\langle\boldsymbol{g}, \boldsymbol{h}\rangle:=\frac{\omega}{k_p} \int_{\mathbb{S}} g_{p}(\boldsymbol{d}) \overline{h_p(\boldsymbol{d})} \mathrm{d} s(\boldsymbol{d})+\frac{\omega}{k_s} \int_{\mathbb{S}} g_{s}(\boldsymbol{d}) \overline{h_{s}(\boldsymbol{d})} \mathrm{d} s(\boldsymbol{d}), \quad \boldsymbol{g}, \boldsymbol{h} \in \mathcal{L}^2(\mathbb{S}).
	$$
	
	\section{Properties of the operator}\label{Section_3}
	In this section, we present the factorization form of the far field operator and discuss the properties of the involved operators to lay a foundation for the theoretical analysis of the factorization method.
	\subsection{Far field operator}
	Given the vector $\bm{g}\in \mathcal{L}^2({\mathbb{S}})$, the superposition of plane waves
	$$
	\boldsymbol{v}^{\mathrm{i}}(\boldsymbol{x}):=\mathrm{e}^{-\frac{\mathrm{i} \pi}{4} } \int_{\mathbb{S}}\left\{\sqrt{\frac{k_p}{\omega}} \boldsymbol{d}e^{\mathrm{i} k_p \boldsymbol{x} \cdot \boldsymbol{d}}  g_{p}(\boldsymbol{d})+\sqrt{\frac{k_s}{\omega}} \boldsymbol{d}^{\perp}e^{\mathrm{i} k_{s} \boldsymbol{x} \cdot \boldsymbol{d}} g_s(\boldsymbol{d})\right\} \mathrm{d} s(\boldsymbol{d}), \quad \boldsymbol{x} \in \mathbb{R}^2,
	$$
	is called an elastic Herglotz wave function with the Herglotz kernel $\boldsymbol{g}$. Herglotz wave functions are clearly entire solutions to the Navier equation. We note that for a given $\bm{g}\in \mathcal{L}^2(\mathbb{S})$ the function
	\begin{equation*}
		\boldsymbol{v}_{\bm{g}}(\boldsymbol{x}):=e^{-\frac{\mathrm{i} \pi}{4}} \int_{\mathbb{S}}\left\{\sqrt{\frac{k_p}{\omega}} \boldsymbol{d}\mathrm{e}^{-\mathrm{i} k_p \boldsymbol{x} \cdot \boldsymbol{d}} g_{p}(\boldsymbol{d})+\sqrt{\frac{k_s}{\omega}} \boldsymbol{d}^{\perp}\mathrm{e}^{-\mathrm{i} k_s\boldsymbol{x} \cdot \boldsymbol{d}} g_{s}(\boldsymbol{d})\right\} \mathrm{d} s(\boldsymbol{d}), \quad \boldsymbol{x} \in \mathbb{R}^2,
	\end{equation*}
	also is a Herglotz wave function. The following lemma establishes a one to one correspondence between Herglotz
	wave functions and their kernels.

	\begin{lemma}\label{g theorem}
		Assume that the Herglotz wave function $\bm{v}_{\bm{g}}$ with kernel $\bm{g}$ vanishes in all of $\mathbb{R}^2$. Then $\bm{g}=\bm{0}$.
	\end{lemma}
	\begin{proof}
		From $\boldsymbol{v}_{\bm{g}}=\boldsymbol{0}$ for all $\mathbb{R}^2$, one can obtain $\bm{z}\cdot \boldsymbol{v}_{\bm{g}}=0$ where $\bm{z}\neq 0$ is a constant vector. We use the polar coordinate forms $\bm{x}=(r\cos\varphi,r\sin\varphi)$ and
		$\bm{d}=(\cos\theta,\sin\theta)$. There are two cases.
		
		\textbf{Case 1:} When $k_p\neq k_s$, we have
		\begin{equation*}
			\sqrt{\frac{k_p}{\omega}} \int_\mathbb{S}\bm{z}\cdot \boldsymbol{d}g_p(\bm{d}) \int_0^{2\pi} e^{-\mathrm{i} k_p \bm{x} \cdot \bm{d}} e^{\mathrm{i} n \varphi} \mathrm{d} \varphi \mathrm{d} s(\bm{d}) + \sqrt{\frac{k_s}{\omega}}
			\int_\mathbb{S} \bm{z}\cdot \boldsymbol{d}^\bot g_s(\bm{d}) \int_0^{2\pi} e^{-\mathrm{i}k_s \bm{x} \cdot \bm{d}} e^{\mathrm{i} n \varphi} \mathrm{d} \varphi \mathrm{d} s(\bm{d})=0.
		\end{equation*}
		It follows from Funk-Hecke formula in two dimension that
		\begin{equation}\label{Funk_equation}
			\begin{aligned}
				\sqrt{\frac{k_p}{\omega}}\frac{2 \pi}{\mathrm{i}^n} J_n(k_pr) \int_\mathbb{S}\bm{z}\cdot \boldsymbol{d}g_p(\bm{d}) e^{-\mathrm{i}n \theta}
				\mathrm{d}s(\bm{d}) + \sqrt{\frac{k_s}{\omega}}\frac{2 \pi}{\mathrm{i}^n} J_n(k_sr)
				\int_\mathbb{S} \bm{z}\cdot \boldsymbol{d}^\bot g_s(\bm{d})e^{-\mathrm{i} n \theta} \mathrm{d} s(\bm{d})=0.
			\end{aligned}
		\end{equation}
		Taking the derivative of equation (\ref{Funk_equation}) with respect to $r$ yields
		\begin{equation*}
			\begin{aligned}
				\sqrt{\frac{k_p}{\omega}}\frac{2 \pi}{\mathrm{i}^n} k_pJ'_n(k_pr) \int_\mathbb{S}\bm{z}\cdot \boldsymbol{d}g_p(\bm{d}) e^{-\mathrm{i} n \theta}
				\mathrm{d}s(\bm{d}) + \sqrt{\frac{k_s}{\omega}}\frac{2 \pi}{\mathrm{i}^n} k_sJ'_n(k_sr)
				\int_\mathbb{S} \bm{z}\cdot \boldsymbol{d}^\bot g_s(\bm{d})e^{-\mathrm{i} n \theta} \mathrm{d} s(\bm{d})=0.
			\end{aligned}
		\end{equation*}
		Since the determinant
		\begin{equation*}
			\left|\begin{array}{ll}
				J_n\left(k_p r\right) & J_n\left(k_s r\right) \\
				k_p J_n^{\prime}\left(k_p r\right) & k_s J_n^{\prime}\left(k_s r\right)
			\end{array}\right| \neq 0,
		\end{equation*}
		we see that
		\begin{equation}\label{zdg}
			\int_\mathbb{S}\bm{z}\cdot \boldsymbol{d} g_p(\bm{d})e^{-\mathrm{i} n \theta}\mathrm{d} s(\bm{d})=0,\quad \int_\mathbb{S} \bm{z}\cdot \boldsymbol{d}^\bot g_s(\bm{d}) e^{-\mathrm{i} n \theta} \mathrm{d} s(\bm{d})=0,
		\end{equation}
		for all basis functions $e^{-in\theta}$ of order $n=0,1, \cdots$. It is follows from the the completeness of basis functions in $L^2(\mathbb{S})$ that $\bm{z}\cdot \boldsymbol{d} g_p=0$ and $\bm{z}\cdot \boldsymbol{d}^\bot g_s=0$. Hence, $\boldsymbol{d}\boldsymbol{g}_p=0$ and $\boldsymbol{d}^\bot\boldsymbol{g}_s=0$ are obtained by
		the arbitrariness of $\bm{z}$.
		
		\textbf{Case 2:} When $k_p=k_s$, we let $k=k_p=k_s$.
		The form of $\bm{v}_{\bm{g}}$ can be written as
		\begin{equation*}
			\boldsymbol{v}_{\bm{g}}(\boldsymbol{x})=e^{-\frac{\mathrm{i} \pi}{4}} \sqrt{\frac{k}{\omega}}\int_{\mathbb{S}} e^{-\mathrm{i} k\boldsymbol{x}\cdot\bm{d}}\left( \boldsymbol{d}g_{p}(\boldsymbol{d})+ \boldsymbol{d}^{\perp}g_{s}(\boldsymbol{d})\right) \mathrm{d} s(\boldsymbol{d})=e^{-\frac{\mathrm{i} \pi}{4}} \sqrt{\frac{k}{\omega}}\int_{\mathbb{S}} e^{-\mathrm{i} k\boldsymbol{x}\cdot\bm{d}}\bm{g}(\bm{d}) \mathrm{d} s(\boldsymbol{d}), \quad \boldsymbol{x} \in \mathbb{R}^2,
		\end{equation*}
		Similarly, we have
		\begin{equation*}
			\int_\mathbb{S}\bm{z}\cdot \boldsymbol{g}(\bm{d}) \int_0^{2\pi} e^{-\mathrm{i} k \bm{x} \cdot \bm{d}} e^{\mathrm{i} n \varphi} \mathrm{d} \varphi \mathrm{d} s(\bm{d}) =0.
		\end{equation*}
		From the Funk-Hecke formula in two dimension again, one can get
		\begin{equation*}
			\begin{aligned}
				\frac{2 \pi}{\mathrm{i}^n} J_n(k_pr) \int_\mathbb{S}\bm{z}\cdot \boldsymbol{g}(\bm{d}) e^{-\mathrm{i} n \theta}
				\mathrm{d}s(\bm{d})=0.
			\end{aligned}
		\end{equation*}
		Using the same treatment as (\ref{zdg}), we obtain $\bm{g}=\bm{0}$.
		
	\end{proof}
	Similarly to the acoustic Herglotz wave function, we have the following lemma.
	\begin{lemma}\label{Herglotz theorem}
		For a given function $\bm{g}\in\mathcal{L}^2(\mathbb{S})$ the solution of the exterior Neumann scattering problem for the incident wave
		\begin{equation*}
			\boldsymbol{v}^{\mathrm{i}}(\boldsymbol{x})=
			e^{-\frac{\mathrm{i} \pi}{4}} \int_{\mathbb{S}}\left\{\sqrt{\frac{k_p}{\omega}} e^{\mathrm{i} k_p \boldsymbol{x} \cdot \boldsymbol{d}} \boldsymbol{d} g_p(\boldsymbol{d})+\sqrt{\frac{k_s}{\omega}} e^{\mathrm{i} k_s \boldsymbol{x} \cdot \boldsymbol{d}} \boldsymbol{d}^{\perp} g_s(\boldsymbol{d})\right\} \mathrm{d} s(\boldsymbol{d}), \quad \boldsymbol{x} \in \mathbb{R}^2,
		\end{equation*}
		is given by
		\begin{equation*}
			\boldsymbol{v}^{\mathrm{s}}(\boldsymbol{x})=
			e^{-\frac{\mathrm{i} \pi}{4}}
			\int_{\mathbb{S}}\sqrt{\frac{k_p}{\omega}}
			\boldsymbol{v}^{\mathrm{s}}(\bm{x};\bm{d},\bm{d})
			g_d(\boldsymbol{d} )  \mathrm{d} s(\boldsymbol{d} )+\int_{\mathbb{S}} \sqrt{\frac{k_s}{\omega}} \boldsymbol{v}^{\mathrm{s}}(\bm{x};\bm{d},\bm{d}^\bot) g_s(\boldsymbol{d} )  \mathrm{d} s(\boldsymbol{d} ),\quad \bm{x}\in \mathbb{R}^2\setminus\overline{D},
		\end{equation*}
		and has the far-field patterns
		\begin{equation*}
			\begin{aligned}
				& \bm{v}^{\infty}_p(\hat{\bm{x}})=e^{-\frac{\mathrm{i} \pi}{4}}\int_{\mathbb{S}} \sqrt{\frac{k_p}{\omega}}\bm{v}^{\infty}_p(\hat{\bm{x}};
				\bm{d}, \bm{d}) g_d(\bm{d})  \mathrm{d} s(\bm{d})
				+\int_{\mathbb{S}}\sqrt{\frac{k_s}{\omega}} \bm{v}^{\infty}_p(\hat{\bm{x}};\bm{d}, \bm{d}^\bot) g_s(\bm{d})  \mathrm{d} s(\bm{d}), \quad \hat{\bm{x}}\in \mathbb{S},\\
				& \bm{v}_s^{\infty}(\hat{\bm{x}})=e^{-\frac{\mathrm{i} \pi}{4}}\int_{\mathbb{S}}\sqrt{\frac{k_p}{\omega}} \bm{v}^{\infty}_s(\hat{\bm{x}};\bm{d}, \bm{d}) g_d(\bm{d})  \mathrm{d} s(\bm{d})
				+\int_{\mathbb{S}}\sqrt{\frac{k_s}{\omega}} \bm{v}^{\infty}_s(\hat{\bm{x}};\bm{d}, \bm{d}^\bot) g_s(\bm{d}) \mathrm{d} s(\bm{d}),\quad \hat{\bm{x}}\in \mathbb{S}.
			\end{aligned}
		\end{equation*}
	\end{lemma}
	\begin{proof}
		The arguments are similar to those used in \cite{CK13} for the Helmholtz equation cases and in \cite{KS01}
		for the Navier equation cases.
	\end{proof}
	
	On the foundation of the above Lemma \ref{Herglotz theorem}, we introduce
	the far field operator
	$\boldsymbol{\mathcal{F}}: \mathcal{L}^2(\mathbb{S}) \longrightarrow \mathcal{L}^2(\mathbb{S}) $
	\begin{equation*}
		\boldsymbol{\mathcal{F} g}(\hat{\boldsymbol{x}}) =e^{-\frac{\mathrm{i} \pi}{ 4} } \int_\mathbb{S}\left\{\sqrt{\frac{k_{p}}{\omega}} \boldsymbol{u}^{\infty}(\hat{\boldsymbol{x}}, \boldsymbol{d};p) g_{\mathrm{p}}(\boldsymbol{d})+\sqrt{\frac{k_{s}}{\omega}} \boldsymbol{u}^{\infty}(\hat{\boldsymbol{x}}, \boldsymbol{d};s) g_{\mathrm{s}}(\boldsymbol{d})\right\} \mathrm{d} s(\boldsymbol{d}),\quad \hat{\bm{x}}\in \mathbb{S},
	\end{equation*}
	it is noted that $\boldsymbol{\mathcal{F}g}$ is the far-field pattern corresponding to the scattered field generated by the elastic Herglotz wave function with kernel $\boldsymbol{g}$ as incident field. The far field operator $\boldsymbol{\mathcal{F}}$ has the following properties (see \cite{A01,AK02}):
	
	\begin{itemize}
		\item [(a)]$\bm{\mathcal{F}}$ is injective if $\omega^2$ is not a Neumann eigenvalue of $-\Delta^*$ in $D$.
		\item [(b)]$\boldsymbol{\mathcal{F}}$ is compact and normal.
		\item [(c)]$\boldsymbol{\mathcal{F}}$ has a countable infinite number of eigenvalues. These eigenvalues lie on the disc with center at $(0,\sqrt{2 \pi / \omega} )$ on the positive imaginary axis and radius $\sqrt{2 \pi / \omega}$.
	\end{itemize}

	\begin{theorem}\label{F_Fac}
		The far field operator $\boldsymbol{\mathcal{F}}: \mathcal{L}^2(\mathbb{S}) \rightarrow \mathcal{L}^2(\mathbb{S})$ has a factorization in the form
		\begin{equation}\label{F factori}
			\boldsymbol{\mathcal{F}}=-\sqrt{8\pi\omega}\boldsymbol
			{\mathcal{G}}
			\boldsymbol{\mathcal{N}}^* \boldsymbol{\mathcal{G}}^*
		\end{equation}
		where the data-to-pattern operator $\boldsymbol{\mathcal{G}}:\left[H^{-1 / 2}(\partial D)\right]^2 \rightarrow \mathcal{L}^2(\mathbb{S}) $ maps $\bm{f} \in [H^{-1 / 2}(\partial D)]^2$ into the far field pattern $\bm{v}^{\infty}=\bm{\mathcal{G}} \bm{f}$ of the exterior Neumann boundary value problem with boundary data $\bm{f}$, and $\boldsymbol{\mathcal{N}}: [H^{1 / 2}(\partial D)]^2 \rightarrow [H^{-1 / 2}(\partial D)]^2$ is the boundary integral operator, defined by
		\begin{equation}\label{N oper}
			(\boldsymbol{\mathcal{N}}\bm{\varphi})(\bm{x})=\bm{T}_{\bm{\nu}(\bm{x})} \int_{\partial D}
			[\bm{T}_{\bm{\nu}(\bm{y})}\bm{\Gamma}(\bm{x},\bm{y})]^{\top}
			\bm{\varphi}(\bm{y}) d s(\bm{y}), \quad \bm{x} \in \partial D,
		\end{equation}
		for $\boldsymbol{\varphi} \in [H^{1 / 2}(\partial D)]^2$.
	\end{theorem}

	\begin{proof}
		We define the operator $\bm{\mathcal{H}}: \mathcal{L}^2(\mathbb{S}) \rightarrow\left[H^{-1 / 2}(\partial D)\right]^2$ by
		\begin{equation*}
			\begin{aligned}
				\bm{\mathcal{H}} \boldsymbol{g}(\boldsymbol{x}):
				=e^{-\frac{\mathrm{i} \pi}{4}}\left\{\int_{\mathbb{S}} \sqrt{\frac{k_{p}}{\omega}} \bm{T}_{\bm{\nu}(\bm{x})}\left(e^{\mathrm{i} k_{p} \boldsymbol{x} \cdot \boldsymbol{d}}\boldsymbol{d} \right) g_{p}(\boldsymbol{d})  \mathrm{d} s(\boldsymbol{d})+\int_S \sqrt{\frac{k_{s}}{\omega}} \bm{T}_{\bm{\nu}(\bm{x})}\left(e^{\mathrm{i} k_{s} \boldsymbol{x} \cdot \boldsymbol{d}} \boldsymbol{d}^{\perp} \right) g_{s}(\boldsymbol{d}) \mathrm{d} s(\boldsymbol{d})\right\},
			\end{aligned}
		\end{equation*}
		for $\bm{x}\in\partial D$. On the one hand, since $\bm{\mathcal{F}}\bm{g}$ represents the far field pattern of the scattered field corresponding to $\bm{v}^{\mathrm{i}}_{\bm{g}}$ as incident field, and $-\bm{\mathcal{H}}\bm{g}$ is the boundary data for
		exterior Neumann boundary value problem, we clearly have
		\begin{equation}\label{FGH}
			\bm{\mathcal{F}}=-\bm{\mathcal{G}}\bm{\mathcal{H}}.
		\end{equation}
		On the other hand, the scattered field is represented by the double-layer potential
		\begin{equation*}
			\bm{v}^{\mathrm{s}}(\bm{x})=\int_{\partial D}\left[
			\bm{T}_{\bm{\nu}(\bm{y})}\bm{\Gamma}(\bm{x},\bm{y})\right]^\top
			\bm{\varphi}(\bm{y}) \mathrm{d} s(\bm{y}), \quad \bm{x} \in \mathbb{R}^2\setminus \overline{D}.
		\end{equation*}
		Using the asymptotic relation (\ref{Gamma_asy}), we have
		\begin{equation*}
			\bm{T}_{\boldsymbol{\nu}(\boldsymbol{y})} \boldsymbol{\Gamma}(\boldsymbol{x}, \boldsymbol{y})=\frac{e^{\mathrm{i} k_p|\boldsymbol{x}|}}{\sqrt{|\boldsymbol{x}|}} \bm{T}_{\boldsymbol{\nu}(\boldsymbol{y})} \boldsymbol{\Gamma}_{p}^{\infty}(\hat{\boldsymbol{x}}, \boldsymbol{y})+\frac{e^{\mathrm{i} k_{s}|\boldsymbol{x}|}}{\sqrt{|\boldsymbol{x}|}} \bm{T}_{\boldsymbol{\nu}(\boldsymbol{y})} \boldsymbol{\Gamma}_{s}^{\infty}(\hat{\boldsymbol{x}}, \boldsymbol{y})+\mathrm{O}\left(|\boldsymbol{x}|^{-3 / 2}\right), \quad|\boldsymbol{x}| \rightarrow \infty,
		\end{equation*}
		which also holds uniformly in all directions $\hat{\boldsymbol{x}} \in\mathbb{S}$. Here,
		\begin{align*}
			\bm{T}_{\boldsymbol{\nu}(\boldsymbol{y})} \boldsymbol{\Gamma}_{p}^{\infty}(\hat{\boldsymbol{x}}, \boldsymbol{y})&=\frac{e^{-\frac{\mathrm{i} \pi} { 4}}}{\lambda+2 \mu} \sqrt{\frac{k_p}{8 \pi}} e^{-\mathrm{i} k_{p} \hat{\boldsymbol{x}} \cdot \bm{y}}[2 \mu(\boldsymbol{\nu}(\boldsymbol{y}) \cdot \hat{\boldsymbol{x}}) \textbf{I}+\lambda \boldsymbol{\nu}(\boldsymbol{y}) \otimes \hat{\boldsymbol{x}}] \cdot \hat{\boldsymbol{x}} \otimes \hat{\boldsymbol{x}}, \\
			\bm{T}_{\boldsymbol{\nu}(\boldsymbol{y})} \boldsymbol{\Gamma}_{s}^{\infty}(\hat{\boldsymbol{x}}, \boldsymbol{y})&=e^{-\frac{\mathrm{i} \pi }{4}} \sqrt{\frac{k_s}{8 \pi}}e^{-\mathrm{i} k_s \hat{\boldsymbol{x}} \cdot \bm{y}}[(\boldsymbol{\nu}(\boldsymbol{y}) \cdot \hat{\boldsymbol{x}}) \textbf{I}+\hat{\boldsymbol{x}} \otimes \boldsymbol{\nu}(\boldsymbol{y})] \cdot
			\hat{\boldsymbol{x}}^{\perp}\otimes \hat{\boldsymbol{x}}^{\perp}.
		\end{align*}
		Thus, the far field pair of the double-layer potential are given by
		\begin{align*}
			&v_p^\infty(\hat{\bm{x}})=\frac{e^{-\frac{\mathrm{i} \pi} { 4}}}{\lambda+2 \mu} \sqrt{\frac{k_p}{8 \pi}}\int_{\partial D}[2 \mu(\boldsymbol{\nu}(\boldsymbol{y}) \cdot \hat{\boldsymbol{x}})(\hat{\boldsymbol{x}}\cdot
			\bm{\varphi}(\bm{y}))+\lambda \boldsymbol{\nu}(\boldsymbol{y})\cdot\bm{\varphi}(\bm{y})] e^{-\mathrm{i} k_{p}\hat{\bm{x}}\cdot\bm{y}} \mathrm{d}s(\bm{y}), \\
			&v_s^\infty(\hat{\bm{x}})=\frac{e^{-\frac{\mathrm{i} \pi }{4}}}{\mu} \sqrt{\frac{k_s}{8 \pi}}\int_{\partial D}
			\mu[(\boldsymbol{\nu}(\boldsymbol{y}) \cdot \hat{\boldsymbol{x}})(\hat{\boldsymbol{x}}^\bot\cdot
			\bm{\varphi}(\bm{y}))+(\hat{\bm{x}}
			\cdot\bm{\varphi}(\bm{y}))(\boldsymbol{\nu}(\boldsymbol{y})
			\cdot\hat{\boldsymbol{x}}^\bot)] e^{-\mathrm{i} k_{s}\hat{\bm{x}}\cdot\bm{y}} \mathrm{d}s(\bm{y}).
		\end{align*}
		A straightforward calculation shows that the Hilbert space adjoint $\boldsymbol{H}^*:\left[H^{1 / 2}(\partial D)\right]^2 \rightarrow \mathcal{L}^2$ of $\boldsymbol{H}$ is given by
		\begin{equation*}
			\begin{aligned}
				\bm{\mathcal{H}}^* \bm{\varphi}(\boldsymbol{\bm{\hat{x}}})
				&=e^{\frac{\mathrm{i}\pi} {4}}\left(\int_{\partial D} \frac{1}{\lambda+2\mu}\sqrt{\frac{\omega}{k_{p}}} \bm{T}_{\bm{\nu}(\bm{y})}e^{-\mathrm{i} k_{p} \bm{\hat{x}}\cdot \boldsymbol{y}}\bm{\hat{x}} \cdot \bm{\varphi}(\boldsymbol{y}) \mathrm{d} s(\boldsymbol{y}), \int_{\partial D} \frac{1}{\mu}\sqrt{\frac{\omega}{k_{s}}} \bm{T}_{\bm{\nu}(\bm{y})}e^{-\mathrm{i} k_{\mathrm{s}} \bm{\hat{x}}\cdot\boldsymbol{y}} \bm{\hat{x}}^{\perp} \cdot \bm{\varphi}(\boldsymbol{y}) \mathrm{d} s(\boldsymbol{y})\right)\nonumber\\
				&=e^{-\frac{\mathrm{i}\pi} {4}}\left( \frac{k_p}{\lambda+2\mu}\sqrt{\frac{\omega}{k_{p}}}
				\int_{\partial D}[2\mu(\boldsymbol{\nu}(\boldsymbol{y}) \cdot \hat{\boldsymbol{x}})(\hat{\boldsymbol{x}}\cdot
				\bm{\varphi}(\bm{y}))+\lambda \boldsymbol{\nu}(\boldsymbol{y})\cdot\bm{\varphi}(\bm{y})] e^{-\mathrm{i} k_{p}\hat{\bm{x}}\cdot\bm{y}} \mathrm{d}s(\bm{y}),\right.\\
				&\quad\quad\quad\left. \frac{k_s}{\mu}\sqrt{\frac{\omega}{k_{s}}}\int_{\partial D}
				\mu[(\boldsymbol{\nu}(\boldsymbol{y}) \cdot \hat{\boldsymbol{x}})(\hat{\boldsymbol{x}}^\bot\cdot
				\bm{\varphi}(\bm{y}))+(\hat{\bm{x}}
				\cdot\bm{\varphi}(\bm{y}))(\boldsymbol{\nu}(\boldsymbol{y})
				\cdot\hat{\boldsymbol{x}}^\bot)] e^{-\mathrm{i} k_{s}\hat{\bm{x}}\cdot\bm{y}} \mathrm{d}s(\bm{y})
				\right)
			\end{aligned}
		\end{equation*}
		for $\hat{\bm{x}}\in\mathbb{S}$. We note that $\frac{1}{\sqrt{8\pi\omega}}\bm{\mathcal{H}}^*\bm{\varphi}$ is the far field pair $(v^\infty_p,v^\infty_s)$ of the scattered field $\bm{v}^{\mathrm{s}}$. Since $\bm{T}_{\bm{\nu}}\bm{v}^{\mathrm{s}}=\bm{\mathcal{N}}\bm{\varphi}$ on $\partial D$, one can obtain $\bm{\mathcal{H}}^{*}\bm{\varphi}=
		\sqrt{8\pi\omega}\bm{\mathcal{G}}\bm{\mathcal{N}}
		\bm{\varphi}$
		and consequently
		\begin{equation}\label{HNG}
			\bm{\mathcal{H}}=\sqrt{8\pi\omega}
			\bm{\mathcal{N}}^*\bm{\mathcal{G}}^*
		\end{equation}
		The statement follows by combining (\ref{FGH}) and (\ref{HNG}).
	\end{proof}
	
	\begin{lemma}\label{G_properties}
		The data-to-pattern operator $\bm{\mathcal{G}}$ is compact, injectivity with dense range in $\mathcal{L}^2(\mathbb{S})$.
	\end{lemma}
	\begin{proof}
		Let $\bm{\mathcal{G}} \boldsymbol{f}=\mathbf{0}$, i.e. $\boldsymbol{v}^\infty=\bm{0}$ on $\mathbb{S}$. It follows from Rellich's lemma that $\boldsymbol{v}^{\mathrm{s}}=
		\bm{0}$  in $\mathbb{R}^2\setminus \overline{D}$.
		We choose a circle
		$B_R:=\{x:|x|<R\}$ centered at the origin which contains $\overline{D}$ in its interior. Then, for $\bm{f} \in [H^{-1 / 2}(\partial D)]^2$, a weak solution of the exterior Neumann boundary value is the function $\bm{v}^{\mathrm{s}}\in \left[H_{\mathrm{l o c}}^1\left(\mathbb{R}^2 \backslash \overline{D}\right)\right]^2$ such that
		\begin{equation*}
			-\int_{(\mathbb{R}^2 \backslash \overline{D})\cap B_R}\mathcal{E}(\bm{w},\bm{v}^{\mathrm{s}})
			-\omega^2\bm{v}^{\mathrm{s}}\bm{w}
			~\mathrm{d}\bm{x}+\int_{\partial B_R}\bm{w}\cdot\bm{T}_{\bm{\nu}}\bm{v}^{\mathrm{s}} ~\mathrm{d} s=\int_{\partial D} \bm{w}\cdot\bm{f} ~\mathrm{d} s
		\end{equation*}
		for all $\bm{w}\in \left[H_{\text {loc }}^1\left(\mathbb{R}^2 \backslash \overline{D}\right)\right]^2$.
		From the weak formulation one can obtain
		$$
		\int_{\partial D} \bm{w}\cdot\bm{f} ~\mathrm{d}s=0
		$$
		for all $\bm{w} \in [H^{1 / 2}(\partial D)]^2$, which give $\bm{f}=\bm{0}$. Thus, $\bm{\mathcal{G}}$ is injective.
		
		In the following we prove the compactness of operator $\bm{\mathcal{G}}$. Using the integral representation (\ref{integral re}), we can decompose $\bm{\mathcal{G}}$ as $\bm{\mathcal{G}}=\bm{\mathcal{G}}_2 \bm{\mathcal{G}}_1$ where $\bm{\mathcal{G}}_1: [H^{-1 / 2}(\partial D)]^2 \rightarrow [C(\partial B_R)]^2 \times [C(\partial B_R)]^2$ and $\bm{\mathcal{G}}_2: [C(\partial B_R)]^2 \times [C(\partial B_R)]^2 \rightarrow \mathcal{L}^2(\mathbb{S})$ are defined by $\bm{\mathcal{G}}_1 \bm{f}=\left(\left.\bm{v}\right|_{\partial B_R}, \left.\bm{T}_{\bm{\nu}}\bm{v} \right|_{\partial B_R}\right)$ and
		\begin{equation*}
			\begin{aligned}
				\bm{\mathcal{G}}_2(\bm{g},\bm{h})(\hat{\bm{x}})
				=\left(\frac{1}{\lambda+2 \mu} \frac{\mathrm{e}^{\mathrm{i} \pi / 4}}{\sqrt{8 \pi k_{p}}}\int_{\partial B_R}
				\left\{\bm{T}_{\bm{\nu}(\bm{y})}e^{-ik_p\hat{\bm{x}}\cdot \bm{y}}\hat{\bm{x}}\cdot \bm{g}(\bm{y}) -\hat{\bm{x}}\cdot \bm{h}(\bm{y})e^{-ik_p\hat{\bm{x}}\cdot \bm{y}}\right\} \mathrm{d} s(\bm{y}),\right.\\ \left.
				\frac{1}{\mu} \frac{\mathrm{e}^{\mathrm{i} \pi / 4}}{\sqrt{8 \pi k_{\mathrm{s}}}}\int_{\partial B_R}
				\left\{\bm{T}_{\nu(\bm{y})}e^{-ik_s\hat{\bm{x}}\cdot \bm{y}}\hat{\bm{x}}^\bot\cdot \bm{g}(\bm{y})
				-\hat{\bm{x}}^\bot\cdot \bm{h}(\bm{y})e^{-ik_s\hat{\bm{x}}\cdot \bm{y}}\right\}\mathrm{d} s(\bm{y})\right)
			\end{aligned}
		\end{equation*}
		respectively. Then $\bm{\mathcal{G}}_1$ is bounded by interior regularity results and $\bm{\mathcal{G}}_2$ is compact which proves compactness of $\bm{\mathcal{G}}$.
		
		To demonstrate that $\bm{\mathcal{G}}$ has dense range,  we rewrite it as an integral operator. It follows from (\ref{ufar_p}) and (\ref{ufar_s}) that the far field representation for a radiating solution $\bm{v}^s$ to the Navier equation can be written in the form
		\begin{align*}
			& \left(v_p^{\infty}, v_{{\mathrm{s}}}^{\infty}\right) \nonumber\\
			& =\left(\frac{1}{\lambda+2 \mu} \frac{e^{\frac{i
						\pi}{4}}}{\sqrt{8 \pi k_p}} \int_{\partial D} \bm{T}_{\bm{\nu}(\bm{y})} \bm{u}^{\mathrm{i}}(\bm{y},-\hat{\bm{x}};p) \cdot \bm{v}^{\mathrm{s}}(\bm{y})
			-\bm{u}^{\mathrm{i}}(\bm{y},-\hat{\bm{x}}
			;p) \cdot \bm{T}_{\bm{\nu}(\bm{y})}\bm{v}^{\mathrm{s}}(\bm{y}) \mathrm{d} s(\bm{y}),\right. \nonumber\\
			& \left.\quad\quad \frac{1}{\mu} \frac{e^{\frac{i \pi}{4}}}{\sqrt{8\pi k_s}} \int_{\partial D} \bm{T}_{\bm{\nu}(\bm{y})} \bm{u}^{\mathrm{i}}(\bm{y},-\hat{\bm{x}};s) \cdot \bm{v}^s(\bm{y})-\bm{u}^{\mathrm{i}}(\bm{y},-\hat{\bm{x}}
			;s) \cdot \bm{T}_{\bm{\nu}(\bm{y})}\bm{v}^{\mathrm{s}}(\bm{y}) \mathrm{d} s(\bm{y})\right).\nonumber
		\end{align*}
		Since both $\bm{u}^{\mathrm{s}}(\bm{y},-\hat{\bm{x}};\hat{\bm{x}})$ and $\bm{u}^{\mathrm{s}}(\bm{y},-\hat{\bm{x}};\hat{\bm{x}}^\bot)$ fulfil the Kupradze radiation condition, there holds
		\begin{align*}
			\int_{\partial D} \bm{T}_{\bm{\nu}(\bm{y})} \bm{u}^{\mathrm{s}}(\bm{y},-\hat{\bm{x}};p) \cdot \bm{v}^{\mathrm{s}}(\bm{y})-\bm{u}^{\mathrm{s}}(\bm{y},
			-\hat{\bm{x}};p) \cdot \bm{T}_{\bm{\nu}(\bm{y})} \bm{v}^{\mathrm{s}}(\bm{y}) d s(\bm{y})&=0,\\
			\int_{\partial D} \bm{T}_{\bm{\nu}(\bm{y})} \bm{u}^{\mathrm{s}}(\bm{y},-\hat{\bm{x}};s) \cdot \bm{v}^{\mathrm{s}}(\bm{y})-\bm{u}^{\mathrm{s}}(\bm{y},-\hat{\bm{x}}
			;s) \cdot \bm{T}_{\bm{\nu}(\bm{y})} \bm{v}^{\mathrm{s}}(\bm{y})d s(\bm{y})& =0,
		\end{align*}
		Adding the two preceding equations and using the Neumman boundary condition $\bm{\Gamma}_{\bm{\nu}}\bm{u}=0$
		on $\partial D$, we can obatin
		\begin{equation*}
			\begin{aligned}
				\left(v_{p}^{ \infty}(\hat{\bm{x}}), v_{\mathrm{s}}^{ \infty}(\hat{\bm{x}})\right) =&\left(-\frac{1}{\lambda+2 \mu} \frac{e^{\frac{i \pi}{4}}}{\sqrt{8\pi k_p}} \int_{\partial D} \bm{u}(\bm{y},-\hat{\bm{x}};p) \cdot \bm{T}_{\bm{\nu}(\bm{y})} \bm{v}^{\mathrm{s}}(\bm{y}) \mathrm{d} s(\bm{y}),\right.  \\
				&\left. -\frac{1}{\mu} \frac{e^{\frac{i\pi}{4}}}{\sqrt{8\pi k_s}} \int_{\partial D} \bm{u}(\bm{y},-\hat{\bm{x}};s) \cdot \bm{T}_{\bm{\nu}(\bm{y})}\bm{v}^s(\bm{y}) \mathrm{d}( \bm{y})\right),\\
			\end{aligned}
		\end{equation*}
		that is
		\begin{equation*}
			\bm{\mathcal{G}}\bm{f}(\bm{d})
			=\left(-\frac{1}{\lambda+2 \mu} \frac{e^{\frac{i \pi}{4}}}{\sqrt{8\pi k_p}} \int_{\partial D} \bm{u}(\bm{y},-\bm{d};p)\cdot\bm{f}(\bm{y})\mathrm{d} s(\bm{y}),-\frac{1}{\mu} \frac{e^{\frac{i\pi}{4}}}{\sqrt{8\pi k_s}} \int_{\partial D} \bm{u}(\bm{y},-\bm{d};s) \cdot \bm{f}(\bm{y}) \mathrm{d}s(\bm{y})\right).
		\end{equation*}
		Thus, a straightforward calculation yields that the adjoint operator $\bm{\mathcal{G}}^*: \mathcal{L}^2 \rightarrow [H^{1/2}(\partial D)]^2$ is given by
		\begin{equation*}
			\begin{aligned}
				\bm{\mathcal{G}}^*\bm{g}(\bm{x})= -\frac{e^{-\frac{i \pi}{4}}}{\sqrt{8 \pi \omega}} \int_\mathbb{S}\left\{ \sqrt{\frac{k_p}{\omega}} \overline{\bm{u}(\bm{x},-\bm{d};p)} g_p(\bm{d})
				+\sqrt{\frac{k_s}{\omega}} \overline{\bm{u}(\bm{x},-\bm{d};s)} g_s(\bm{d})\right\}\mathrm{d} s(\bm{d}),\quad \bm{x}\in \partial D.
			\end{aligned}
		\end{equation*}
		Define by
		\begin{equation*}
			\begin{aligned}
				\widetilde{\bm{v}}_{\bm{g}}^{\mathrm{i}}(\bm{x})&=\mathrm{e}^{-\frac{\mathrm{i} \pi}{4}} \int_{\mathbb{S}}\left\{\sqrt{\frac{k_{p}}{\omega}} \boldsymbol{d e}^{-\mathrm{i} k_{p} \boldsymbol{d} \cdot \boldsymbol{x}} \overline{g_{p}(\boldsymbol{d})}+\sqrt{\frac{k_{s}}{\omega}} \boldsymbol{d}^{\perp} \mathrm{e}^{-\mathrm{i} k_{s} \boldsymbol{d} \cdot \boldsymbol{x}} \overline{g_{s}(\boldsymbol{d})}\right\} \mathrm{d} s(\boldsymbol{d})\\
				&=\mathrm{e}^{-\frac{\mathrm{i} \pi}{4}} \int_{\mathbb{S}}\left\{\sqrt{\frac{k_{\mathrm{p}}}{\omega}} \bm{u}^{\mathrm{i}}(\bm{x},-\bm{d};p) \overline{g_{p}(\boldsymbol{d})}+
				\sqrt{\frac{k_{s}}{\omega}} \bm{u}^{\mathrm{i}}(\bm{x},-\bm{d};s)  \overline{g_{s}(\boldsymbol{d})}\right\} \mathrm{d} s(\boldsymbol{d})
				, \quad \bm{x} \in \mathbb{R}^2,
			\end{aligned}
		\end{equation*}
		be the Herglotz wave function with kernel $\bm{g}\in\mathcal{L}^2(\mathbb{S})$. According to
		Lemma \ref{Herglotz theorem}, we can deduce that
		\begin{equation*}
			\widetilde{\bm{v}}_{\bm{g}}(\bm{x})=\mathrm{e}^{-\frac{\mathrm{i} \pi}{4}} \int_S\left\{\sqrt{\frac{k_{p}}{\omega}} \bm{u}(\bm{x},-\bm{d};p)  \overline{g_{p}(\boldsymbol{d})}+
			\sqrt{\frac{k_{s}}{\omega}} \bm{u}(\bm{x},-\bm{d};s) \overline{g_{s}(\boldsymbol{d})}\right\} \mathrm{d} s(\boldsymbol{d})
			, \quad \bm{x} \in \mathbb{R}^2,
		\end{equation*}
		is the total field that is scattered by $\widetilde{\bm{v}}_{\bm{g}}^{\mathrm{i}}$ from $D$. Hence,
		\begin{equation}\label{G adjoint}
			\overline{\bm{\mathcal{G}}^*\bm{g}}=
			-\frac{\mathrm{i}}{\sqrt{8\pi\omega}} \widetilde{\bm{v}}_{\bm{g}},\quad \text{on}~\partial D.
		\end{equation}
		Now let $\bm{g}$ satisfy $\bm{\mathcal{G}}^* \bm{g}=\bm{0}$. Then equation (\ref{G adjoint}) implies that $\widetilde{\bm{v}}_{\bm{g}}=\bm{0}$ on $\partial D$. It follows
		from $\bm{T}_{\bm{\nu}}\bm{u}=\bm{0}$ on $\partial D$ that
		\begin{equation*}
			\bm{T}_{\bm{\nu}(\bm{x})}\widetilde{\bm{v}}_{\bm{g}}(\bm{x})
			=\mathrm{e}^{-\frac{\mathrm{i} \pi}{4}} \int_{\mathbb{S}}\left\{\sqrt{\frac{k_{\mathrm{p}}}{\omega}} \bm{T}_{\bm{\nu}(\bm{x})}\bm{u}(\bm{x},-\bm{d};p) \overline{g_{p}(\boldsymbol{d})}+
			\sqrt{\frac{k_{\mathrm{s}}}{\omega}} \bm{T}_{\bm{\nu}(\bm{x})}\bm{u}(\bm{x},-\bm{d};s) \overline{g_{\mathrm{s}}(\boldsymbol{d})}\right\} \mathrm{d} s(\boldsymbol{d})=\bm{0},
		\end{equation*}
		for $\bm{x} \in \partial D$. Therefore, by Holmgren's Theorem, it can be deduced that $\widetilde{\bm{v}}_{\bm{g}}=\bm{0}$ in $\mathbb{R}^2 \backslash \overline{D}$. Then, we can get $\widetilde{\bm{v}}^{\mathrm{i}}_{\bm{g}}
		=-\widetilde{\bm{v}}^{\mathrm{s}}_{\bm{g}}$ in   $\mathbb{R}^2 \backslash \overline{D}$.
		The entire solution $\widetilde{\bm{v}}_{\bm{g}}^{\mathrm{i}}$ satisfies the radiation condition , and consequently it must vanish identically. From Lemma \ref{g theorem}, we have $\bm{g}=\bm{0}$, which implies the injectivity of
		$\bm{\mathcal{G}}^*$. Hence $\bm{\mathcal{G}}$ has dense range by Theorem 4.6 in \cite{CK13}.
	\end{proof}

	\subsection{Boundary integral operators}
	We first briefly review  elastic single- and double-
	layer potentials and jumps relations. Given an integral
	function $\bm{\varphi}$, the integrals
	\begin{equation*}
		\widetilde{\bm{u}}(\bm{x}):= \int_{\partial D}
		\bm{\Gamma}(\bm{x},\bm{y})
		\bm{\varphi}(\bm{y}) \mathrm{d} s(\bm{y}), \quad \bm{x} \in \mathbb{R}^2\setminus\partial D,
	\end{equation*}
	and
	\begin{equation*}
		\widetilde{\bm{v}}(\bm{x}):= \int_{\partial D}
		[\bm{T}_{\bm{\nu}(\bm{x})}\bm{\Gamma}(\bm{x},\bm{y})]^\top
		\bm{\varphi}(\bm{y}) \mathrm{d} s(\bm{y}),  \quad \bm{x} \in \mathbb{R}^2\setminus\partial D,
	\end{equation*}
	are called elastic single-layer potential and elastic double-layer potentials with density $\bm{\varphi}$, respectively. As we all know, the behaviour of the
	surface potentials at the boundary $\partial D$ is
	describe by the following jump relations \cite{K79}. For the single-layer potential $\widetilde{\bm{u}}$ with density $\bm{\varphi}$, we have the jump relation
	\begin{equation*}
		\bm{T}_{\bm{\nu}}\widetilde{\bm{u}}_{+}(\bm{x})-
		\bm{T}_{\bm{\nu}}\widetilde{\bm{u}}_{-}(\bm{x})
		=-\bm{\varphi}(\bm{x}),\quad \bm{x}\in \partial D,
	\end{equation*}
	where
	\begin{equation*}
		\bm{T}_{\bm{\nu}}\widetilde{\bm{u}}_{\pm}(\bm{x}):=
		\lim_{h\rightarrow+0}\mu\left(\nabla\widetilde{\bm{u}}(\bm{x}\pm h\bm{\nu}(\bm{x}))+\nabla^\top\widetilde{\bm{u}}(\bm{x}\pm h\bm{\nu}(\bm{x}))\right)\bm{\nu}
		(\bm{x})+\lambda\nabla\cdot\widetilde{\bm{u}}(\bm{x}\pm h\bm{\nu}(\bm{x}))\bm{\nu}(\bm{x}).
	\end{equation*}
	For the double-layer potential $\widetilde{\bm{v}}$ with density $\bm{\varphi}$, we have the jump relation
	\begin{equation}\label{Double_jump}
		\widetilde{\bm{v}}_{+}(\bm{x})-
		\widetilde{\bm{v}}_{-}(\bm{x})=\bm{\varphi}(\bm{x}),\quad \bm{x}\in \partial D,
	\end{equation}
	where
	\begin{equation*}
		\widetilde{\bm{v}}_{\pm}(\bm{x}):=
		\lim_{h\rightarrow+0}\widetilde{\bm{v}}(\bm{x}\pm h\bm{\nu}(\bm{x})).
	\end{equation*}
	
	For the direct values of the single- and double-layer potentials on the boundary $\partial D$, we have more regularity. Thus, we introduce boundary integral operators defined by
	\begin{align}
		\bm(\mathcal{S}\bm{\varphi})(\bm{x})&:= \int_{\partial D}
		\bm{\Gamma}(\bm{x},\bm{y})
		\bm{\varphi}(\bm{y}) \mathrm{d} s(\bm{y}), \quad \bm{x} \in \partial D, \label{S_BIO}\\
		(\boldsymbol{\mathcal{K}}'\bm{\varphi})(\bm{x})&:= \int_{\partial D}
		[\bm{T}_{\bm{\nu}(\bm{x})}\bm{\Gamma}(\bm{x},\bm{y})]
		\bm{\varphi}(\bm{y}) \mathrm{d} s(\bm{y}), \quad \bm{x} \in \partial D, \label{K'_BIO}\\
		(\boldsymbol{\mathcal{K}}\bm{\varphi})(\bm{x})&:= \int_{\partial D}
		[\bm{T}_{\bm{\nu}(\bm{y})}\bm{\Gamma}(\bm{x},\bm{y})]^{\top}
		\bm{\varphi}(\bm{y}) \mathrm{d} s(\bm{y}), \quad \bm{x} \in \partial D, \label{K_BIO}
	\end{align}
	The operator $\bm{\mathcal{S}}$ is compact from $[H^{-1/2}(\partial D)]^2$ to $[H^{1/2}(\partial D)]^2$ since its kernel is weakly singular; see Appendix Proposition \ref{S_kernel}. However, in contrast to the acoustic case, the operators $\bm{\mathcal{K}}'$ and $\bm{\mathcal{K}}$ are not compact from $[H^{-1/2}(\partial D)]^2$ to $[H^{1/2}(\partial D)]^2$ because their kernels are no longer
	weakly singular; see Appendix Proposition \ref{K'_kernel}
	and \ref{K_kernel}. Let
	$$
	\bm{\Gamma_0}(\boldsymbol{x}, \boldsymbol{y})=\frac{\lambda+3 \mu}{4 \pi \mu(\lambda+2 \mu)}\left(\ln\frac{1}{|\bm{x}-\bm{y}|}\mathbf{I}
	+\frac{\lambda+\mu}{\lambda+3 \mu}
	\frac{(\boldsymbol{\bm{x}-\bm{y}}) (\boldsymbol{\bm{x}-\bm{y}})^{\top}}{
		|\boldsymbol{\bm{x}-\bm{y}}|^2}\right)
	$$
	be the fundamental solution of Navier equation (\ref{navier equation}) with $\omega=0$, and define boundary integral operators $\boldsymbol{\mathcal{S}}_0$, $\boldsymbol{\mathcal{K}}'_0$, $\boldsymbol{\mathcal{K}}_0$ and $\boldsymbol{\mathcal{N}}_0$ in the same manner by
	replacing $\bm{\Gamma}(\bm{x},\bm{y})$ by $\bm{\Gamma}_0(\bm{x},\bm{y})$ in (\ref{S_BIO})-(\ref{K_BIO}) and (\ref{N oper}), respectively.
	It can be shown that the operators
	$\boldsymbol{\mathcal{K}}-\boldsymbol{\mathcal{K}}_0, \boldsymbol{\mathcal{K}}'-\boldsymbol{\mathcal{K}}_0'$ and $\boldsymbol{\mathcal{N}}-\boldsymbol{\mathcal{N}}_0$ have weakly singular kernels; see \textbf{Appendix} Proposition \ref{K_kernel}, \ref{K_kernel} and \ref{N_Ker}.  Thus, they are compact.

	\begin{lemma}\label{N properties}
		Assume that $\omega^2$ is not a Neumann eigenvalue of $-\Delta^*$ in $D$. Then there hold that
		\begin{itemize}
			\item[(a)] $\bm{\mathcal{N}}$ is an isomorphism from $[H^{1 / 2}(\partial D)]^2$ onto $[H^{-1 / 2}(\partial D)]^2$.
			\item[(b)] $\operatorname{Im}\left( \bm{\mathcal{N}} \bm{\varphi}, \bm{\varphi}\right)>0$ for all $ \bm{\varphi} \in [H^{1 / 2}(\partial D)]^2$ with $\bm{\varphi} \neq 0$. Here, $(\cdot, \cdot)$ denotes the duality pair in $\left( [H^{-1 / 2}(\partial D)]^2, [H^{1 / 2}(\partial D)]^2\right)$;
			\item[(c)] Let $\bm{\mathcal{N}}_{\mathrm{i}}$ be the boundary operator (\ref{N oper}) corresponding to the frequency $\omega=\mathrm{i}$. The operator $\bm{\mathcal{N}}_{\mathrm{i}}$ is self-adjoint in $[L^2(\partial D)]^2$ with respect to the bilinear from $\langle\bm{\varphi},\bm{\psi}\rangle_{[L^2(\partial D)]^2}:=\int_{\partial D}\bm{\varphi}\cdot\bm{\psi}~\mathrm{d}s$.
			Moreover, $-\bm{\mathcal{N}}_{\mathrm{i}}$ is coercive as an operator from $[H^{1 / 2}(\partial D)]^2$ onto $[H^{-1 / 2}(\partial D)]^2$, i.e.,
			$$
			-\left(\bm{\mathcal{N}}_{\mathrm{i}} \bm{\varphi}, \bm{\varphi}\right) \geq c_0\|\bm{\varphi}\|_{[H^{1 / 2}(\partial D)]^2} \quad \text { for all } \bm{\varphi} \in [H^{1 / 2}(\partial D)]^2 .
			$$
			\item[(d)] The difference $\bm{\mathcal{N}}-\bm{\mathcal{N}}_{\mathrm{i}}$ is compact from $[H^{1 / 2}(\partial D)]^2$ into $[H^{-1 / 2}(\partial D)]^2$.
		\end{itemize}
	\end{lemma}
	
	\begin{proof}
		\begin{itemize}
			\item[(a)]
			According to Theorem 7.17 in \cite{W00}, we know that $\bm{\mathcal{N}}$ is a Fredholm operator of index zero from
			$\left[H^{1 / 2}(\partial D)\right]^2$ to $\left[H^{-1 / 2}(\partial D)\right]^2$. Then we only need to prove injectivity of $\bm{\mathcal{N}}$. Let $\bm{\varphi}
			\in [H^{1 / 2}(\partial D)]^2$ be a solution to the homogeneous equation $\bm{\mathcal{N}}\bm{\varphi}=0$. Then the double-layer potential $\widetilde{\bm{v}}$ with denstiy $\bm{\varphi}$ belongs to $[H^{1}(D)]^2$ and $[H^{1}_{\mathrm{loc}}(\mathbb{R}^2\setminus \overline{D})]^2$, and satisfies $\bm{T}\widetilde{\bm{v}}=0$ on $\partial D$. The uniqueness for the exterior Neumann problem yields
			$\widetilde{\bm{v}}=0$ in $\mathbb{R}^2\setminus\overline{D}$, and the assumption  based on $\omega^2$ yields $\widetilde{\bm{v}}=0$ in $D$. Now it follows from the jump relations (\ref{Double_jump}) that $\bm{\varphi}=0$. The application of the Fredholm alternative (\cite{W00} Theorem 2.27) to the operator $\bm{\mathcal{N}}$ yields
			that $\bm{\mathcal{N}}$ is an isomorphism from  $\left[H^{1 / 2}(\partial D)\right]^2$ to $\left[H^{-1 / 2}(\partial D)\right]^2$.

			\item[(b)]Define the double potential $\widetilde{\bm{v}}$ in the same way as above. It
			follows from jump relations that $\bm{\varphi}=\widetilde{\bm{v}}_{+}-\widetilde{\bm{v}}_{-}$ and $\bm{T}_{\bm{\nu}}\widetilde{\bm{v}_{-}}=$ $\bm{T}_{\bm{\nu}}\widetilde{\bm{v}}_{+}=
			\bm{\mathcal{N}} \bm{\varphi}$ on $\partial D$. Therefore, applying the Betti's formula in $D$ and in $D_R:=\left\{\bm{x} \in \mathbb{R}^2 \backslash \overline{D}:|\bm{x}|<R\right\}$ obtains
			\begin{equation*}
				\begin{aligned}
					(\bm{\mathcal{N}}\bm{\varphi}, \boldsymbol{\varphi}) & =\left(\boldsymbol{T}_{\bm{\nu}} \widetilde{\boldsymbol{v}},  \widetilde{\boldsymbol{v}}_{+}
					-\widetilde{\boldsymbol{v}}_{-}\right) \\
					&=\int_{\partial D}\boldsymbol{T}_{\bm{\nu}}
					\widetilde{\boldsymbol{v}}\cdot\overline{ \widetilde{\boldsymbol{v}}_{+}}\mathrm{d}  s-\int_{\partial D}\boldsymbol{T}_{\bm{\nu}}\widetilde{\boldsymbol{v}}\cdot\overline{ \widetilde{\boldsymbol{v}}_{-}}\mathrm{d}  s,\\
					&=\int_{\partial D_R} \boldsymbol{T}_{\bm{\nu}} \widetilde{\boldsymbol{v}} \cdot \overline{\widetilde{\boldsymbol{v}}} \mathrm{d} s-\int_{D\cup D_{R}}\left\{\mathcal{E}
					(\widetilde{\boldsymbol{v}}, \overline{\widetilde{\boldsymbol{v}}})
					+\overline{\widetilde{\boldsymbol{v}}}\cdot \Delta^*\widetilde{\boldsymbol{v}}\right\}\mathrm{d} \bm{x}\\
					& =\int_{\partial D_R} \boldsymbol{T}_{\bm{\nu}} \widetilde{\boldsymbol{v}} \cdot \overline{\widetilde{{\boldsymbol{v}}}} \mathrm{d} s-\int_{D \cup D_R}\left\{\mathcal{E}(\widetilde{\boldsymbol{v}}, \overline{\widetilde{\boldsymbol{v}}})-
					\omega^2|\widetilde{\boldsymbol{v}}|^2\right\} \mathrm{d} \bm{x}
				\end{aligned}
			\end{equation*}
			where
			\begin{equation*}
				\begin{aligned} \mathcal{E}(\widetilde{\boldsymbol{v}},
					\overline{\widetilde{\boldsymbol{v}}})
					=2\mu\left(\left|\operatorname{grad} \widetilde{v}_1\right|^2+\left|\operatorname{grad} \widetilde{v}_2\right|^2\right)
					+\lambda|\operatorname{div} \widetilde{\boldsymbol{v}}|^2-\mu\left|\operatorname{div}^{\perp} \widetilde{\boldsymbol{v}}\right|^2.
				\end{aligned}
			\end{equation*}
			Since $\widetilde{\boldsymbol{v}}$ satisfies Kupradze's radiation condition, we let $R$ approach infinity and use asymptotic behaviour (\ref{scattered_asy}) and (\ref{Tscattered_asy}) to derive
			\begin{equation}\label{N_psi}
				(\bm{\mathcal{N}} \boldsymbol{\varphi}, \boldsymbol{\varphi})=\mathrm{i} \omega\left\langle\widetilde{\boldsymbol{v}}^{\infty}, \widetilde{\boldsymbol{v}}^{\infty}\right\rangle-
				\int_{\mathbb{R}^2}\left\{\mathcal{E}
				(\widetilde{\boldsymbol{v}}, \overline{\widetilde{\boldsymbol{v}}})
				-\omega^2|\widetilde{\boldsymbol{v}}|^2\right\} \mathrm{d} \bm{x}.
			\end{equation}
			Taking the imaginary part of both sides of (\ref{N_psi}) yields
			\begin{equation*}
				\mathrm{Im}(\bm{\mathcal{N}} \boldsymbol{\varphi}, \boldsymbol{\varphi})=
				\omega\left\langle\widetilde{\boldsymbol{v}}^{\infty}, \widetilde{\boldsymbol{v}}^{\infty}\right\rangle=
				\int_{\mathbb{S}}\frac{\omega^2}{k_p}
				\left|\widetilde{v}_p^{\infty}\right|^2+\frac{\omega^2}{k_s}
				\left|\widetilde{v}_s^{\infty}\right|^2\mathrm{d}s\geq0.
			\end{equation*}
			Thus, if $\mathrm{Im}(\bm{\mathcal{N}} \boldsymbol{\varphi}, \boldsymbol{\varphi})=\bm{0}$ for some $\boldsymbol{\varphi}\in [H^{1/2}(\partial D)]^2$, we have $\widetilde{\boldsymbol{v}}^{\infty}=\bm{0}$. By Rellich's lemma, we conclude that $\widetilde{\boldsymbol{v}}$ vanishes outside of $D$.
			Therefore, $\boldsymbol{\mathcal{N}} \boldsymbol{\varphi}=\bm{0}$ on $\partial D$ by the trace theorem. Since $\bm{\mathcal{N}}$ is an isomorphism, $\bm{\varphi}$ must vanish.
			\item[(c)]Let $\widetilde{\bm{v}}_1$ and $\widetilde{\bm{v}}_2$ denote the double-layer potentials with the densities $
			\bm{\varphi}$ and $\bm{\psi}$, respectively. By the radiation condition (\ref{scattered_asy}), the jump relation (\ref{Double_jump}), and Betti's formula, we find
			\begin{equation}\notag
				\begin{aligned}
					\int_{\partial D} \bm{\mathcal{N}}_{\mathrm{i}} \bm{\varphi}\cdot\bm{\psi}\mathrm{~d} s&=\int_{\partial D} \bm{T}_{\bm{\nu}}\widetilde{\bm{v}}_1\cdot
					\left(\widetilde{\bm{v}}_{2,+}
					-\widetilde{\bm{v}}_{2,-}\right) \mathrm{d} s\\
					&=\int_{\partial D}\left(\widetilde{\bm{v}}_{1,+}-
					\widetilde{\bm{v}}_{1,-}\right)\cdot \bm{T}_{\bm{\nu}}\widetilde{\bm{v}}_2
					\mathrm{d}s=\int_{\partial D} \bm{\varphi} \cdot\bm{\mathcal{N}}_{\mathrm{i}} \bm{\psi} \mathrm{~d} s,
				\end{aligned}
			\end{equation}
			Thus, $\bm{\mathcal{N}}_{\mathrm{i}}$ is self-adjoint.
			For $\omega=\mathrm{i}$, the same arguments as above in Lemma \ref{N properties} (b) and Lemma 2.1 in \cite{A01} yield
			\begin{equation}
				\begin{aligned}
					-(\bm{\mathcal{N}}_{\mathrm{i}}\bm{\varphi}, \boldsymbol{\varphi})&=\int_{D \cup D_R}\left\{\mathcal{E}(\overline{\widetilde{\boldsymbol{v}}}, \widetilde{\boldsymbol{v}})
					+|\widetilde{\boldsymbol{v}}|^2\right\} \mathrm{d} \boldsymbol{x}-\int_{\partial D_R} \boldsymbol{T}_{\boldsymbol{\nu}} \widetilde{\boldsymbol{v}} \cdot \overline{\widetilde{\boldsymbol{v}}} \mathrm{d}s\\
					&=-\int_{D \cup D_R} \overline{\widetilde{\boldsymbol{v}}}\cdot
					\Delta^*\widetilde{\boldsymbol{v}} \mathrm{d} \boldsymbol{x}+\int_{D \cup D_R} |\widetilde{\boldsymbol{v}}|^2 \mathrm{d} \boldsymbol{x}\\
					&=\int_{D \cup D_R} \mathcal{E}_{\widetilde{\mu}}(\overline{\widetilde{\boldsymbol{v}}}
					,\widetilde{\boldsymbol{v}})\mathrm{d} \boldsymbol{x}
					-\int_{\partial D_R} \boldsymbol{P}_{\boldsymbol{\nu}} \widetilde{\boldsymbol{v}} \cdot \overline{\widetilde{\boldsymbol{v}}} \mathrm{d}s+\int_{D \cup D_R} |\widetilde{\boldsymbol{v}}|^2 \mathrm{d} \boldsymbol{x}.
				\end{aligned}
			\end{equation}
			Using the fact that the field $\widetilde{\bm{v}}$ decays exponentially for $\omega=\mathrm{i}$, we deduce
			\begin{equation}
				-(\bm{\mathcal{N}}_{\mathrm{i}}\bm{\varphi}, \boldsymbol{\varphi})=\int_{\mathbb{R}^2} \left\{\mathcal{E}_{\widetilde{\mu}}(\overline{\widetilde{\boldsymbol{v}}}
				,\widetilde{\boldsymbol{v}})
				+|\widetilde{\boldsymbol{v}}|^2\right\}
				\mathrm{d} \boldsymbol{x},
			\end{equation}
			where
			\begin{equation*}
				\begin{aligned} \mathcal{E}_{\widetilde{\mu}}(\overline{
						\widetilde{\boldsymbol{v}}},
					\widetilde{\boldsymbol{v}})
					=(\mu+\widetilde{\mu})\left(\left|\operatorname{grad} \widetilde{v}_1\right|^2+\left|\operatorname{grad} \widetilde{v}_2\right|^2\right)
					+\left(\mu+\lambda-\widetilde{\mu}\right)
					|\operatorname{div} \widetilde{\boldsymbol{v}}|^2
					-\widetilde{\mu}\left|\operatorname{div}^{\perp} \widetilde{\boldsymbol{v}}\right|^2.
				\end{aligned}
			\end{equation*}
			we let $\widetilde{\mu}=0$ and thus obtain
			\begin{equation}
				\int_{\mathbb{R}^2}\left\{\mathcal{E}_0(\overline{\widetilde{\boldsymbol{v}}}, \widetilde{\boldsymbol{v}})+|\widetilde{\boldsymbol{v}}|^2\right\} \mathrm{d} \boldsymbol{x}\geq C\|\boldsymbol{\widetilde{v}}\|_{\left[H^1(\mathbb{R}^2)
					\right]^2}.
			\end{equation}
			The trace theorem and the boundedness of $\bm{\mathcal{N}}_\mathrm{i}^{-1}$ yields
			\begin{equation*}
				-(\bm{\mathcal{N}}_{\mathrm{i}}\bm{\varphi}, \boldsymbol{\varphi})\geq C\|\boldsymbol{\bm{T}_{\bm{\nu}}\widetilde{\bm{v}}}\|
				_{\left[H^{-1/2}\left(\partial D\right)\right]^2}=C\left\|\bm{\mathcal{N}}_\mathrm{i} \boldsymbol{\varphi}\right\|_{\left[H^{-1 / 2}(\partial D)\right]^2} \geqslant C\|\boldsymbol{\varphi}\|_{\left[H^{1 / 2}(\partial D)\right]^2},
			\end{equation*}
			where $C_1>0$ and $C>0$ are constants which are independent of $\bm{\varphi}\in\left[H^{1 / 2}(\partial D)\right]^2$.
			
			\item[(d)] The kernel of the operator $\bm{\mathcal{N}}-\bm{\mathcal{N}}_{\mathrm{i}}$ can be rewritten as
			\begin{equation*}
				\boldsymbol{T}_{\nu(\boldsymbol{x})}\left\{\boldsymbol{T}_{\nu(\boldsymbol{y})}\left[\boldsymbol{\Gamma}(\boldsymbol{x}, \boldsymbol{y})-\boldsymbol{\Gamma}_{\mathrm{i}}(\boldsymbol{x}, \boldsymbol{y})\right]\right\}^{\top}=\sum_{j=1}^2 \sum_{k=0}^2\left\{\gamma_j^{(k)}
				(|\boldsymbol{x}-\boldsymbol{y}|)-\gamma_{\mathrm{i}, j}^{(k)}(|\boldsymbol{x}-\boldsymbol{y}|)\right\} \mathbf{N}_j^{(k)}(\boldsymbol{x}, \boldsymbol{y}),
			\end{equation*}
			where $\bm{\Gamma}_{\mathrm{i}}(\bm{x},\bm{y})$ is the fundamental solution of Navier equation (\ref{navier equation}) with $\omega=\mathrm{i}$, and $\gamma_{\mathrm{i}, j}^{(k)}$ is defined by letting $\omega=\mathrm{i}$ in $\gamma_{j}^{(k)}$, $j=1,2$ and $k=0,1,2$. See \textbf{Appendix} for the specific expressions of $\gamma_{j}^{(k)}(|\boldsymbol{x}-\boldsymbol{y}|)$ and $\mathbf{N}_j^{(k)}(\boldsymbol{x}, \boldsymbol{y})$. By Proposition \ref{N_Ker}, the strong singularity of the kernel of the operator $\bm{\mathcal{N}}$ is reflected in (\ref{chi3}), (\ref{chi1}) and (\ref{chi12}). The singular terms of (\ref{chi3}), (\ref{chi1}) and (\ref{chi12}) can be simplified as
			\begin{equation}\label{singular_term}
				\frac{k_p^2+k_s^2}{4 \pi \omega^2v^2}=\frac{\lambda+3\mu}{4\pi\mu(\lambda+2\mu)v^2},
				\quad \frac{k_s^2-k_p^2}{4 \pi \omega^2v^2}=
				\frac{\lambda+\mu}{4\pi\mu(\lambda+2\mu)v^2}.
			\end{equation}
			One finds that the strong singular terms (\ref{singular_term}) is independent of frequency $\omega$. Thus, the difference  $\bm{\mathcal{N}}-\bm{\mathcal{N}}_{\mathrm{i}}$ eliminates these strong singular terms. The operator $\bm{\mathcal{N}}-\bm{\mathcal{N}}_{\mathrm{i}}$ is compact from $\left[H^{1 / 2}(\partial D)\right]^2$ to $\left[H^{-1 / 2}(\partial D)\right]^2$ since its kernel is weakly singular.
		\end{itemize}
	\end{proof}
	
	\begin{remark}
		Here, we give an alternative proof to prove Lemma \ref{N properties} (a). The operators
		$\boldsymbol{\mathcal{S}}, \boldsymbol{\mathcal{K}}, \boldsymbol{\mathcal{K}}^{\prime}$ are bounded from $[H^{-1 / 2}(\partial D)]^2$ to $[H^{1 / 2}(\partial D)]^2$ and $\boldsymbol{\mathcal{N}}$ is bounded from $[H^{1 / 2}(\partial D)]^2$ to $[H^{-1 / 2}(\partial D)]^2$. They satisfy the following commutation relation (see \cite{W00}):
		\begin{equation*}
			\boldsymbol{\mathcal{N}} \boldsymbol{\mathcal{S}}=\left(\boldsymbol{\mathcal{K}}
			^{\prime}-\frac{1}{2} \boldsymbol{\mathcal{I}}\right)\left(\boldsymbol{\mathcal{K}}^{\prime
			}+\frac{1}{2} \boldsymbol{\mathcal{I}}\right),
		\end{equation*}
		If $\omega^2$ is not a Dirichlet eigenvalue of $-\Delta^*$ in $D$, the operator $\bm{\mathcal{S}}$
		are invertible from $[H^{-1/2}(\partial D)]^2$ to $[H^{1/2}(\partial D)]^2$. If $\omega^2$ is not a Neumann eigenvalue of $-\Delta^*$ in $D$, the operator $\pm\frac{1}{2} \boldsymbol{\mathcal{I}}+\boldsymbol{\mathcal{K}}^{\prime}$
		are invertible in $[H^{-1/2}(\partial D)]^2$; see also \cite{AAL99,ABGKLA15}.
		We conclude that $\bm{\mathcal{N}}$ is an isomorphism from $[H^{1 / 2}(\partial D)]^2$ to $[H^{-1 / 2}(\partial D)]^2$. We can write
		\begin{equation*}
			\bm{\mathcal{N}}^{-1}=\boldsymbol{\mathcal{S}}\left(\boldsymbol{\mathcal{K}}^{\prime}
			+\frac{1}{2} \boldsymbol{\mathcal{I}}\right)^{-1}\left(\boldsymbol{\mathcal{K}}^{\prime}
			-\frac{1}{2} \boldsymbol{\mathcal{I}}\right)^{-1}.
		\end{equation*}
	\end{remark}
	\begin{remark}
		It follows from the process of proof for Lemma \ref{N properties}(d) that the kernel of $\bm{\mathcal{N}}_w-\bm{\mathcal{N}}_{w_0}$ is weakly singular for any two different frequencies $w$ and $w_0$. Thus, $\bm{\mathcal{N}}_w-\bm{\mathcal{N}}_{w_0}$ is compact from $[H^{1 / 2}(\partial D)]^2$ into $[H^{-1 / 2}(\partial D)]^2$.
	\end{remark}

	\section{The factorization method}\label{Section_4}
	The basic idea of the linear sampling method is to find a
	function
	$\bm{g}(\cdot, \bm{z};\bm{p})\in\mathcal{L}^{2}(\mathbb{S})$, $\bm{z}\in D$, such that
	\begin{equation}\label{Fg}
		\bm{\mathcal{F}} \bm{g}(\cdot, \bm{z};\bm{p})=\bm{\Gamma}^{\infty}(\cdot, \bm{z};\bm{p}),
	\end{equation}
	where $\bm{\Gamma}^{\infty}(\cdot, \bm{z};\bm{p})=\bm{\Gamma}^{\infty}(\cdot, \bm{z})\bm{p}$, and $\bm{p}\in\mathbb{S}$ is the polarization direction of the point source. However, the equation (\ref{Fg}) may not be solvable in general and it is unclear what will happen when $\bm{z}\notin D$.
	Such drawbacks can been overcome by
	replacing (\ref{Fg}) by
	\begin{equation}\label{FFg_equ}
		\left(\bm{\mathcal{F}}^*\bm{\mathcal{F}}\right)^{1 / 4} \bm{g}(\cdot, \bm{z};\bm{p})=\bm{\Gamma}^{\infty}(\cdot, \bm{z};\bm{p}),
	\end{equation}
	which is known as the factorization method; see \cite{A01,AK02}.
	
	\subsection{Reconstruction for FF case}
	We first consider the FF case.
	\begin{lemma}\label{z_range}
		For any $\bm{z} \in \mathbb{R}^2$ and $\bm{p}\in \mathbb{S}$, $\bm{z} \in D$ if and only if $\bm{\Gamma}^{\infty}(\cdot,\bm{z};\bm{p})$ belongs to the range of $\bm{\mathcal{G}}$.
	\end{lemma}
	\begin{proof}
		Assume that $\boldsymbol{z} \in D$ and denote $\boldsymbol{v}^{\mathrm{s}}(\boldsymbol{x}):=
		\bm{\Gamma}(\boldsymbol{x}, \boldsymbol{z}) \boldsymbol{p}$. Then $\boldsymbol{v}^s$ solves the exterior Neumann boundary value problem with boundary
		data $\bm{f}=\left.\bm{T}_{
			\bm{\nu}}\boldsymbol{v}^{\mathrm{s}}\right|_{\partial D} \in\left[H^{-1 / 2}(\partial D)\right]^2$. We know that
		$\bm{\Gamma}^{\infty}(\cdot, \boldsymbol{z} ; \boldsymbol{p})$ is the the far-field pattern of $\boldsymbol{v}^{\mathrm{s}}$. Thus, $\bm{\mathcal{G}}\bm{f}=\bm{\Gamma}^{\infty}(\cdot, \boldsymbol{z} ; \boldsymbol{p})$, i.e.,
		$\bm{\Gamma}^{\infty}(\cdot, \boldsymbol{z} ; \boldsymbol{p}) \in \mathcal{R}(\boldsymbol{G})$.
		
		For $\bm{z}\notin D$, assume that there exists $\boldsymbol{f} \in\left[H^{-1 / 2}(\partial D)\right]^2$ such that $\boldsymbol{\mathcal{G}} \boldsymbol{f}=\bm{\Gamma}^{\infty}(\cdot, \boldsymbol{z} ; \boldsymbol{p})$. Let $\boldsymbol{v}^s$ denote the solution of the exterior Neumann problem with boundary data $\boldsymbol{f}$. Then $\boldsymbol{v}^{\infty}={\boldsymbol{\mathcal{G}}} \boldsymbol{f}=\bm{\Gamma}^{\infty}(\cdot, \boldsymbol{z} ; \boldsymbol{p})$. Since $\bm{\Gamma}^{\infty}(\cdot, \boldsymbol{z} ; \boldsymbol{p})$ is the far-field pattern of $\bm{\Gamma}(\cdot, \boldsymbol{z}) \boldsymbol{p}$, it follows from Rellich's lemma that $\boldsymbol{v}^{\mathrm{s}}(\boldsymbol{x})=\bm{\Gamma}(\bm{x}, \boldsymbol{z})\bm{p}$ in $\mathbb{R}^2 \backslash(\overline{D} \cup\{\bm{z}\})$. However, the solution of the exterior Neumann problem is analytic in $\mathbb{R}^2 \backslash \overline{D}$ while $\bm{\Gamma}(\cdot, \boldsymbol{z}) \boldsymbol{p}$ has a singularity at $\boldsymbol{z}$. This leads to the contradiction for $\bm{z}\in\mathbb{R}^2 \setminus\overline {D}$.
		
		Assume $\boldsymbol{z} \in \partial D$. Then, from the boundary condition, $\boldsymbol{f}(\boldsymbol{x})=\bm{T}_{
			\bm{\nu}}\left(\bm{\Gamma}(\bm{x}, \boldsymbol{z}) \boldsymbol{p}\right), \boldsymbol{x} \in \partial D$, i.e. $\left.(\bm{T}_{
			\bm{\nu}}\left(\bm{\Gamma}(\bm{x}, \boldsymbol{z}) \boldsymbol{p}\right))\right|_{\partial D} \in\left[H^{-1 / 2}(\partial D)\right]^2$. However, it follows from Appendix Proposition \ref{K_kernel} that
		$\bm{T}_{
			\bm{\nu}}\left(\bm{\Gamma}(\bm{x}, \boldsymbol{z})\right)=\mathrm{O}
		(1/|\boldsymbol{x}-\boldsymbol{z}|^2)$ for $|\boldsymbol{x}-\boldsymbol{z}|\rightarrow 0$. Thus, this leads a contradiction to
		$\bm{T}_{
			\bm{\nu}}\left(\bm{\Gamma}(\bm{x}, \boldsymbol{z})\right)\notin [L_{\mathrm{loc}}^2(\mathbb{R}^2\setminus\overline{D})]^2$ for $\bm{z}\in \partial D$.
	\end{proof}

	\begin{lemma}\label{G_range}
		Assume that $\omega^2$ is not a Neumann eigenvalue of $-\Delta^*$ in $D$. The ranges of  $\bm{\mathcal{G}}$ and $\left(\bm{\mathcal{F}}^*\bm{\mathcal{F}}\right)^{1 / 4} $ coincide. Furthermore,
		the operators  $\left(\bm{\mathcal{F}}^*\bm{\mathcal{F}}\right)^{-1 / 4} \bm{\mathcal{G}}$ and $\bm{\mathcal{G}}^{-1}\left(\bm{\mathcal{F}}^*
		\bm{\mathcal{F}}\right)^{1 / 4} $ are isomorphism from $[H^{1 / 2}(\partial D)]^2$ onto $\mathcal{L}^2(\mathbb{S})$ and
		from $\mathcal{L}^2(\mathbb{S})$ onto $[H^{1 / 2}(\partial D)]^2$,
		respectively.
	\end{lemma}
	\begin{proof}
		By taking $H=\mathcal{L}^2\left(\mathbb{S}\right)$, $X=[H^{-1 / 2}(\partial D)]^2$, $B=\bm{\mathcal{G}}$ and $A=-\sqrt{8 \pi \omega}\bm{\mathcal{N}}^*$ in Theorem 1.23 in \cite{KAG07}, and using Lemma \ref{N properties}, we immediately complete the proof.
	\end{proof}
	
	Combining Lemmas \ref{G_properties}, \ref{z_range} and \ref{G_range} yields the
	following final characterization of the scatterer $D$.

	\begin{theorem}\label{FF_Th}
		Let $\bm{\mathcal{F}}$ be the far field operator
		and assume that $\omega^2$ is not a Neumann eigenvalue of $-\Delta^*$ in $D$. Then
		a point $\bm{z} \in \mathbb{R}^2$ belongs to $D$ if and only if the series
		$$
		\sum_{j=1}^{\infty} \frac{\left|\left\langle\bm{\Gamma}^{\infty}(\cdot,\bm{z},\bm{p}), \bm{\psi}_j\right\rangle\right|^2_{\mathcal{L}^2(\mathbb{S})}}{\left|\lambda_j\right|}
		$$
		converges, i.e., if and only if
		\begin{equation}\label{W_z}
			W(\bm{z}):=\left[\sum_{j=1}^{\infty} \frac{\left|\left\langle\bm{\Gamma}^{\infty}(\cdot,\bm{z},\bm{p}), \bm{\psi}_j\right\rangle_{\mathcal{L}^2(\mathbb{S})}\right|^2}
			{\left|\lambda_j\right|}\right]^{-1}>0,
		\end{equation}
		where $\lambda_j \in \mathbb{C}$ are the eigenvalues of the operator $\bm{\mathcal{F}}$ with corresponding eigenfunctions $\bm{\psi}_j \in \mathcal{L}^2(\mathbb{S})$.
	\end{theorem}
	\begin{proof}
		Since the operator $\bm{\mathcal{F}}$ is compact and normal, by the spectral theorem we have the expansion
		\begin{equation}
			\bm{\mathcal{F}} \bm{g}=\sum_{j=1}^{\infty} \lambda_j\left\langle\bm{g}, \bm{\psi}_j\right\rangle \bm{\psi}_j, \quad \bm{g} \in \mathcal{L}^2(\mathbb{S}),
		\end{equation}
		which gives
		\begin{equation}
			\left(\bm{\mathcal{F}}^*\bm{\mathcal{F}}\right)^{1/4} \bm{g}=\sum_{j=1}^{\infty}\sqrt{|\lambda_j|}\left\langle\bm{g}, \bm{\psi}_j\right\rangle \bm{\psi}_j, \quad \bm{g} \in \mathcal{L}^2(\mathbb{S}).
		\end{equation}
		Thus, we know that $\{\sqrt{|\lambda_j|}, \bm{\psi}_j\}$
		is the eigensystem of $\left(\bm{\mathcal{F}}^*\bm{\mathcal{F}}\right)^{1/4} $.
		It follows from Lemmas \ref{z_range} and \ref{G_range} that $\bm{z}\in\mathbb{R}^2$ belongs to $D$ if and only if
		the equation $(\ref{FFg_equ})$ is solvable in $\mathcal{L}^2(\mathbb{S})$. Utilizing the Picard theorem (see Theorem 4.8 \cite{CK13}), we can obtain that $(\ref{FFg_equ})$ is solvable in $\mathcal{L}^2(\mathbb{S})$ if and only if the series
		\begin{equation}\label{4.5_1}
			\sum_{j=1}^{\infty} \frac{|\left\langle\boldsymbol{\Gamma}^{\infty}(\cdot, \boldsymbol{z}, \boldsymbol{p}), \boldsymbol{\psi}_j\right
				\rangle_{\mathcal{L}^2(\mathbb{S})}|^2}
			{\left|\lambda_j\right|}<\infty,
		\end{equation}
		in this case, the solution of equation $(\ref{FFg_equ})$
		is given by
		\begin{equation*}
			\bm{g}=\sum_{j=1}^{\infty} \frac{\left\langle\boldsymbol{\Gamma}^{\infty}(\cdot, \boldsymbol{z}, \boldsymbol{p}), \boldsymbol{\psi}_j\right
				\rangle_{\mathcal{L}^2(\mathbb{S})}}
			{\sqrt{\left|\lambda_j\right|}}\boldsymbol{\psi}_j.
		\end{equation*}
		The equivalence of (\ref{4.5_1}) and (\ref{W_z}) is obvious.
	\end{proof}
	
	\begin{remark}
		It is derived from the equation $\bm{\mathcal{G}}\left(\bm{\Gamma}(\cdot,\bm{z};\bm{p})\right)
			=\bm{\Gamma}^{\infty}(\cdot,\bm{z};\bm{p})$ that
			$\bm{g}=(\bm{\mathcal{F}}^*\bm{\mathcal{F}})^{-1/4}
			\bm{\mathcal{G}}\left
			(\bm{\Gamma}(\cdot,\bm{z};\bm{p})\right)$.
			The uniqueness of the inverse scattering problem can be directly derived as a result of the isomorphism of the operator $\left(\bm{\mathcal{F}}^*\bm{\mathcal{F}}\right)^{-1 / 4} \bm{\mathcal{G}}$. Therefore, the factorization
			method can be considered a novel proof for demonstrating the result of uniqueness.
	\end{remark}
	
	\subsection{Reconstruction for PP case and SS case}
	In this subsection, we consider the PP and SS cases. Following the idea presented in \cite{HKS12} for Dirichlet boundary problem, we define two projection spaces
	\begin{align*}
		\mathcal{L}_p^2(\mathbb{S}):&=\left\{g_p\in L^2(\mathbb{S}):\ g_p(\bm{d})= \bm{g}(\bm{d})\bm{d},~\bm{g}\in \mathcal{L}^2(\mathbb{S}),\bm{d}\in \mathbb{S}\right\},
		\\
		\mathcal{L}_s^2(\mathbb{S}):&=\left\{g_s\in L^2(\mathbb{S}):\ g_s(\bm{d})= \bm{g}(\bm{d})\bm{d}^\bot,~\bm{g}\in \mathcal{L}^2(\mathbb{S}),\bm{d}\in \mathbb{S}\right\},
	\end{align*}
	and define the orthogonal projection operators $\bm{\mathcal{P}}_p: \mathcal{L}^2(\mathbb{S})\rightarrow \mathcal{L}_p^2(\mathbb{S})$ and
	$\bm{\mathcal{P}}_s: \mathcal{L}^2(\mathbb{S})\rightarrow \mathcal{L}_s^2(\mathbb{S})$ by
	\begin{equation*}
		\bm{\mathcal{P}}_p \bm{g}(\bm{d}):=\bm{g}_p(\bm{d}),\quad \bm{\mathcal{P}}_s\bm{g}(\bm{d}):=\bm{g}_s(\bm{d}).
	\end{equation*}
	Thus, we have far field operator $\bm{\mathcal{F}}_m:\mathcal{L}_m^2(\mathbb{S})\rightarrow \mathcal{L}_m^2(\mathbb{S})$, defined by
	\begin{equation*}
		\bm{\mathcal{F}}_m:=\bm{\mathcal{P}}_m\bm{\mathcal{F}} \bm{\mathcal{P}}_m^*,\quad m\in\{p,s\},
	\end{equation*}
	where
	$\bm{\mathcal{P}}_m^*:\mathcal{L}_m^2(\mathbb{S})\rightarrow \mathcal{L}^2(\mathbb{S})$ is the adjoint operator of $\bm{\mathcal{P}}_m$. From (\ref{F factori}), $\bm{\mathcal{F}}_m$ has the following factorization form
	\begin{equation}\label{Fm_factor}
		\bm{\mathcal{F}}_m=-\sqrt{8\pi\omega}\left(\bm{\mathcal{P}}_m \bm{\mathcal{G}}\right) \bm{\mathcal{N}}^*\left(\bm{\mathcal{P}}_m\bm{\mathcal{G}}\right)^* ,\quad m\in\{p,s\}.
	\end{equation}
	Since the operator $\bm{\mathcal{F}}_m$ fails to be normal, Theorem 1.23 in \cite{KAG07} is not applicable.
	We factorize the auxiliary self-adjoint positive operator
	$$
	\bm{\mathcal{F}}_m^{\#}=|\operatorname{Re}\bm{\mathcal{F}}_m
	|+|\operatorname{Im} \bm{\mathcal{F}}_m|
	$$
	rather than directly utilizing $\bm{\mathcal{F}}_m$; see also \cite{KAG07}. Here, the real and imaginary parts of $\bm{\mathcal{F}}_m$ are self-adjoint operators defined by
	$$
	\operatorname{Re} \bm{\mathcal{F}}_m=\frac{1}{2}\left(\bm{\mathcal{F}}_m
	+\bm{\mathcal{F}}_m^*\right) \quad \text { and } \quad \operatorname{Im} \bm{\mathcal{F}}_m=\frac{1}{2 \mathrm{i}}\left(\bm{\mathcal{F}}_m
	-\bm{\mathcal{F}}_m^*\right),\quad m\in\{p,s\}.
	$$
	
	\begin{lemma}\label{zm_range}
		For any $\bm{z} \in \mathbb{R}^2$,
		$\bm{p}\in\mathbb{S}$, and $m\in\{p,s\}$, $\bm{z} \in D$ if and only if
		$\bm{\Gamma}_m^{\infty}(\cdot,\bm{z},\bm{p})
		$ belongs to the range of $\bm{\mathcal{P}}_m\bm{\mathcal{G}}$.
		
	\end{lemma}
	\begin{proof}
		Assume that $\boldsymbol{z} \in D$, and define $\boldsymbol{v}^{\mathrm{s}}(\boldsymbol{x}):
		=\boldsymbol{\Gamma}(\boldsymbol{x}, \boldsymbol{z}) \boldsymbol{p}$ and $\boldsymbol{f}:=\left.\boldsymbol{T}_{\boldsymbol{\nu}} \boldsymbol{v}^{\mathrm{s}}\right|_{\partial D}$.
		It follows from the same argument as the proof of Lemma \ref{z_range} that $\bm{f}$
		belongs to $\left[H^{-1 / 2}(\partial D)\right]^2$ and
		the far-field pattern of $\boldsymbol{v}^{\mathrm{s}}$
		is $\boldsymbol{\Gamma}^{\infty}(\cdot, \boldsymbol{z} ; \boldsymbol{p})$. From the definition of
		$\bm{\mathcal{P}}_m$, we have $\left(\bm{\mathcal{P}}_m\bm{\mathcal{G}}\right)\bm{f}=
		\bm{\mathcal{P}}_m\boldsymbol{\Gamma}^{\infty}(\cdot, \boldsymbol{z} ; \boldsymbol{p})=\boldsymbol{\Gamma}_m^{\infty}(\cdot, \boldsymbol{z} ; \boldsymbol{p})$, i.e., $\boldsymbol{\Gamma}_m^{\infty}(\cdot, \boldsymbol{z} ; \boldsymbol{p})\in \mathcal{R}\left(\bm{\mathcal{P}}_m\bm{\mathcal{G}}\right)$.

		For $\boldsymbol{z}\in \mathbb{R}^2 \backslash \overline{D}$, assume that there exists $\tilde{\boldsymbol{f}}\in\left[H^{-1 / 2}(\partial D)\right]^2$ such that $\left(\bm{\mathcal{P}}_m\bm{\mathcal{G}}\right) \tilde{\boldsymbol{f}}=\boldsymbol{\Gamma}_m^{\infty}(\cdot, \boldsymbol{z} ; \boldsymbol{p})$. Let $\boldsymbol{v}^s$ be the solution of the exterior Neumann problem with boundary data $\tilde{\boldsymbol{f}}$, and $\bm{v}^{\infty}$ the far field pattern of $\boldsymbol{v}^s$. Denote by $\bm{v}^{\infty}=\bm{v}_p^\infty+\bm{v}_s^\infty$. One immediately has $\boldsymbol{v}_m^{\infty}=\left(\boldsymbol{\mathcal { P}}_m\bm{\mathcal{G }} \right) \tilde{\boldsymbol{f}}=\boldsymbol{\Gamma}_m^{\infty}(\cdot, \boldsymbol{z} ; \boldsymbol{p})$.
		Since $\boldsymbol{\Gamma}_m^{\infty}(\cdot, \boldsymbol{z} ; \boldsymbol{p})$ is the far-field pattern of $\boldsymbol{\Gamma}_m(\cdot, \boldsymbol{z}) \boldsymbol{p}$, it follows from the Rellich's lemma that $\boldsymbol{v}_m^s(\boldsymbol{x})=$ $\boldsymbol{\Gamma}_m(\bm{x}, \boldsymbol{z}) \boldsymbol{p}$ in $\mathbb{R}^2 \backslash(\overline{D} \cup\{\bm{z}\})$. However, the scattered field $\boldsymbol{v}_m^s(\boldsymbol{x})$ is analytic in $\mathbb{R}^2 \backslash \overline{D}$ while $\boldsymbol{\Gamma}_m(\bm{x}, \boldsymbol{z}) \boldsymbol{p}$ has a singularity at $\boldsymbol{z}$. This leads to a contradiction.
		
		For $\boldsymbol{z} \in \partial D$, using the same argument as the proof of Lemma \ref{z_range}
		also yields the contradiction.
	\end{proof}
	
	\begin{lemma}\label{PG_range}
		Assume that $\omega^2$ is not a Neumann eigenvalue of $-\Delta^*$ in $D$. For $m\in\{p,s\}$, the ranges of $\bm{\mathcal{P}}_m \bm{\mathcal{G}}$ and $\left(\bm{\mathcal{F}}_{m}^{\#}\right)^{1 / 2}$ coincide, where $\bm{\mathcal{F}}_{m}^{\#}:=\left|\operatorname{Re} \bm{\mathcal{F}}_m\right|+\left|\operatorname{Im} \bm{\mathcal{F}}_m\right|$.
	\end{lemma}
	\begin{proof}
		By taking $Y=\mathcal{L}_m^2(\mathbb{S})$, $X=\left[H^{-1 / 2}(\partial D)\right]^2$, $G=\bm{\mathcal{P}}_m\bm{\mathcal{G}}$ and $T=-\sqrt{8 \pi \omega}\bm{\mathcal{N}}^*$ in Theorem 2.15 in \cite{KAG07}, and using Lemmas \ref{G_properties} and \ref{N properties}, we immediately complete the proof.
	\end{proof}
	
	Combining Lemmas \ref{G_properties}, \ref{zm_range}, and \ref{PG_range} together we arrive at the
	following final characterization of the scatterer $D$.

	\begin{theorem}\label{Fm_Th}
		Assume that $\omega^2$ is not a Neumann eigenvalue of $-\Delta^*$ in $D$. Then for $m\in\{p,s\}$,
		a point $\bm{z} \in \mathbb{R}^2$ belongs to $D$ if and only if the series
		\begin{align}
			\sum_{j=1}^{\infty} \frac{\left|\left\langle\bm{\Gamma}^{\infty}_m(\cdot,\bm{z},\bm{p}), \bm{\psi}_{m,j}\right\rangle_{\mathcal{L}_m^2(\mathbb{S})}
				\right|^2}{\lambda_{m,j}}
		\end{align}
		converges, i.e., if and only if
		\begin{align}\label{4.7_1}
			W_m(\bm{z}):=\left[\sum_{j=1}^{\infty} \frac{\left|\left\langle\bm{\Gamma}^{\infty}_m(\cdot,\bm{z},\bm{p}), \bm{\psi}_{m,j}\right\rangle_{\mathcal{L}_m^2(\mathbb{S})}
				\right|^2}{\lambda_{m,j}}\right]^{-1}>0,
		\end{align}
		where $\lambda_{m,j} \in \mathbb{C}$ are the eigenvalues of the operator $\bm{\mathcal{F}}_m^{\#}$ with corresponding eigenfunctions $\bm{\psi}_{m,j} \in \mathcal{L}_m^2(\mathbb{S})$.
	\end{theorem}
	
	\begin{proof}
		From Lemmas \ref{zm_range} and \ref{PG_range}, we know that $\bm{z}\in \mathbb{R}^2$ belongs to $D$ if and only if the equation
		\begin{equation}\label{F_p}
			\left(\bm{\mathcal{F}}_{m}^{\#}\right)^{1 / 2}\bm{g}_m(\cdot,\bm{z};\bm{p})=
			\bm{\Gamma}_m^\infty(\cdot,\bm{z};\bm{p})
		\end{equation}
		has a solution in $\mathcal{L}_m^2(\mathbb{S})$.
		Using the Picard theorem (see Theorem 4.8 \cite{CK13}), we can obtain that equation (\ref{F_p}) is solvable in $\mathcal{L}_m^2(\mathbb{S})$ if and only if the series
		\begin{equation}\label{4.7_2}
			\sum_{j=1}^{\infty} \frac{\left\langle\boldsymbol{\Gamma}_m^{\infty}(\cdot, \boldsymbol{z}, \boldsymbol{p}), \boldsymbol{\psi}_{m,j}\right
				\rangle_{\mathcal{L}_m^2(\mathbb{S})}^2}
			{\lambda_{m,j}}<\infty.
		\end{equation}
		In this case, the solution of equation (\ref{F_p}) is given by
		$$
		\boldsymbol{g}=\sum_{j=1}^{\infty} \frac{\left\langle\boldsymbol{\Gamma}_m^{\infty}(\cdot, \boldsymbol{z}, \boldsymbol{p}), \boldsymbol{\psi}_{m,j}\right
			\rangle_{\mathcal{L}_m^2(\mathbb{S})}}
		{\sqrt{\lambda_{m,j}}} \boldsymbol{\psi}_{m,j}.
		$$
		The equivalence of (\ref{4.7_1}) and (\ref{4.7_2}) is obvious.
	\end{proof}

	\subsection{Reconstruction for the limited aperture case}
	In this subsection, we consider the limited aperture case
	which only uses the partial data of the far field pattern, that is, $\bm{u}^{\infty}(\hat{\bm{x}},\bm{d})$ for $\hat{\bm{x}},\bm{d}\in \mathbb{U}$. Here, $\mathbb{U}$ denotes a non-empty subset of $\mathbb{S}$.
	
	We define
	the limited aperture far field operator $\bm{\mathcal{F}}_{la}: \mathcal{L}^2(\mathbb{U}) \longrightarrow \mathcal{L}^2(\mathbb{U})$ by
	\begin{equation}
		\bm{\mathcal{F}}_{la} \boldsymbol{g}(\hat{\boldsymbol{x}})=e^{-\frac{\mathrm{i} \pi}{4}} \int_{\mathbb{U}}\left\{\sqrt{\frac{k_p}{\omega}} \boldsymbol{u}^{\infty}(\hat{\boldsymbol{x}} ; \boldsymbol{d}; p) g_{p}(\boldsymbol{d})+\sqrt{\frac{k_s}{\omega}} \boldsymbol{u}^{\infty}\left(\hat{\boldsymbol{x}} ; \boldsymbol{d};s\right) g_{s}(\boldsymbol{d})\right\} \mathrm{d} s(\boldsymbol{d}), \quad \hat{\boldsymbol{x}} \in \mathbb{U},
	\end{equation}
	and the restriction operator $\bm{\mathcal{P}}_{la}:\mathcal{L}^2(\mathbb{S})\rightarrow
	\mathcal{L}^2(\mathbb{U})$ by
	\begin{equation*}
		\bm{\mathcal{P}}_{la} \bm{g}(\bm{d}):=\bm{g}(\bm{d})|_{\mathbb{U}}.
	\end{equation*}
	The adjoint $\bm{\mathcal{P}}_{la}^*:\mathcal{L}^2(\mathbb{U})\rightarrow \mathcal{L}^2(\mathbb{S})$ is given by the extension
	such that $\bm{\mathcal{P}}_{la}^*\bm{g}(\bm{d})=\bm{g}(\bm{d})$ for $\bm{d}\in \mathbb{U}$, and $\bm{\mathcal{P}}_{la}^*\bm{g}(\bm{d})=\bm{0}$ otherwise.
	Thus, the definition of $\bm{\mathcal{F}}_{la}$ can be rewritten as
	\begin{equation}
		\bm{\mathcal{F}}_{la}:=\bm{\mathcal{P}}_{la} \bm{\mathcal{F}}\bm{\mathcal{P}}_{la}^*.
	\end{equation}
	From (\ref{F factori}), $\bm{\mathcal{F}}_{la}$ has the following factorization form
	\begin{equation}\label{Fm_factor}
		\bm{\mathcal{F}}_{la}=-\sqrt{8\pi\omega}
		\left(\bm{\mathcal{P}}_{la}\bm{\mathcal{G}}\right) \bm{\mathcal{N}}^*\left(\bm{\mathcal{P}}_{la}
		\bm{\mathcal{G}}\right)^* .
	\end{equation}
	It follows from the injectivity, compactness and denseness of the ranges of $\bm{\mathcal{G}}$ that $\bm{\mathcal{P}}_{la}\bm{\mathcal{G}}$ is injective, compact and have dense ranges in $\bm{\mathcal{L}}^2(\mathbb{U})$.
	
	\begin{lemma}\label{4.8}
		For any $\boldsymbol{z} \in \mathbb{R}^2$ and $\boldsymbol{p} \in \mathbb{S}$, $\boldsymbol{z} \in D$ if and only if $\boldsymbol{\Gamma}^{\infty}(\cdot, \boldsymbol{z}, \boldsymbol{p})|_{\mathbb{U}}=\bm{\mathcal{P}}_{la}
		\boldsymbol{\Gamma}^{\infty}(\cdot, \boldsymbol{z}, \boldsymbol{p})$ belongs to the range of $\bm{\mathcal{P}}_{la}\bm{\mathcal{G}}$.
	\end{lemma}
	\begin{proof}
		From Lemma \ref{z_range}, we know that $\bm{z} \in D$ if and only if $\boldsymbol{\Gamma}^{\infty}(\cdot, \boldsymbol{z}, \boldsymbol{p})$ belongs to the range of
		$\bm{\mathcal{G}}$. This immediately implies that $\bm{\mathcal{P}}_{la}
		\boldsymbol{\Gamma}^{\infty}(\cdot, \boldsymbol{z}, \boldsymbol{p})$ belongs to the range of $\bm{\mathcal{P}}_{la}\bm{\mathcal{G}}$ if $\bm{z} \in D$.
		Conversely, if $\bm{\mathcal{P}}_{la}
		\boldsymbol{\Gamma}^{\infty}(\cdot, \boldsymbol{z}, \boldsymbol{p})$ belongs to the range of $\bm{\mathcal{P}}_{la}\bm{\mathcal{G}}$, one has
		$\boldsymbol{\Gamma}^{\infty}(\cdot, \boldsymbol{z}, \boldsymbol{p})=\bm{v}^\infty$ on $\mathbb{U}$
		for far field pattern $\bm{v}^{\infty}=\bm{\mathcal{ G}}\bm{f}$. Thus, the analyticity of $\boldsymbol{\Gamma}^{\infty}(\cdot, \boldsymbol{z}, \boldsymbol{p})-\bm{\mathcal{ G}}\bm{f}$ implies that $\bm{z}\in D$.
	\end{proof}

	\begin{lemma}\label{4.9}
		Assume that $\omega^2$ is not a Neumann eigenvalue of $-\Delta^*$ in $D$. The ranges of $\bm{\mathcal{P}}_{la} \bm{\mathcal{G}}$ and $\left(\bm{\mathcal{F}}_{la}^{\#}\right)^{1 / 2}$ coincide, where $\bm{\mathcal{F}}_{la}^{\#}:=\left|\operatorname{Re} \bm{\mathcal{F}}_{la}\right|+\left|\operatorname{Im} \bm{\mathcal{F}}_{la}\right|$.
	\end{lemma}
	\begin{proof}
		By taking $Y=\mathcal{L}^2(\mathbb{U})$, $X=\left[H^{-1 / 2}(\partial D)\right]^2$, $G=\bm{\mathcal{P}}_{la}\bm{\mathcal{G}}$ and $T=-\sqrt{8 \pi \omega}\bm{\mathcal{N}}^*$ in Theorem 2.15 in \cite{KAG07}, and using Lemmas \ref{G_properties} and \ref{N properties}, we immediately complete the proof.
	\end{proof}
	
	Combining Lemmas \ref{G_properties}, \ref{4.8} and \ref{4.9} together yields the
	following final characterization of the scatterer $D$.
	
	\begin{theorem}
		Assume that $\omega^2$ is not a Neumann eigenvalue of $-\Delta^*$ in $D$. Then a point $\boldsymbol{z} \in \mathbb{R}^2$ belongs to $D$ if and only if the series
		$$
		\sum_{j=1}^{\infty} \frac{\left|\left\langle\boldsymbol{\Gamma}^{\infty}(\cdot, \boldsymbol{z}, \boldsymbol{p}), \boldsymbol{\psi}_{la, j}\right\rangle_{\mathcal{L}^2(\mathbb{U})}
			\right|^2}{\lambda_{la, j}}
		$$
		converges, i.e., if and only if
		$$
		W_m(\boldsymbol{z}):=\left[\sum_{j=1}^{\infty} \frac{\left|\left\langle\boldsymbol{\Gamma}^{\infty}(\cdot, \boldsymbol{z}, \boldsymbol{p}), \boldsymbol{\psi}_{la, j}\right\rangle_{\mathcal{L}^2(\mathbb{U})}\right|^2}
		{\lambda_{la, j}}\right]^{-1}>0,
		$$
		where $\lambda_{la, j} \in \mathbb{C}$ are the eigenvalues of the operator $\bm{\mathcal{F}}_{la}^{\#}$ with corresponding eigenfunctions $\boldsymbol{\psi}_{la, j} \in$ $\mathcal{L}^2(\mathbb{U})$.
	\end{theorem}
	\begin{proof}
		The proof is similar to the proof of Theorem \ref{Fm_Th}.
	\end{proof}

	\section{Numerical experiments}\label{Section_5}
	In this section, several numerical examples are presented to demonstrate the performance of the proposed reconstruction algorithm. We select a point set $\mathbb{Z}\subset \mathbb{R}^2$ such that the unknown cavity is contained within the convex hull of these points. We employ Nystr\"{o}m-type
	discretization to approximate the boundary integral equation, which enables us to generate synthetic far-field data. Subsequently, we add noise to each exact measurement by
	$$
	\textbf{F}^\delta=\textbf{F}+\delta\|\textbf{F}\| \frac{\textbf{R}_1+\textbf{R}_2 \mathrm{i}}{\left\|\textbf{R}_1+\textbf{R}_2 \mathrm{i}\right\|},
	$$
	where $\textbf{R}_1$ and $\textbf{R}_2$ are two random matrices produced by the MATLAB function $randn (N, N)$, and $\delta>0$ is the relative noise level. Let $\mathbb{D}=\{\bm{d}_{j}=(\cos\theta_j,\sin\theta_j)^\top;
	\ \theta_j=2\pi j/N, j=1,2,\dots,N\}$ be a set of $N$ equidistantly distributed directions.
	
	
	\textbf{FF case:} We sequentially use a plane compression wave and a plane shear wave as incident waves in each direction $\bm{d}_{j}$, $j=1,2,\dots,N$, i.e., a total of $2N$ incident waves. For the illumination in direction $\bm{d}_j$, we collect far-field data
	\begin{equation*}
		\textbf{F}_{k,j}=\left(\begin{array}{cc}
			u_p^{\infty}(\bm{d}_k, \bm{d}_j;p) & u_p^{\infty}(\bm{d}_k, \bm{d}_j;s) \\
			u_s^{\infty}(\bm{d}_k, \bm{d}_j;p) & u_s^{\infty}(\bm{d}_k, \bm{d}_j;s)
		\end{array}\right),
	\end{equation*}
	where $u_p^{\infty}(\bm{d}_k, \bm{d}_j;p)\bm{d}_k+u_s^{\infty}(\bm{d}_k, \bm{d}_j;p)\bm{d}_k^\bot=\bm{u}^{\infty}(\bm{d}_k, \bm{d}_j;p)$ and $u_p^{\infty}(\bm{d}_k, \bm{d}_j;s)\bm{d}_k+u_s^{\infty}(\bm{d}_k, \bm{d}_j;s)\bm{d}_k^\bot=\bm{u}^{\infty}(\bm{d}_k, \bm{d}_j;s)$.
	
	Using the trapezoidal rule to  discrete equation (\ref{Fg}) yields
	\begin{equation}\notag
		\begin{aligned}
			&\frac{2 \pi e^{-\frac{\mathrm{i} \pi}{4}}}{N}  \sum_{j=1}^N\left\{\sqrt{\frac{k_p}{\omega}} u_p^{\infty}\left(\bm{d}_k, \bm{d}_j ; p\right) g_p\left( \bm{d}_j,\bm{z};\bm{p}\right)
			+\sqrt{\frac{k_s}{\omega}} u_p^{\infty}\left( \bm{d}_k,  \bm{d}_j ; s\right) g_s\left(\bm{d}_j,\bm{z};\bm{p}\right)\right\}\\
			=&\frac{1}{2 \mu+\lambda} \frac{e^{\frac{\mathrm{i} \pi}{4}}}{\sqrt{8 \pi k_p}} e^{-\mathrm{i} k_p \bm{d}_k\cdot \bm{z}} \bm{d}_k \cdot \bm{p}, \quad k=1, \cdots, N,
		\end{aligned}
	\end{equation}
	\begin{equation}\notag
		\begin{aligned}
			&\frac{2 \pi e^{-\frac{\mathrm{i} \pi}{4}}}{N}  \sum_{j=1}^N\left\{\sqrt{\frac{k_p}{\omega}} u_s^{\infty}\left(\bm{d}_k, \bm{d}_j ; p\right) g_p\left( \bm{d}_j,\bm{z};\bm{p}\right)+\sqrt{\frac{k_s}{\omega}} u_s^{\infty}\left( \bm{d}_k,  \bm{d}_j ; s\right) g_s\left(\bm{d}_j,\bm{z};\bm{p}\right)\right\}\\
			=&\frac{1}{\mu} \frac{e^{\frac{\mathrm{i} \pi}{4}}}{\sqrt{8 \pi k_s}} e^{-\mathrm{i} k_s \bm{d}_k\cdot \bm{z}} \bm{d}^\bot_k \cdot \bm{p}, \quad k=1, \cdots, N.
		\end{aligned}
	\end{equation}
	Rewrite the above equations in matrix form as
	\begin{equation}
		\textbf{F}\textbf{g}^{(\bm{z},\bm{p})}
		=\bm{h}^{(\bm{z},\bm{p})},
	\end{equation}
	where
	\begin{equation*}
		\textbf{F}=\frac{2 \pi e^{-\frac{\mathrm{i} \pi}{4}}}{N}\left(\begin{array}{ll}
			\sqrt{\frac{k_p}{\omega}}\left(u_p^{\infty}\left(\bm{d}_k, \bm{d}_j;p\right)\right)_{k, j} & \sqrt{\frac{k_s}{\omega}}\left(u_p^{\infty}\left(\bm{d}_k, \bm{d}_j;s\right)\right)_{k, j} \\
			\sqrt{\frac{k_p}{\omega}}\left(u_s^{\infty}\left(\bm{d}_k, \bm{d}_j;p\right)\right)_{k, j} & \sqrt{\frac{k_s}{\omega}}\left(u_s^{\infty}\left(\bm{d}_k, \bm{d}_j;s\right)\right)_{k, j}
		\end{array}\right),\quad k,j=1,2, \ldots, N.
	\end{equation*}
	and
	\begin{equation*}
		\textbf{g}^{(\bm{z},\bm{p})}=\left(\begin{array}{c}
			g_p(\bm{d}_1,\bm{z};\bm{p})\\
			\vdots\\
			g_p(\bm{d}_N,\bm{z};\bm{p})\\
			g_s(\bm{d}_1,\bm{z};\bm{p})\\
			\vdots\\
			g_s(\bm{d}_N,\bm{z};\bm{p})
		\end{array}\right),
		\quad
		\bm{h}^{(\bm{z},\bm{p})}:=\left(\begin{array}{c}
			\frac{1}{2 \mu+\lambda} \frac{e^{\frac{\mathrm{i} \pi}{4}}}{\sqrt{8 \pi k_p}} e^{-\mathrm{i} k_p \bm{d}_k\cdot \bm{z}} \bm{d}_1 \cdot \bm{p}\\
			\vdots\\
			\frac{1}{2 \mu+\lambda} \frac{e^{\frac{\mathrm{i} \pi}{4}}}{\sqrt{8 \pi k_p}} e^{-\mathrm{i} k_p \bm{d}_k\cdot \bm{z}} \bm{d}_N \cdot \bm{p}\\
			\frac{1}{\mu} \frac{e^{\frac{\mathrm{i} \pi}{4}}}{\sqrt{8 \pi k_s}} e^{-\mathrm{i} k_s \bm{d}_k\cdot \bm{z}} \bm{d}^\bot_1 \cdot \bm{p}\\
			\vdots\\
			\frac{1}{\mu} \frac{e^{\frac{\mathrm{i} \pi}{4}}}{\sqrt{8 \pi k_s}} e^{-\mathrm{i} k_s \bm{d}_k\cdot \bm{z}} \bm{d}^\bot_N \cdot \bm{p}
		\end{array}\right).
	\end{equation*}
	We compute the singular value decomposition $\textbf{F}=\textbf{U}\textbf{S}\textbf{V}^*$, where the
	diagonal elements of $\textbf{S}$ are denoted by $\lambda_j$, $j=1,2,\ldots,2N$. For each $\bm{z}\in\mathbb{Z}$ , we compute
	\begin{equation*}
		W(\bm{z}):=\left[\sum_{j=1}^N \frac{\left|\rho_{j}^{(\bm{z},\bm{p})}\right|^2}{\left|\lambda
			_{j}\right|}\right]^{-1},
	\end{equation*}
	where $\bm{\rho}^{(\bm{z,p})}=\textbf{V}^*\bm{h}^{(\bm{z},\bm{p})}$.
	The values of $W(\bm{z})$ for $\bm{z} \notin D$ should be much smaller than that for $\bm{z}\in D$.
	
	\textbf{PP case or SS case:} We use a plane compression wave or a plane shear wave as incident wave in each direction $\boldsymbol{d}_j, j=1,2, \ldots, N$, i.e., a total of $N$ incident waves. For the illumination in direction $\boldsymbol{d}_j$, we collect far-field data
	$
	\mathbf{F}_{k, j}=\left(
	u_p^{\infty}\left(\boldsymbol{d}_k, \boldsymbol{d}_j ; p\right)\right)^N_{k,j=1}
	$
	or
	$
	\mathbf{F}_{k, j}=\left(
	u_s^{\infty}\left(\boldsymbol{d}_k, \boldsymbol{d}_j ; p\right)\right)^N_{k,j=1}
	$.
	
	From the definition of $\bm{\mathcal{F}}_m$, we have
	\begin{equation}\label{eq5.2}
		\bm{\mathcal{F}}_m\bm{g}_m(\cdot,\bm{z};\bm{p})=\boldsymbol{\Gamma}
		_m^{\infty}(\cdot, \boldsymbol{z}; \boldsymbol{p}),\quad m=p\ \text{or}\ s.
	\end{equation}
	Using the trapezoidal rule to discrete equation (\ref{eq5.2}) yields
	\begin{equation}
		\begin{aligned}
			\frac{2 \pi e^{-\frac{\mathrm{i} \pi}{4}}}{N} \sum_{j=1}^N\sqrt{\frac{k_p}{\omega}} u_p^{\infty}\left(\boldsymbol{d}_k, \boldsymbol{d}_j ; p\right) g_p\left(\boldsymbol{d}_j, \boldsymbol{z} ; \boldsymbol{p}\right)=\frac{1}{2 \mu+\lambda} \frac{e^{\frac{\mathrm{i} \pi}{4}}}{\sqrt{8 \pi k_p}} e^{-\mathrm{i} k_p \boldsymbol{d}_k \cdot \boldsymbol{z}} \boldsymbol{d}_k \cdot \boldsymbol{p}, \quad k=1, \cdots, N,
		\end{aligned}
	\end{equation}
	or
	\begin{equation}
		\begin{aligned}
			\frac{2 \pi e^{-\frac{\mathrm{i} \pi}{4}}}{N} \sum_{j=1}^N\sqrt{\frac{k_s}{\omega}} u_s^{\infty}\left(\boldsymbol{d}_k, \boldsymbol{d}_j ; s\right) g_s\left(\boldsymbol{d}_j, \boldsymbol{z} ; \boldsymbol{p}\right)
			=  \frac{1}{\mu} \frac{e^{\frac{\mathrm{i} \pi}{4}}}{\sqrt{8 \pi k_s}} e^{-\mathrm{i} k_s \boldsymbol{d}_k \cdot \boldsymbol{z}} \boldsymbol{d}_k^{\perp} \cdot \boldsymbol{p}, \quad k=1, \cdots, N .
		\end{aligned}
	\end{equation}
	Rewrite the above equations in matrix form as
	\begin{equation}
		\textbf{F}_m\textbf{g}_m^{(\bm{z},\bm{p})}=
		\bm{h}_m^{(\bm{z},\bm{p})},\quad m=p\ \text{or}\ s,
	\end{equation}
	where $\textbf{g}_m^{(\bm{z},\bm{p})}=
	\left(g_m(\bm{d}_1,\bm{z};\bm{p}),\cdots
	,g_m(\bm{d}_N,\bm{z};\bm{p})\right)^\top$,
	\begin{equation*}
		\textbf{F}_p:=\frac{2 \pi e^{-\frac{\mathrm{i} \pi}{4}}}{N}\left(\left(\sqrt{\frac{k_p}{\omega}}u_p^{\infty}\left(\bm{d}_k, \bm{d}_j;\bm{d}_j\right)\right)_{k, j}\right),\quad \textbf{F}_s:=\frac{2 \pi e^{-\frac{\mathrm{i} \pi}{4}}}{N}\left(\left(\sqrt{\frac{k_s}{\omega}}u_s^{\infty}
		\left(\bm{d}_k, \bm{d}_j;\bm{d}_j^\bot\right)\right)_{k,j}\right),
	\end{equation*}
	and
	\begin{equation*}
		\boldsymbol{h}_p^{(\boldsymbol{z}, \boldsymbol{p})}:=\frac{1}{2 \mu+\lambda} \frac{e^{\frac{\mathrm{i} \pi}{4}}}{\sqrt{8 \pi k_p}} \left(\begin{array}{c}
			e^{-\mathrm{i} k_p \boldsymbol{d}_1 \cdot \boldsymbol{z}} \boldsymbol{d}_1 \cdot \boldsymbol{p}\\
			\vdots\\
			e^{-\mathrm{i} k_p \boldsymbol{d}_N \cdot \boldsymbol{z}} \boldsymbol{d}_N \cdot \boldsymbol{p}
		\end{array}\right),
		\quad
		\boldsymbol{h}_s^{(\boldsymbol{z}, \boldsymbol{p})}:=\left(\begin{array}{c}
			\frac{1}{\mu} \frac{e^{\frac{\mathrm{i} \pi}{4}}}{\sqrt{8 \pi k_s}} e^{-\mathrm{i} k_s \boldsymbol{d}_1 \cdot \boldsymbol{z}} \boldsymbol{d}_1^{\perp} \cdot \boldsymbol{p}\\
			\vdots\\
			\frac{1}{\mu} \frac{e^{\frac{\mathrm{i} \pi}{4}}}{\sqrt{8 \pi k_s}} e^{-\mathrm{i} k_s \boldsymbol{d}_N \cdot \boldsymbol{z}} \boldsymbol{d}_N^{\perp} \cdot \boldsymbol{p}
		\end{array}\right).
	\end{equation*}
	We compute singular value decomposition $\mathbf{F}_m=\textbf{U}_m \textbf{S}_m \textbf{V}_m^*$, where the diagonal elements of $\mathbf{S}_m$ are denoted by $\eta_{m,j}$, $j=1,2, \ldots, N$. Since the operators $\bm{\mathcal{F}}_m$ fails to be normal, we take $\lambda_{m,j}=|\mathrm{Re}\eta_{m,j}|
	+|\mathrm{Im}\eta_{m,j}|$; see also
	\cite{HKS12}.
	For each $\bm{z}\in \mathbb{Z}$ and $m=p$ or $s$, we compute
	$$
	W_m(\boldsymbol{z}):=\left[\sum_{j=1}^N \frac{\left|\rho_{m,j}^{(\boldsymbol{z}, \boldsymbol{p})}\right|^2}
	{\lambda_{m,j}}\right]^{-1},
	$$
	where $\boldsymbol{\rho}_m^{(\boldsymbol{z}, \boldsymbol{p})}=\textbf{V}_m^* \boldsymbol{h}_m^{(\boldsymbol{z}, \boldsymbol{p})}$.
	
	\textbf{Limited aperture case:}
	Denote by $[a,b]\subset[0,2\pi]$ a non-empty interval.
	Let
	$$\mathbb{D}_{la}=\left\{\boldsymbol{d}_j=\left(\cos \theta_j, \sin \theta_j\right)^{\top}:\ \theta_j=(b-a)j / N, j=1,2, \ldots, N\right\}$$
	be a set of $N$ equidistantly distributed directions.
	We sequentially use a plane compression wave and a plane shear waves as incident wave in each direction $\boldsymbol{d}_j, j=1,2, \ldots, N$, i.e., a total of $2N$ incident waves. For the illumination in direction $\boldsymbol{d}_j\in \mathbb{D}_{la}$ , we collect far field data
	\begin{equation}
		\mathbf{F}_{k, j}=\left(\begin{array}{cc}
			u_p^{\infty}\left(\boldsymbol{d}_k, \boldsymbol{d}_j ; p\right) & u_p^{\infty}\left(\boldsymbol{d}_k, \boldsymbol{d}_j ; s\right) \\
			u_s^{\infty}\left(\boldsymbol{d}_k, \boldsymbol{d}_j ; p\right) & u_s^{\infty}\left(\boldsymbol{d}_k, \boldsymbol{d}_j ; s\right)
		\end{array}\right)
	\end{equation}
	in each direction $\boldsymbol{d}_k \in \mathbb{D}_{la}$.
	
	From the definition of $\bm{\mathcal{F}}_{la}$, we have
	\begin{equation}\label{Fmgm}
		\bm{\mathcal{F}}_{la}\bm{g}(\hat{\bm{x}},\bm{z};\bm{p})
		=\boldsymbol{\Gamma}
		^{\infty}(\hat{\bm{x}}, \boldsymbol{z} ; \boldsymbol{p}),\quad \hat{\bm{x}}\in \mathbb{U}.
	\end{equation}
	We use the trapezoidal rule to discrete equation (\ref{Fmgm}) and rewrite the equation in matrix form as
	\begin{equation}
		\textbf{F}_{la}\textbf{g}_{la}^{(\bm{z},\bm{p})}=
		\bm{h}_{la}^{(\bm{z},\bm{p})},
	\end{equation}
	where $\textbf{F}_{la}=\textbf{F}|_{\mathbb{D}_{la}}$,
	$\textbf{g}_{la}^{(\bm{z},\bm{p})}
	=\textbf{g}^{(\bm{z},\bm{p})}|_{\mathbb{D}_{la}}$,
	and $\boldsymbol{h}_{la}^{(\boldsymbol{z}, \boldsymbol{p})}=\boldsymbol{h}^{(\boldsymbol{z}, \boldsymbol{p})}|_{\mathbb{D}_{la}}$.
	We compute a singular value decomposition $\mathbf{F}_{la}=\textbf{U}_{la}\textbf{S}_{la} \textbf{V}_{la}^*$, where the diagonal elements of $\mathbf{S}_{la}$ are denoted by $\eta_{la,j}$, $j=1,2, \ldots, 2N$. Since the operator $\bm{\mathcal{F}}_{la}$ fails to be normal, we take $\lambda_{la,j}=|\mathrm{Re}\eta_{la,j}|
	+|\mathrm{Im}\eta_{la,j}|$.
	For each $\bm{z}\in \mathbb{Z}$, we compute
	$$
	W_{la}(\boldsymbol{z}):=\left[\sum_{j=1}^N \frac{\left|\rho_{la,j}^{(\boldsymbol{z}, \boldsymbol{p})}\right|^2}
	{\lambda_{la,j}}\right]^{-1},
	$$
	where $\boldsymbol{\rho}_{la}^{(\boldsymbol{z}, \boldsymbol{p})}=\textbf{V}_{la}^* \boldsymbol{h}_{la}^{(\boldsymbol{z}, \boldsymbol{p})}$.
	
	In our numerical experiments, we use the integral equation approach to solve the forward problem, and the integration nodes are always taken to be $128$. In the inverse problem, both incident and observation directions are taken to be $64$. The set of sampling points is
	$$\mathbb{Z}=\{(x_j,y_k):\ x_j=-3+\frac{3}{50}j,y_k=-3+\frac{3}{50}k,\ j,k=0,1,
	\ldots,100\}.$$
	In the displayed figures, the actual boundary $\partial D$ is always drawn with the black dashed line.

	\subsection{Influence of polarization direction}
	In this subsection, we examine the effect of the polarization direction $\bm{p}$ defined parametrically by
	$(\cos\alpha,\sin\alpha)$, $\alpha\in[0,2\pi]$, on the numerical results.
	
\begin{figure}[htbp]
		\centering
		\subfigure[]{
			\includegraphics[width=0.25\textwidth]{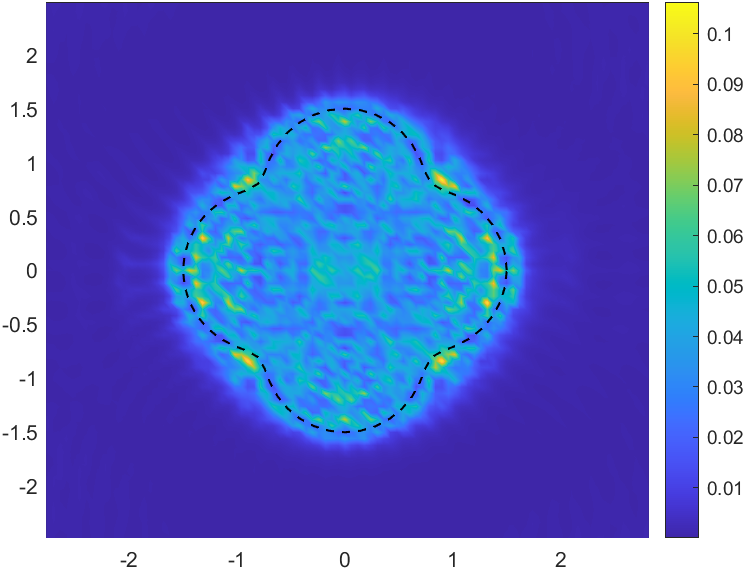}
		}
		\subfigure[]{
			\includegraphics[width=0.25\textwidth]{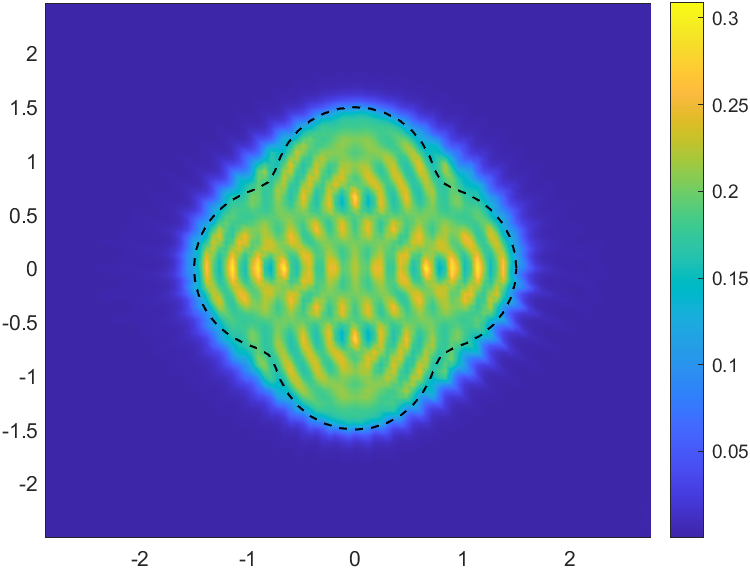}
		}
		\subfigure[]{
			\includegraphics[width=0.25\textwidth]{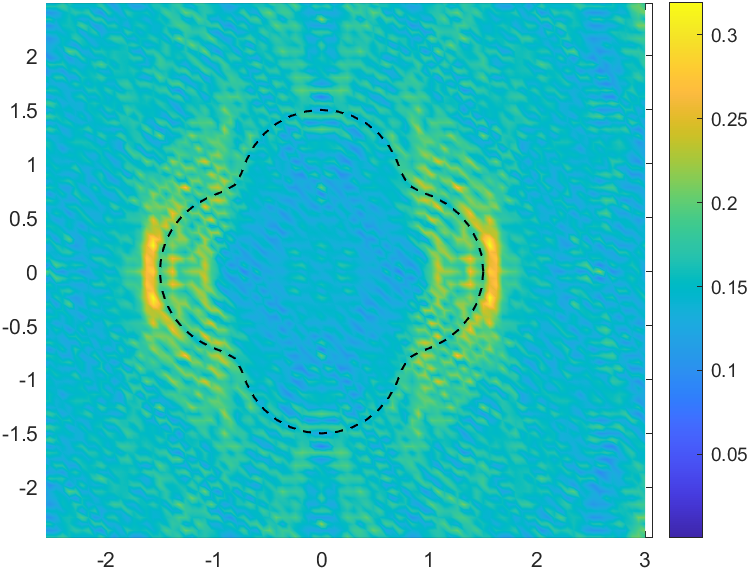}
		}
		
		\subfigure[]{
			\includegraphics[width=0.25\textwidth]{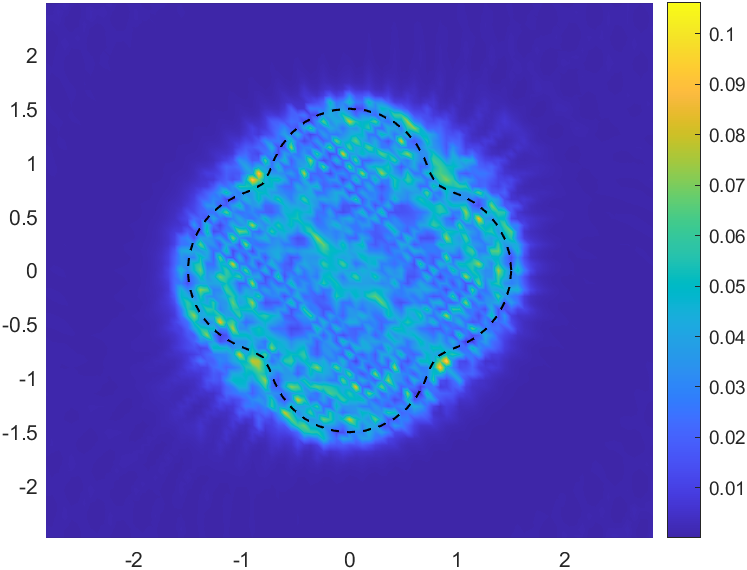}
		}
		\subfigure[]{
			\includegraphics[width=0.25\textwidth]{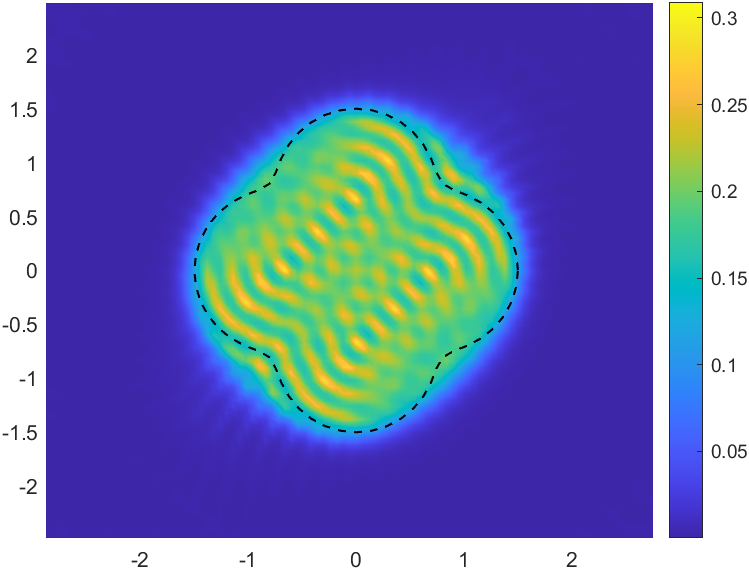}
		}
		\subfigure[]{
			\includegraphics[width=0.25\textwidth]{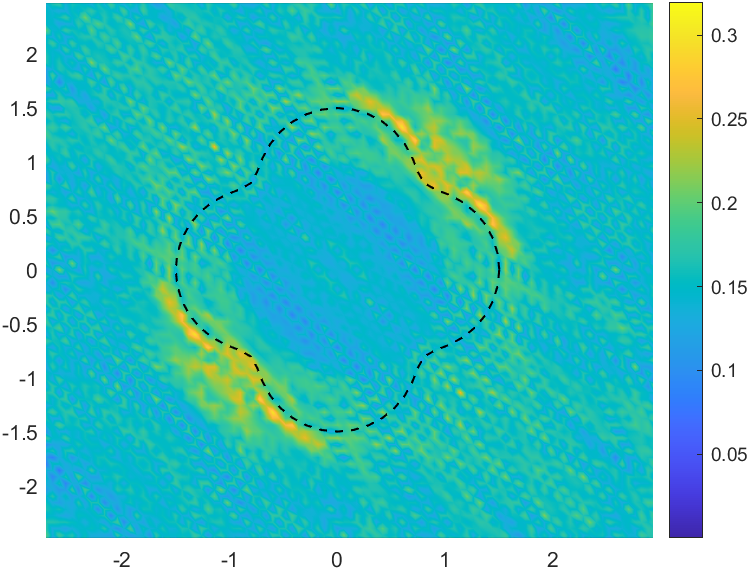}
		}
		
		\subfigure[]{
			\includegraphics[width=0.25\textwidth]{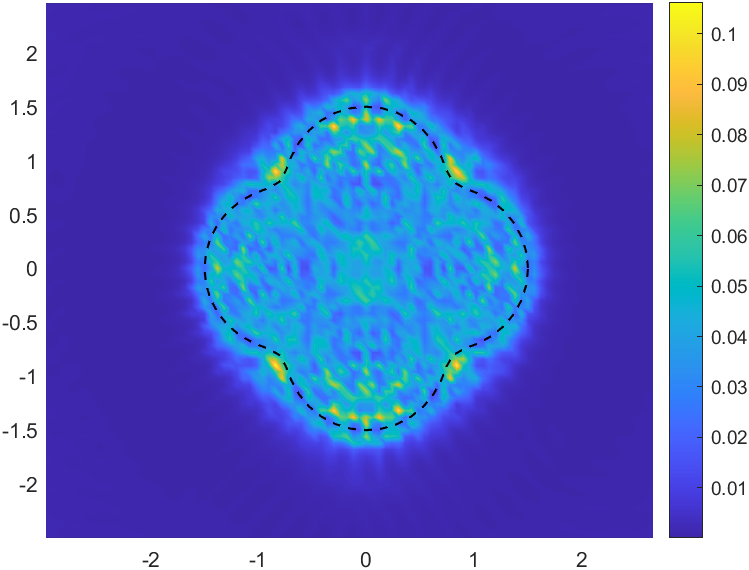}
		}
		\subfigure[]{
			\includegraphics[width=0.25\textwidth]{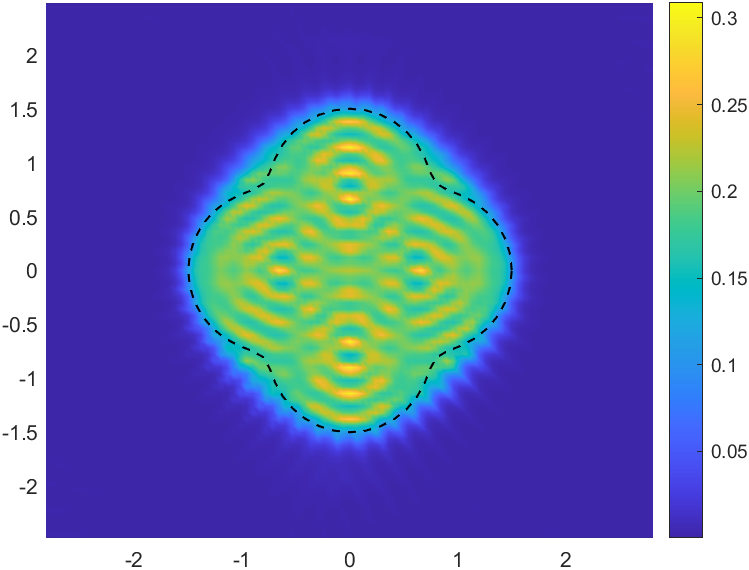}
		}
		\subfigure[]{
			\includegraphics[width=0.25\textwidth]{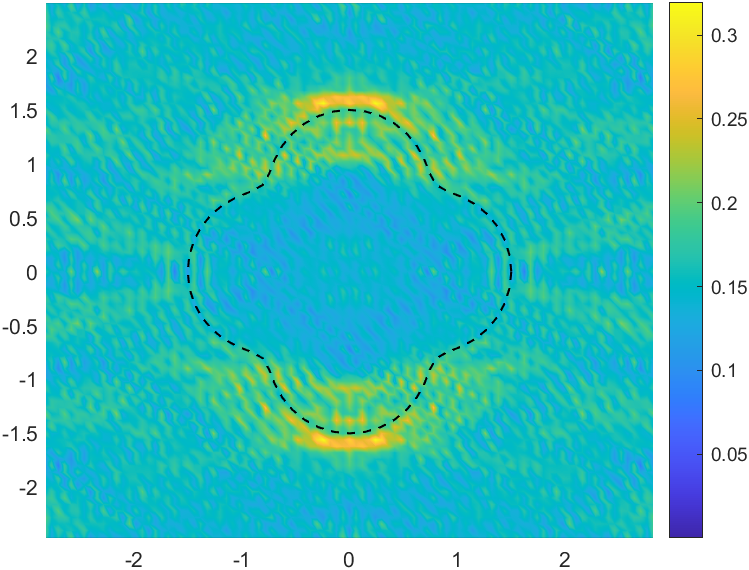}
		}
		
		\subfigure[]{
			\includegraphics[width=0.25\textwidth]{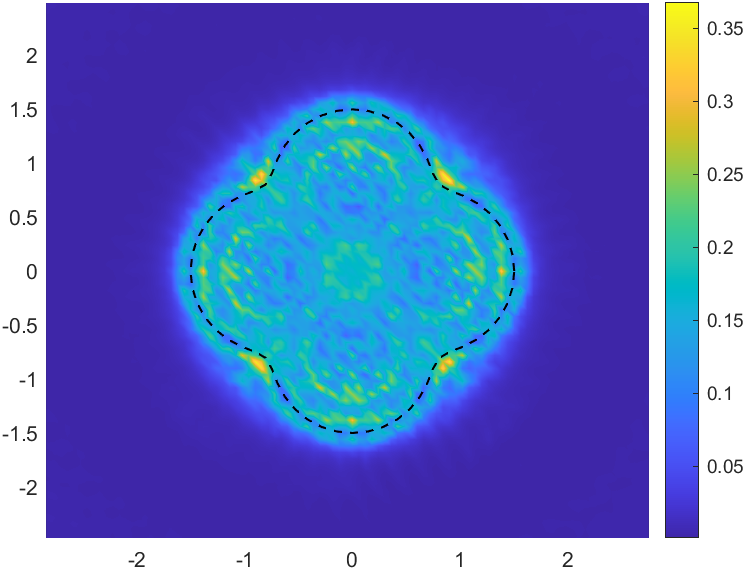}
		}
		\subfigure[]{
			\includegraphics[width=0.25\textwidth]{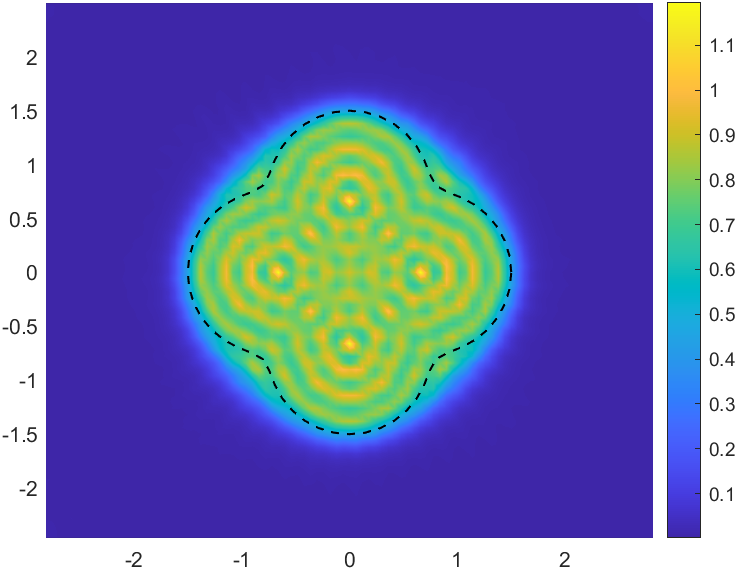}
		}
		\subfigure[]{
			\includegraphics[width=0.25\textwidth]{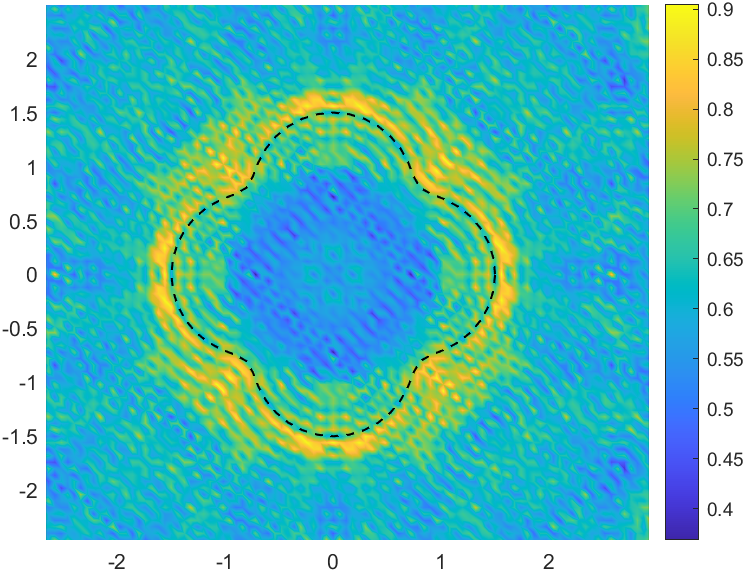}
		}
		\caption{Reconstructions of rounded rectangle shaped obstacle for different polarization directions with $\delta=10\%$.
			Top row: $\alpha=0$; Second row: $\alpha=\pi/4$; Third row: $\alpha=\pi/2$; Bottom row: Superposition of four polarization directions: $\alpha\in \{0,\pi/4,\pi/2,3\pi/4\}$. Left column: FF case; Middle column: PP case; Right column: SS case.}
		\label{F_1}
	\end{figure}

	We consider a rounded rectangle-shaped cavity parameterized as
	\begin{equation}\label{5.1.1}
		\bm{x}(t)=1.5\sqrt[10]{\cos ^{10} t+\sin ^{10} t}\ (\cos t, \sin t),\quad t\in[0,2\pi].
	\end{equation}
	The reconstructions with four different polarization angles are shown in Figure \ref{F_1}, where we have taken $\mu=1,\lambda=1,\omega=8\pi$.
	In Figure \ref{F_1}, columns $1$, $2$, and $3$ correspond to the FF, PP, and SS cases, respectively, and the first three rows relate to the polarization angle $\alpha=0$, $\pi/4$, $\pi/2$, respectively.
	The results in the fourth row of Figure \ref{F_1} are obtained by inverting the contours of the cavity using the four polarization angles $\alpha= 0$, $\pi/4$, $\pi/2$, and $3\pi/4$ in turn, and superimposing their results. From Figure \ref{F_1} $(\mathrm{a})$-$(\mathrm{i})$, it can be seen that for all three cases, FF, PP, and SS, sharper contours are obtained along the polarization direction and in its opposite direction than in the other directions.
	Such phenomenon is most evident for the SS case. This observation motivates us to explore a novel indicator function based on the superposition of multiple polarization direction indicator functions, the results of which are shown in Figure \ref{F_1} $(\mathrm{j})$-$(\mathrm{l})$. The novel indicator function significantly improve the quality of the reconstructed profiles compared with the original indicator function.
	
	Consequently, in the subsequent numerical examples, a superposition of four polarization directions $\alpha\in \{0,\pi/4,\pi/2,3\pi/4\}$ will be used as the default polarization direction unless otherwise specified.
	
	\subsection{Influence of frequency and wavenumber}
\begin{figure}[htbp]
		\centering
		\subfigure[]{
			\includegraphics[width=0.25\textwidth]{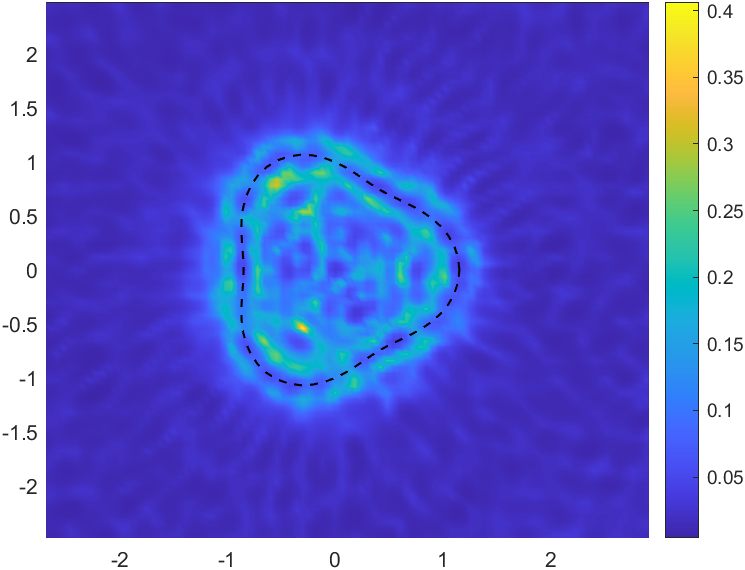}
		}
		\subfigure[]{
			\includegraphics[width=0.25\textwidth]{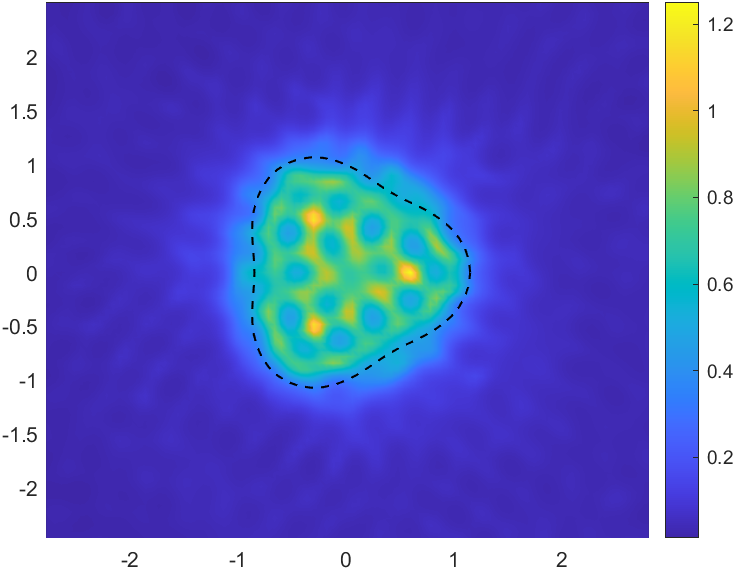}
		}
		\subfigure[]{
			\includegraphics[width=0.25\textwidth]{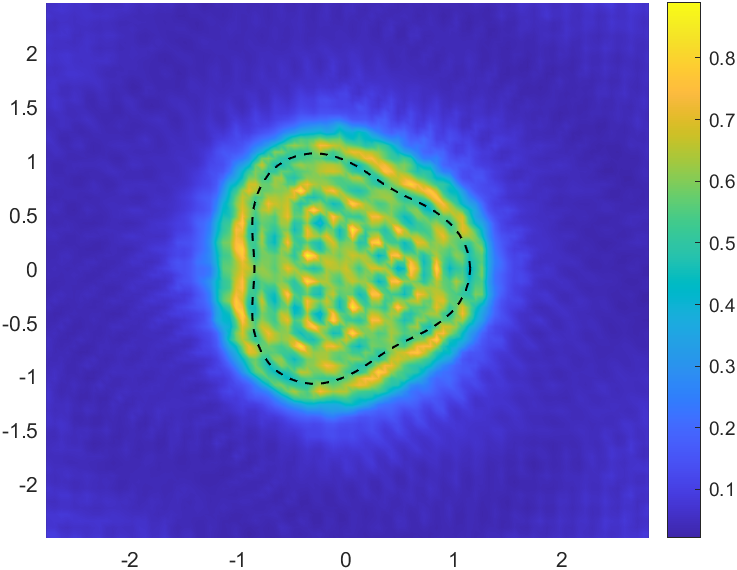}
		}
		
		\subfigure[]{
			\includegraphics[width=0.25\textwidth]{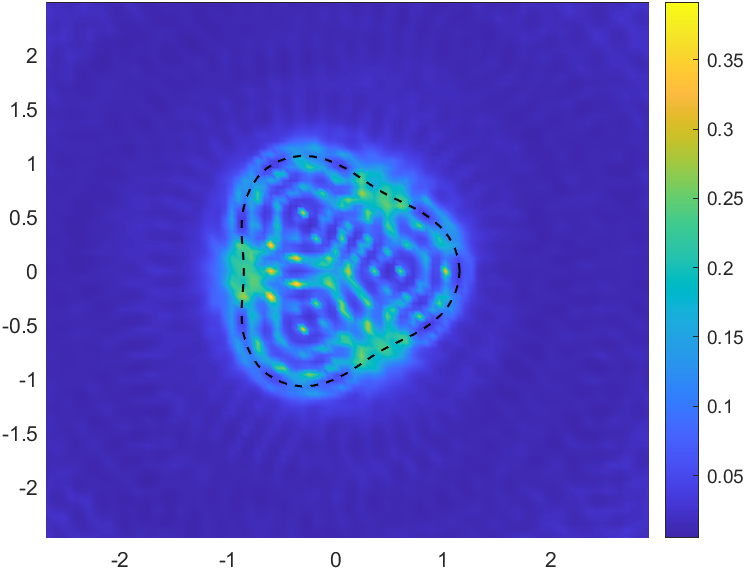}
		}
		\subfigure[]{
			\includegraphics[width=0.25\textwidth]{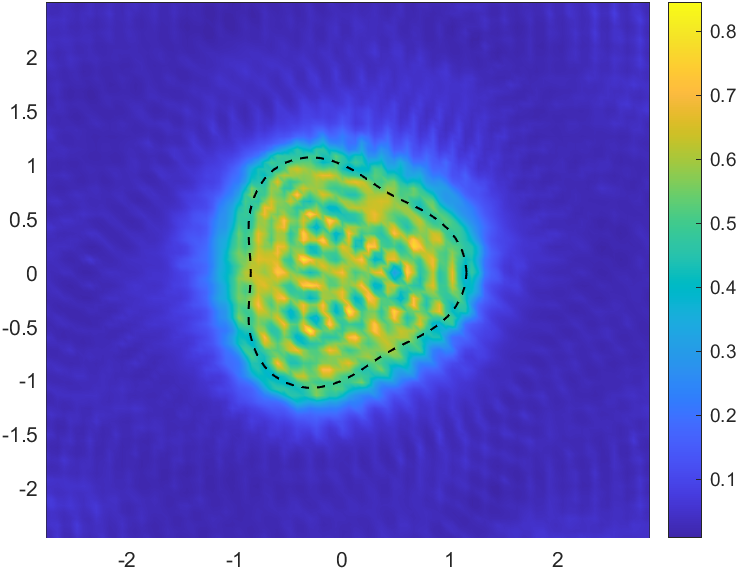}
		}
		\subfigure[]{
			\includegraphics[width=0.25\textwidth]{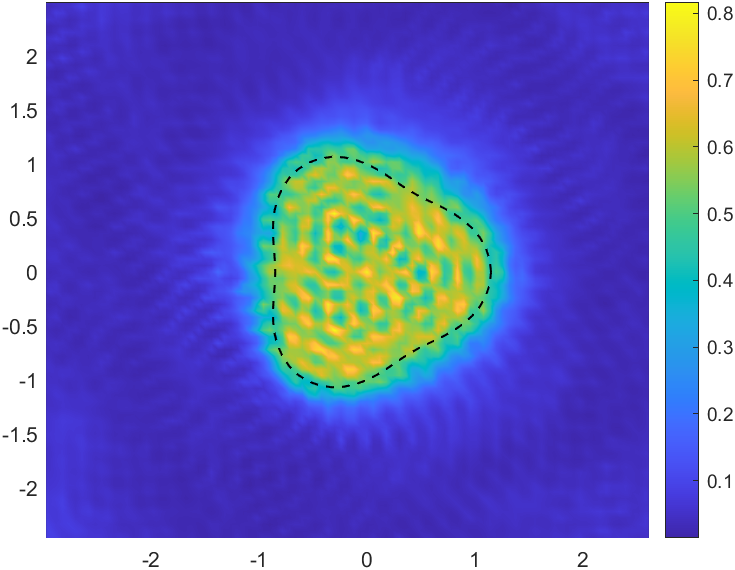}
		}
		
		\subfigure[]{
			\includegraphics[width=0.25\textwidth]{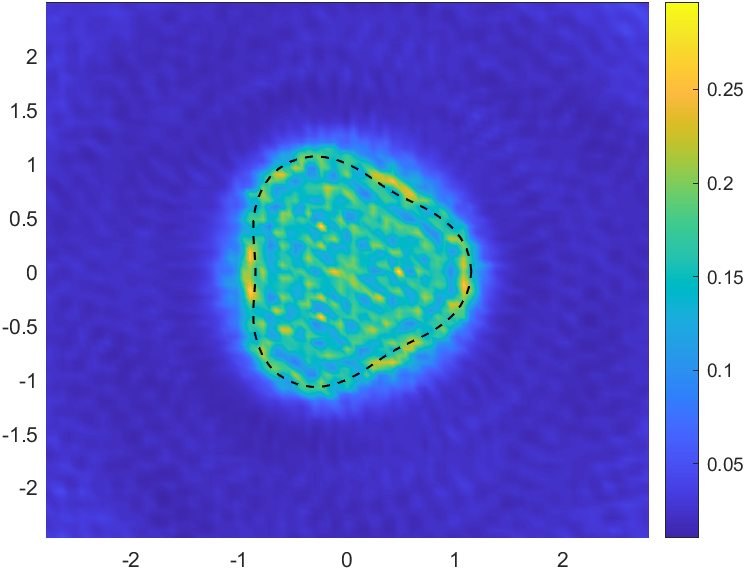}
		}
		\subfigure[]{
			\includegraphics[width=0.25\textwidth]{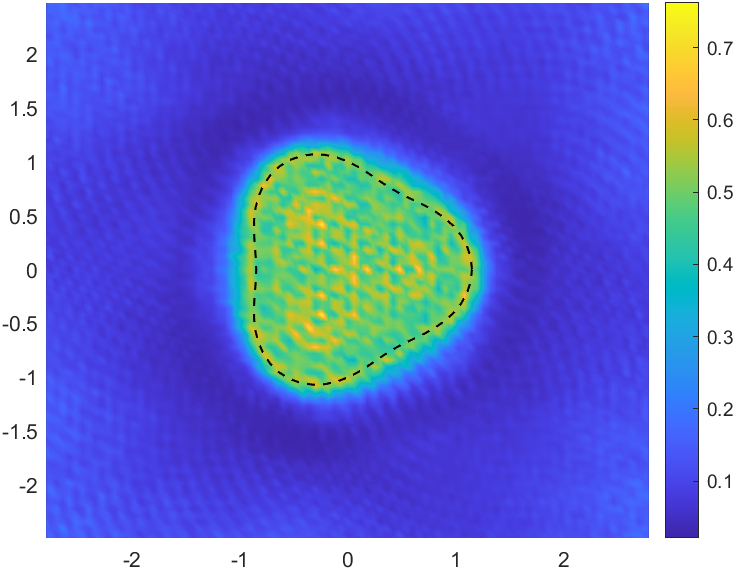}
		}
		\subfigure[]{
			\includegraphics[width=0.25\textwidth]{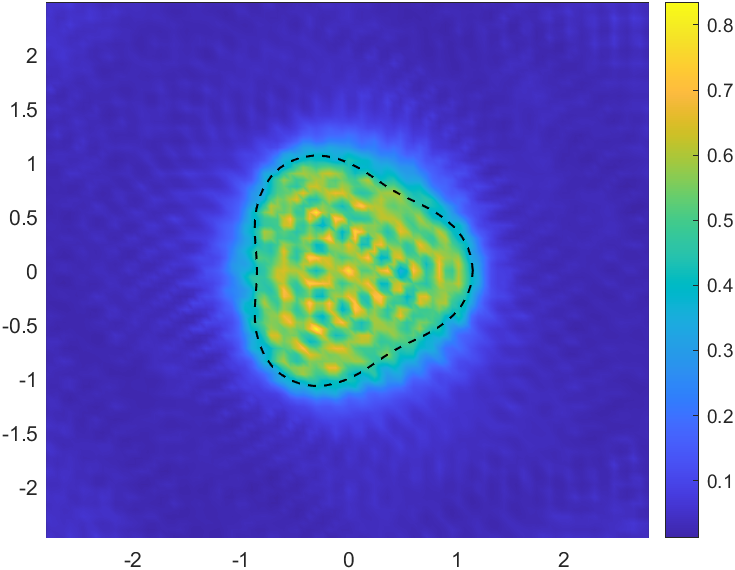}
		}
		\caption{Reconstructions of  pear-shaped obstacle for different $\lambda$ with $\delta=10\%$.
			Top row: $\lambda=1$; Second row: $\lambda=-1$; Third row: $\lambda=-1.5$. Left column: FF case; Middle column: PP case; Right column: SS case.}
		\label{F_pear}
	\end{figure}

	Consider a pear-shaped cavity parameterized as
	\begin{equation*}
		x(t)=(1+0.15\cos(3t))\ (\cos t, \sin t),\quad t\in[0,2\pi].
	\end{equation*}
	Choose $\omega=5\pi$ and $\mu=1$. The reconstructions are shown in Figure \ref{F_pear}. Columns 1, 2, and 3 correspond to the FF, PP, and SS cases, respectively, and rows 1, 2, and 3 relate to $\lambda=1$, $-1$, $-1.5$, respectively. It looks like the FF case is always the worst. When $\lambda=1$, $k_p<k_s$, it looks like the SS case is the best. When $\lambda=-1.5$, $k_p>k_s$, it looks like the PP case is the best. In addition, when $\lambda=-1$, $k_p=k_s$, the reconstruction results are closer for the PP and SS cases.
	
	%
	%
	
\begin{figure}[htbp]
		\centering
		\subfigure[]{
			\includegraphics[width=0.25\textwidth]{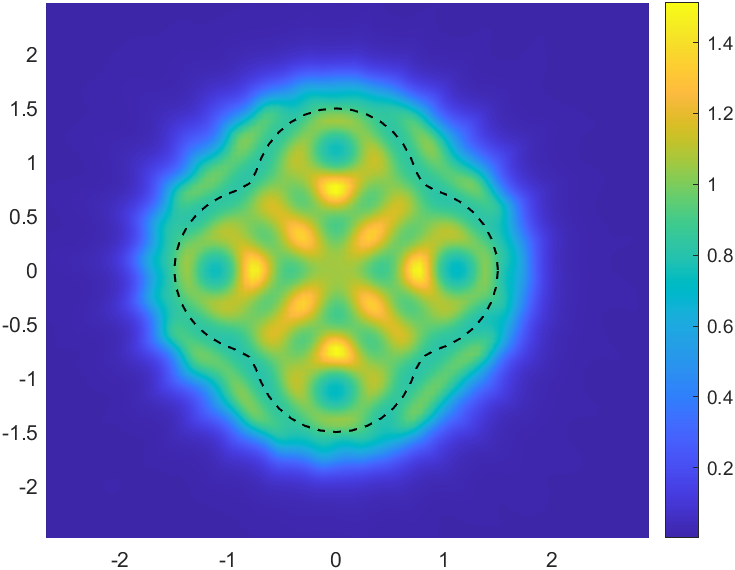}
		}
		\subfigure[]{
			\includegraphics[width=0.25\textwidth]{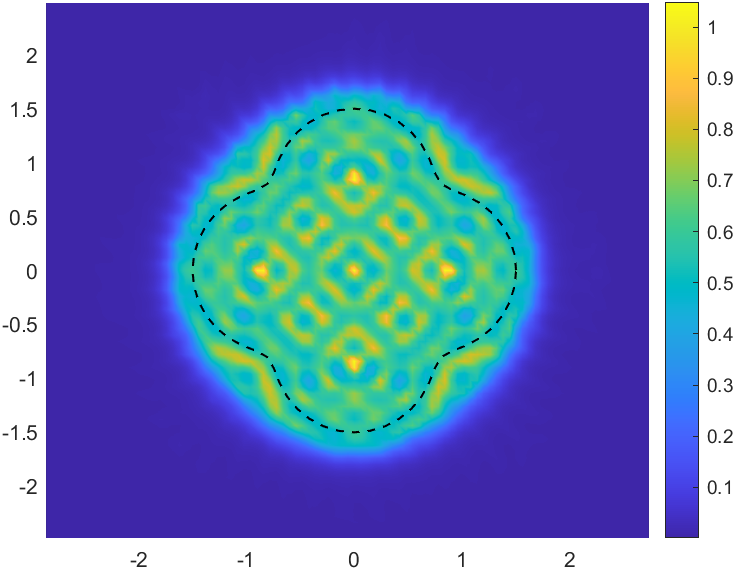}
		}
		\subfigure[]{
			\includegraphics[width=0.25\textwidth]{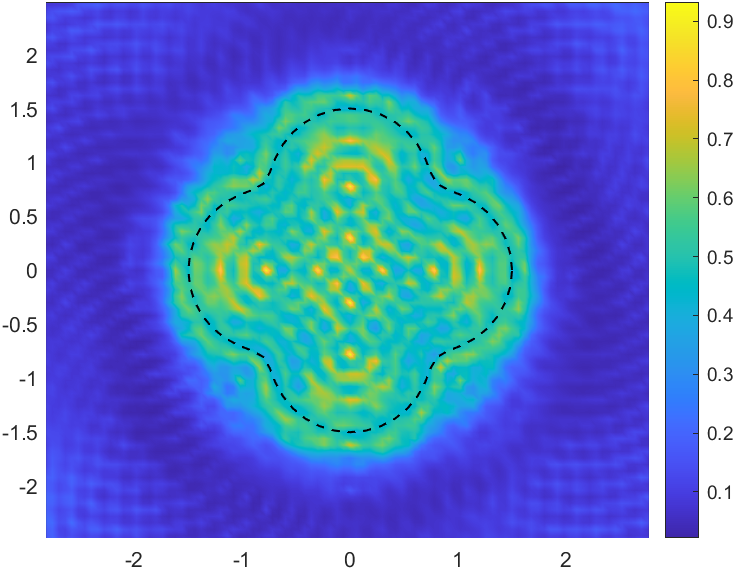}
		}
		
		\subfigure[]{
			\includegraphics[width=0.25\textwidth]{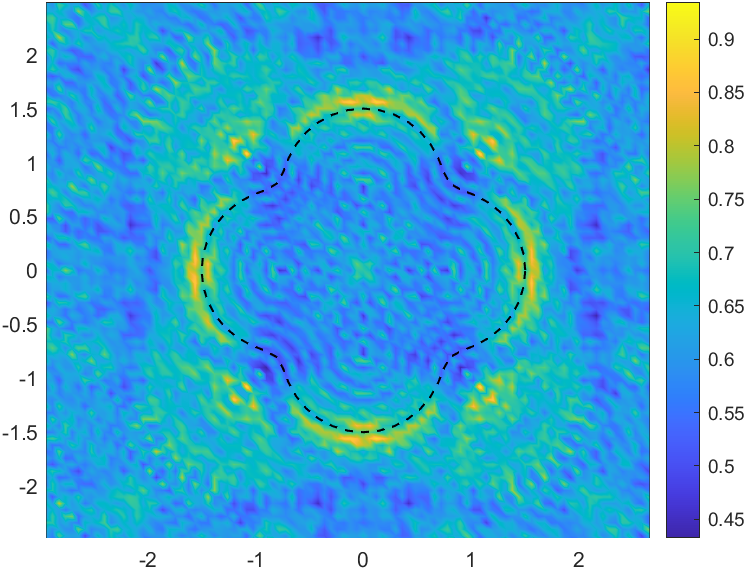}
		}
		\subfigure[]{
			\includegraphics[width=0.25\textwidth]{F0pi_Sp.png}
		}
		\subfigure[]{
			\includegraphics[width=0.25\textwidth]{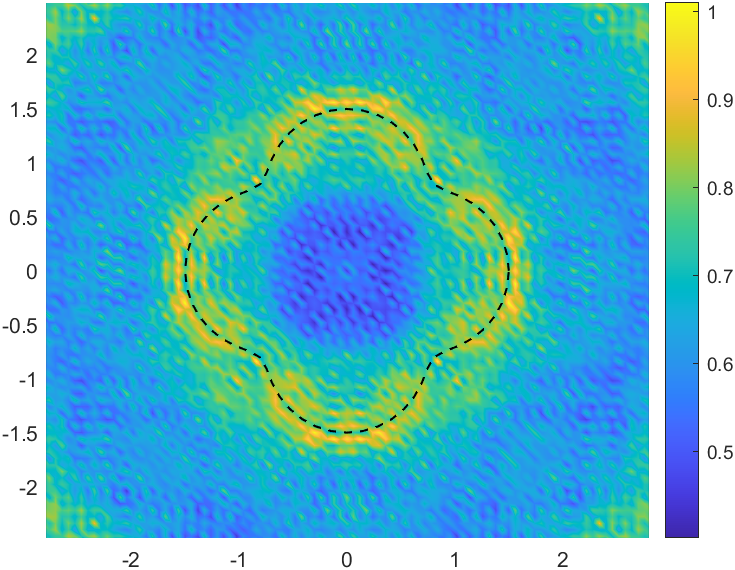}
		}
		\caption{Reconstructions of rounded rectangle-shaped cavity for different frequencies $\omega$ with  $\delta=10\%$.
			$\omega=2\pi, 4\pi, 5\pi, 7\pi,
			8\pi, 9\pi$ for
			$(\mathrm{a})$ to $(\mathrm{f})$, respectively.}
		\label{F_Round}
	\end{figure}

	We consider the influence of frequency $\omega$ on reconstruction. The reconstruction results for $\omega=2\pi, 4\pi, 5\pi, 7\pi,
	8\pi, 9\pi$, with $\lambda=1$ and $\mu=1$ as parameters, are
	shown in Figure \ref{F_Round}. Based on the relationship between wavenumbers and reconstruction results discussed earlier, we only consider SS case here. It can be observed from these results that the rough outline of the scatterer can be obtained at low frequencies and more localized details of the boundary of the scatterer can be captured at high frequencies.
	
	Neglecting the oscillations generated at higher frequencies, it seems that the larger the frequency is, the better the detail of the reconstructed contour will be. To gather different features of low and high frequencies simultaneously, it is often recommended to use multi-frequency data to enhance better reconstruction quality.
	
	\subsection{Tests on multiple scatterers}
\begin{figure}[htbp]
		\centering
		\subfigure[]{
			\includegraphics[width=0.25\textwidth]{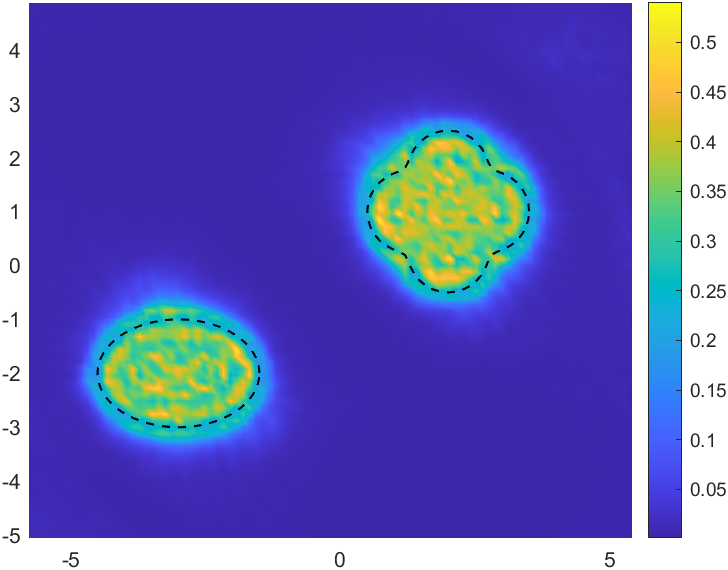}
		}
		\subfigure[]{
			\includegraphics[width=0.25\textwidth]{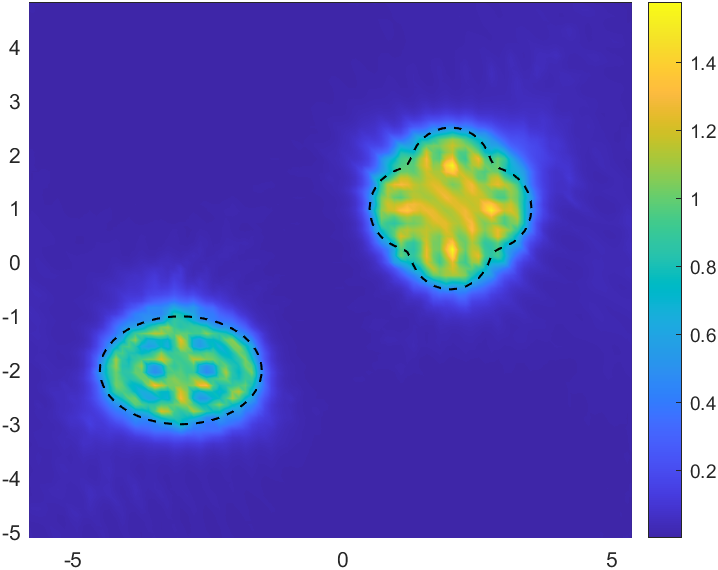}
		}
		\subfigure[]{
			\includegraphics[width=0.25\textwidth]{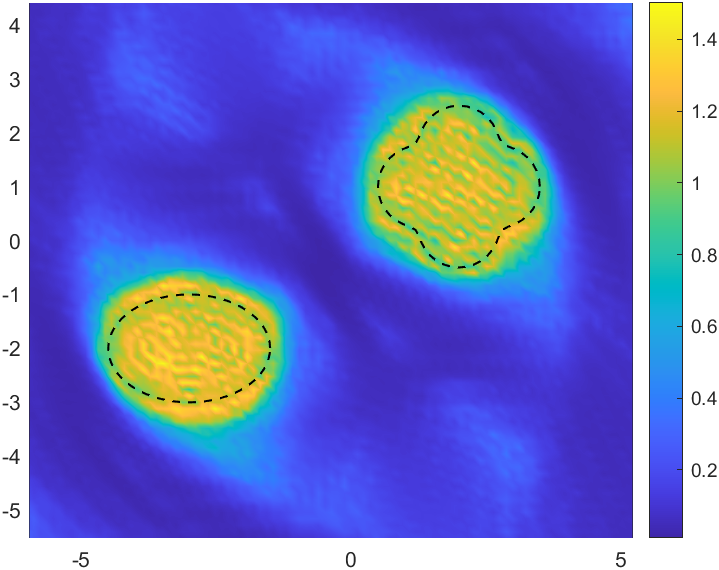}
		}
		
		\subfigure[]{
			\includegraphics[width=0.25\textwidth]{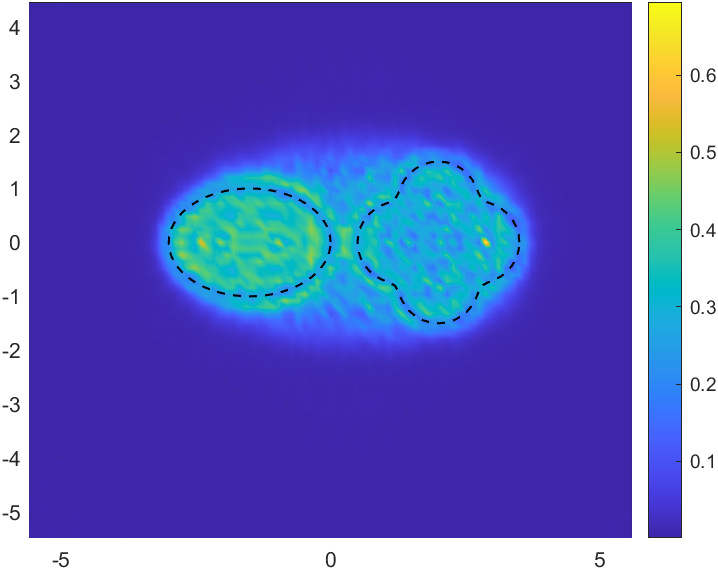}
		}
		\subfigure[]{
			\includegraphics[width=0.25\textwidth]{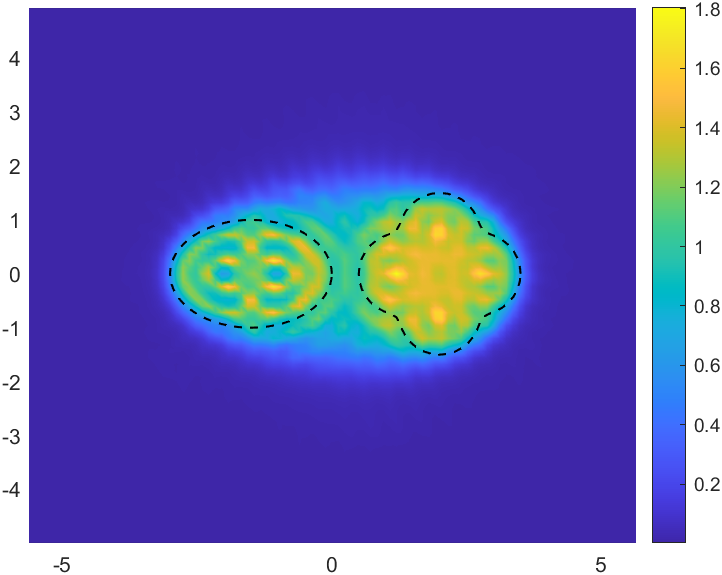}
		}
		\subfigure[]{
			\includegraphics[width=0.25\textwidth]{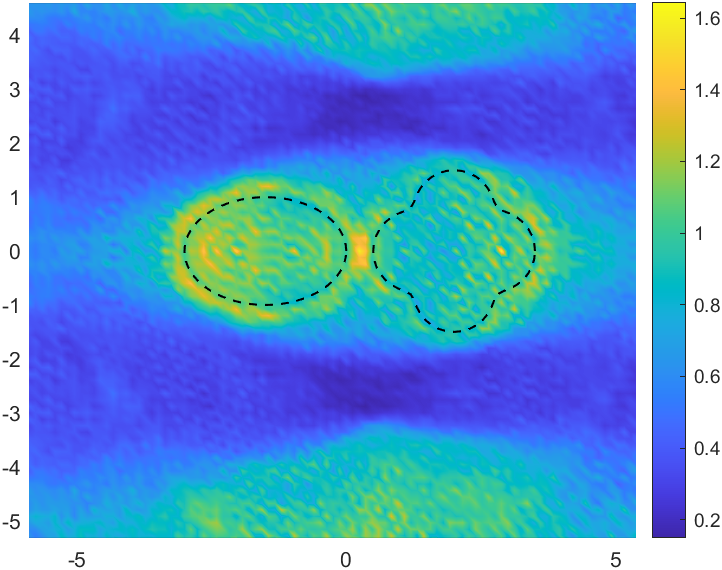}
		}
		\caption{Reconstructions of two scatterers for $\omega=9\pi, \mu=5,\lambda=5$ with $\delta=10\%$. Left column: FF case; Middle column: PP case; Right column: SS case.}
		\label{F_Two}
	\end{figure}

	In this subsection, denote by $D=D_1\cup D_2$ the union of two disjoint cavities.
	The set of sampling points is
	$$\mathbb{Z}=\{(x_j,y_k):\ x_j=-6+\frac{3}{25}j,y_k=-6+\frac{3}{25}k,\ j,k=0,1,
	\ldots,100\}.$$
	$D_1$ is chosen to be an ellipse  with axes 1.5 and 1, respectively, and $D_2$ is a rounded rectangular cavity with parameter (\ref{5.1.1}). When $D_1$ and $D_2$ are far apart (the centers of $D_1$ and $D_2$ are taken at $(-3,2)$ and $(2,1)$, respectively), the reconstruction results are shown in Figure \ref{F_Two} $(a)$-$(c)$.
	For the FF, PP, and SS cases, it is possible to get a rough outline of $D_1$ and $D_2$ and to separate them clearly. In particular, it seems to look best in the FF case.
	
	When $D_1$ and $D_2$ are close to each other (the centers of $D_1$ and $D_2$ are taken at $(-1.5,1)$ and $(1.5,-1)$, respectively), the reconstruction results are displayed in Figure \ref{F_Two}
	$(d)$-$(f)$. For the FF, PP and SS cases, the contours of D1 and D2 can be obtained roughly, but the reconstruction results are worse than the above cases, especially the FF case seems to be the worst. We find that none of the three cases can clearly separate $D_1$ and $D_2$ when the two scatterers are too close together. It is worth pointing out that for the SS case, the imaging results oscillate to some extent regardless of whether the two scatterers are farther or closer apart.
	
	%
	
	\subsection{Influence of multi-scale scatterers}
	In this subsection, suppose that $D$ consists of two scatterers $D_1$ and $D_2$ at different scales.
	We take
	$$\mathbb{Z}=\{(x_j,y_k):\ x_j=-8+\frac{4}{25}j,y_k=-8+\frac{4}{25}k,\ j,k=0,1,
	\ldots,100\}$$
	as the set of sampling points.
	Two examples are provided to demonstrate the influence
	of multi-scale scatterers on reconstruction results.
	
	\textbf{Example 1:} Let $D_1$ be a rounded rectangle-shaped cavity centered at (-2,0) whose diameter is three times the diameter of the original shape (\ref{5.1.1}), and $D_2$ be a pear-shaped cavity centered at $(6,0)$ which is parameterized by
	\begin{equation}
		\bm{x}(t)=0.5(1+0.2\cos(3t))(\cos t, \sin t),\quad t\in[0,2\pi].
	\end{equation}
	The size of $D_1$ is about $8$ times the size of $D_2$.
	
\begin{figure}[htbp]
		\centering
		\subfigure[]{
			\includegraphics[width=0.25\textwidth]{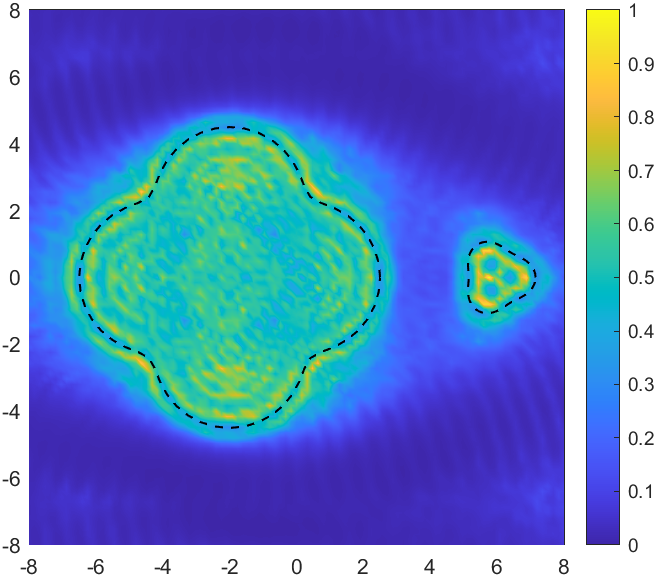}
		}
		\subfigure[]{
			\includegraphics[width=0.25\textwidth]{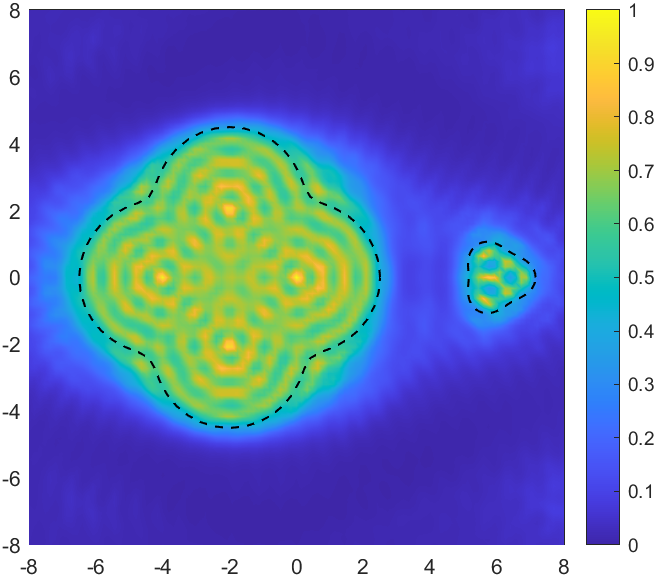}
		}
		\subfigure[]{
			\includegraphics[width=0.25\textwidth]{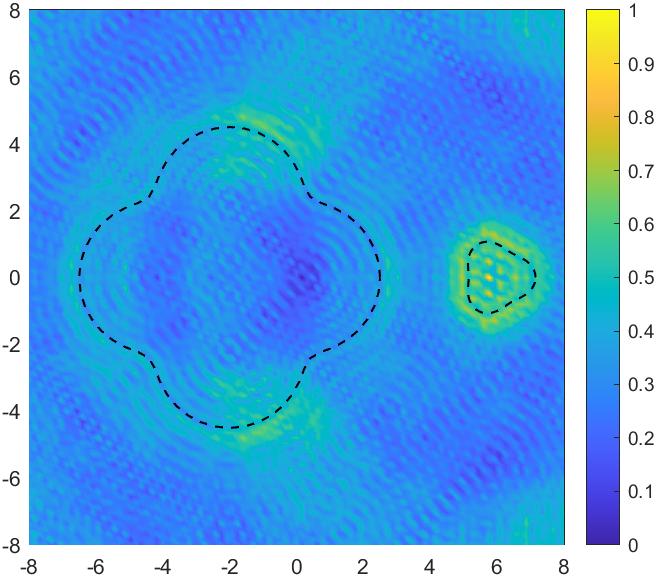}
		}
		\caption{Reconstructions of two scatterers for $\omega=6\pi,\mu=5,\lambda=5$ with $\delta=10\%$. Left column: FF case; Middle column: PP case; Right column: SS case.}
		\label{F_three}
	\end{figure}

	The reconstruction results are shown in Figure \ref{F_three}.
	We can observe that both the FF and PP cases can effectively reconstruct these two different sized scatterers. In particular, the large-sized cavity seems to be better in the PP case, while the small-sized cavity seems to look better in the FF case. However, in the SS case, neither large nor small-sized cavities are satisfactorily reconstructed, and the imaging results produce some oscillations.
	
\begin{figure}[htbp]
		\centering
		\subfigure[]{
			\includegraphics[width=0.25\textwidth]{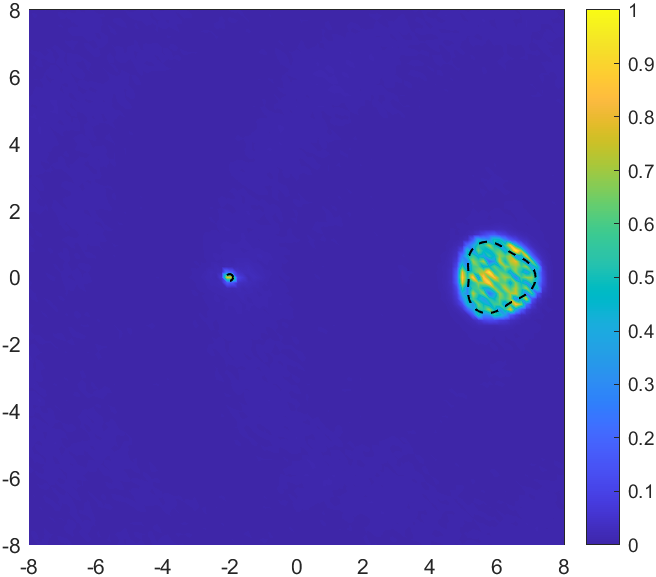}
		}
		\subfigure[]{
			\includegraphics[width=0.25\textwidth]{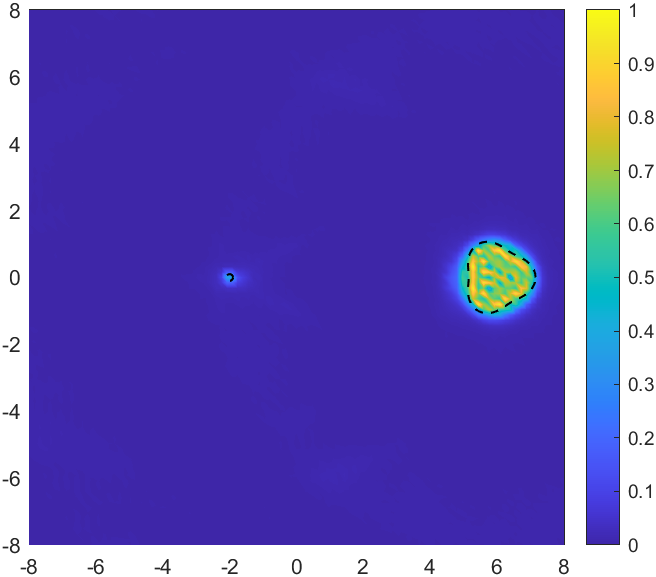}
		}
		\subfigure[]{
			\includegraphics[width=0.25\textwidth]{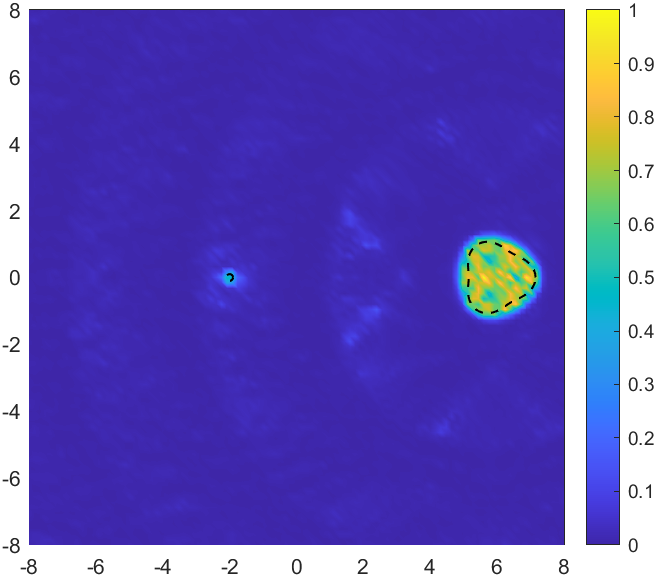}
		}
		\caption{Reconstructions of two scatterers with $\omega=12\pi,\mu=5,\lambda=5$ and $\delta=10\%$. Left column: FF case; Middle column: PP case; Right column: SS case.}
		\label{F_Five}
	\end{figure}

	\textbf{Example 2:} We replace $D_1$ in Example 1 with a circular cavity of radius $0.1$ centered at $(-2,0)$, and keep $D_2$ in the same shape and position as in Example 1. Now $D_1$ is one-fifth the size of $D_2$.
	From the reconstruction results in Figure \ref{F_Five}, it can be seen that not only regular-sized cavities can be reconstructed in the FF, PP, and SS cases, but also tiny cavities can be reconstructed at the same time. However, in the SS case, the imaging results still have
	some small oscillations.

	\subsection{Tests on limited aperture}
\begin{figure}[htbp]
		\centering
		\subfigure[$(0,\pi/2)$]{
			\includegraphics[width=0.2\textwidth]{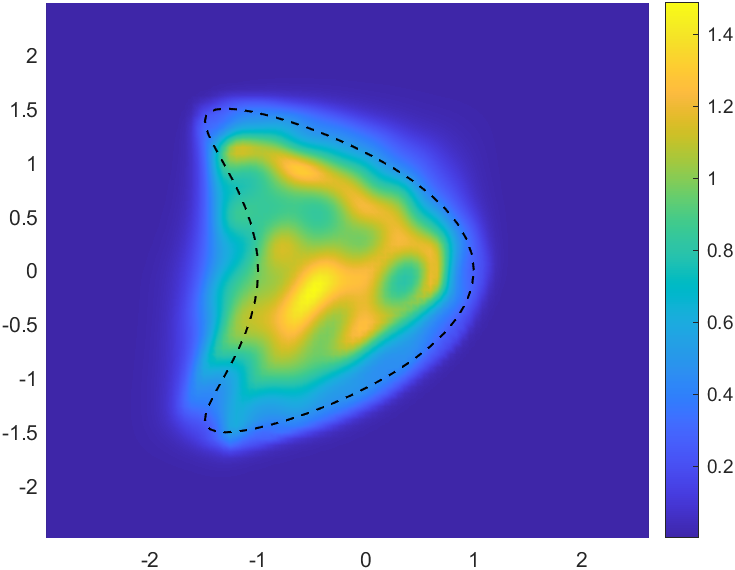}
		}
		\subfigure[$(\pi/2,\pi)$]{
			\includegraphics[width=0.2\textwidth]{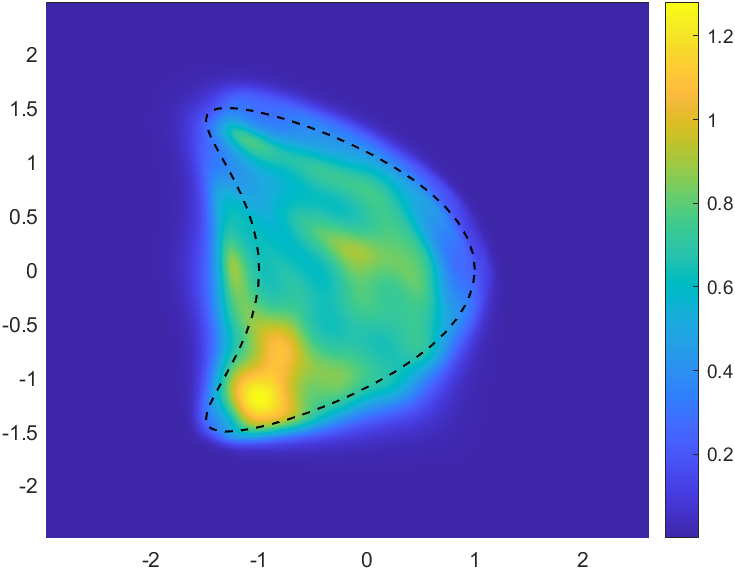}
		}
		\subfigure[$(\pi,3\pi/2)$]{
			\includegraphics[width=0.2\textwidth]{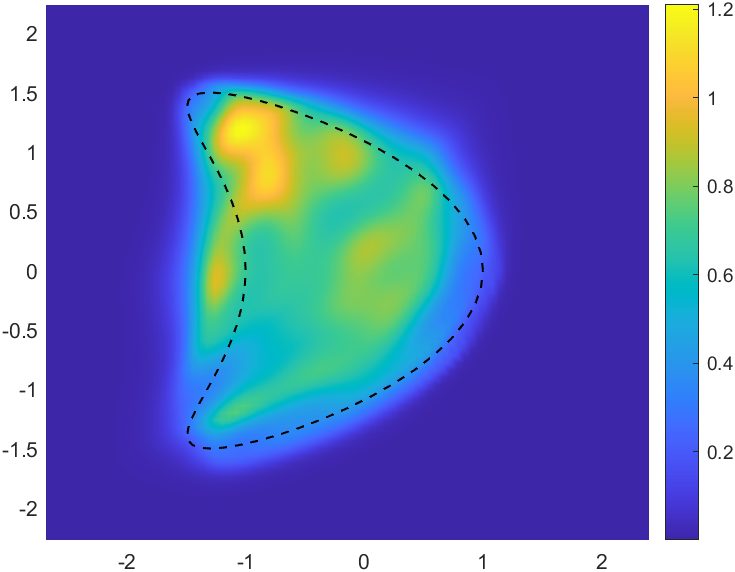}
		}
		\subfigure[$(3\pi/2,2\pi)$]{
			\includegraphics[width=0.2\textwidth]{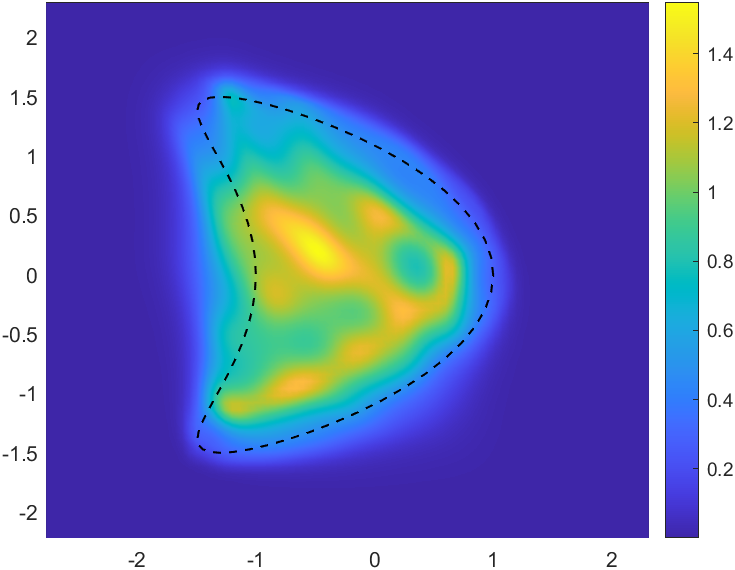}
		}
		
		\subfigure[$(0,\pi/2)$]{
			\includegraphics[width=0.2\textwidth]{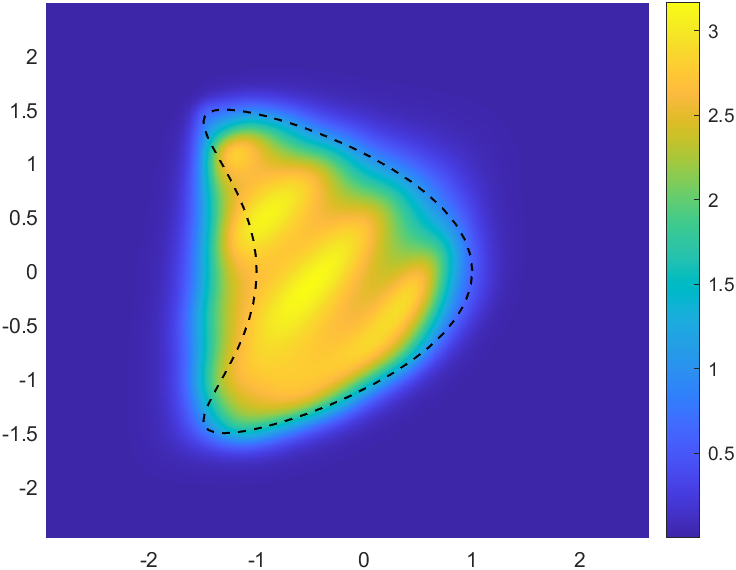}
		}
		\subfigure[$(\pi/2,\pi)$]{
			\includegraphics[width=0.2\textwidth]{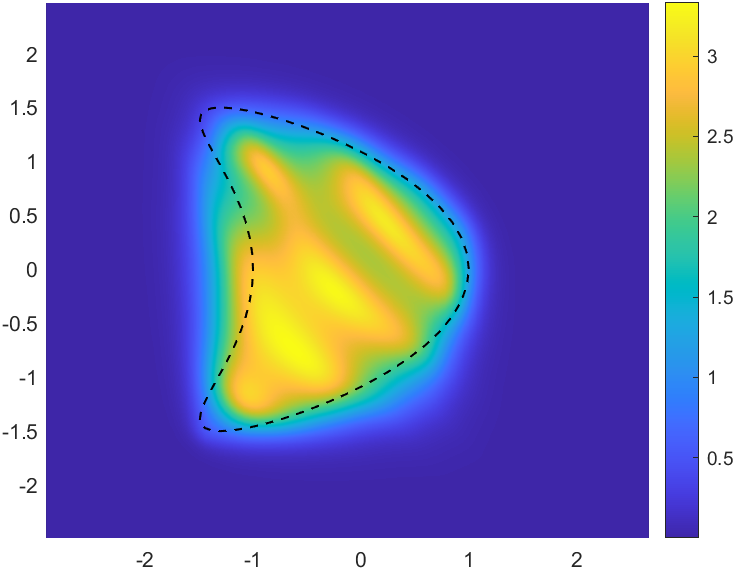}
		}
		\subfigure[$(\pi,3\pi/2)$]{
			\includegraphics[width=0.2\textwidth]{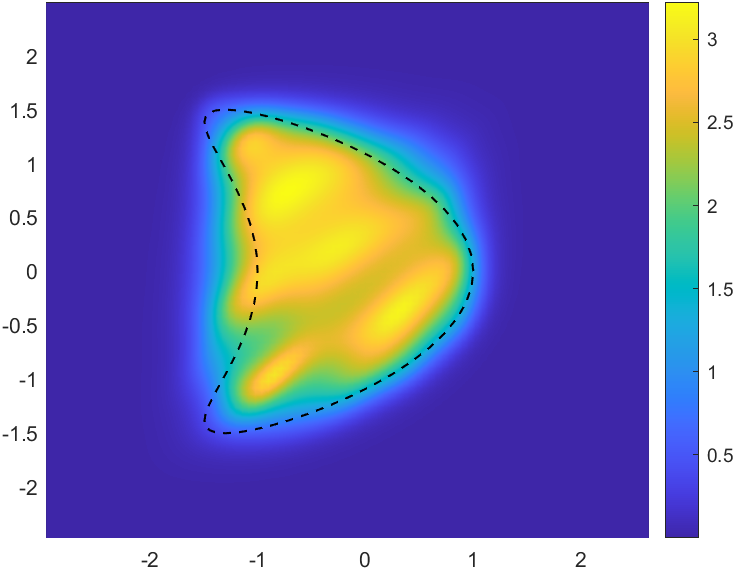}
		}
		\subfigure[$(3\pi/2,2\pi)$]{
			\includegraphics[width=0.2\textwidth]{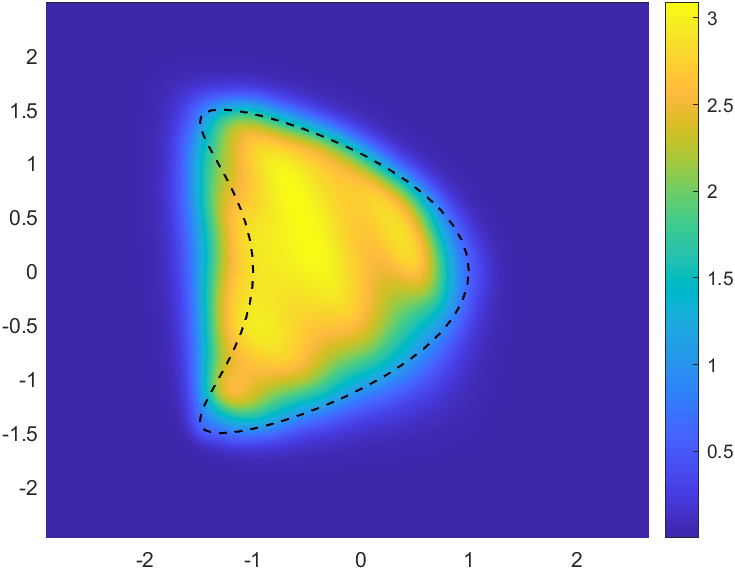}
		}

		\subfigure[$(0,\pi/2)$]{
			\includegraphics[width=0.2\textwidth]{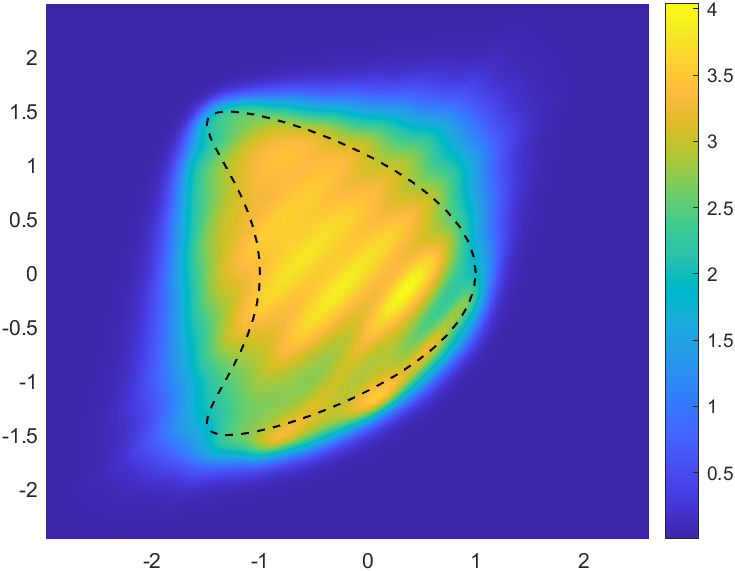}
		}
		\subfigure[$(\pi/2,\pi)$]{
			\includegraphics[width=0.2\textwidth]{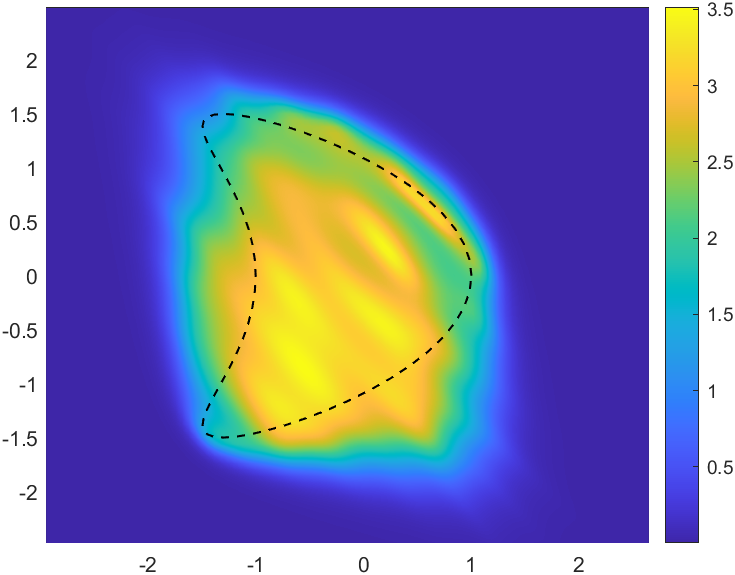}
		}
		\subfigure[$(\pi,3\pi/2)$]{
			\includegraphics[width=0.2\textwidth]{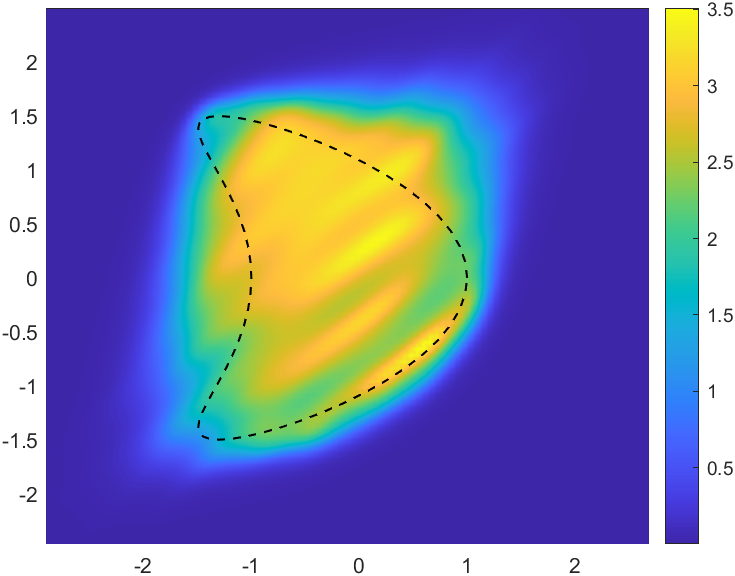}
		}
		\subfigure[$(3\pi/2,2\pi)$]{
			\includegraphics[width=0.2\textwidth]{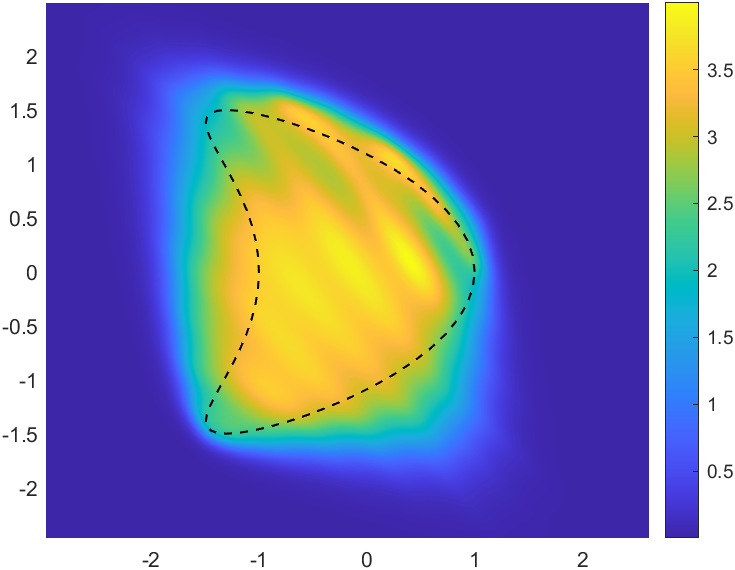}
		}
		\caption{Reconstructions of kite scatterers for
			the limited aperture with $\omega=3\pi, \mu=1, \lambda=2$.  Top row: FF case; Middle row: PP case; Bottom row: SS case.}
		\label{F_L}
	\end{figure}

	In this subsection, we investigate the reconstruction of the cavity in the limited aperture case. We consider a kite-shaped cavity parameterized as
	\begin{equation*}
		\bm{x}(t)=(\cos t+0.65\cos(2t)-0.65,1.5\sin t),\quad
		t\in[0,2\pi],
	\end{equation*}
	and take $\omega=3\pi$, $\mu=1$ and $\lambda=2$.
	The reconstructions are depicted in Figure \ref{F_L}.
	Rows 1-3 correspond to the FF, PP and SS cases,  respectively, and columns 1-4 relate to observation intervals $[a,b]$ as $[0,\pi/2]$, $[\pi/2,\pi]$,
	$[\pi,3\pi/2]$ and $[3\pi/2,2\pi]$, respectively.
	As can be seen, even if we use limited aperture data, the rough contours can still be obtained in FF, PP and SS cases, and in particular, the PP case seems to look
	better. The illuminated part in the PP case is reconstructed better than the unilluminated part.

	\section{Conclusions and future works}\label{Section_6}
	For the inverse scattering problem of elastic waves with Neumann boundary condition, we present a theoretical factorization of the far-field operator and propose a numerical algorithm to reconstruct the shape and location of the elastic cavity. Numerical experiments are performed to demonstrate the effectiveness of our method. We can see that the reconstructions in the SS case seem to look worse in most examples.
	
	It is worth mentioning that our method does not seem to be directly applicable to impedance boundary condition
	\begin{equation*}
		\bm{T}_{\bm{\nu}}\bm{u}+\mathrm{i}\eta\bm{u}=\bm{0},\quad \text{on}~\partial D,
	\end{equation*}
	where $\eta=\omega c$ and $c$ is a positive constant \cite{S05}. Indeed, define by
	$\bm{\mathcal{F}}_{imp}$ the far-field operator for the impedance boundary condition. Using the similar approach in Theorem \ref{F_Fac} yields the factorization of the far-field operator
	$$
	\bm{\mathcal{F}}_{imp }=-\sqrt{8\pi\omega}\bm{\mathcal{G}}_{i m p}\bm{\mathcal{T}}_{i m p}^*\bm{\mathcal{G}}_{i m p}^*,
	$$
	where $\bm{\mathcal{G}}_{imp}: [H^{-1 / 2}(\partial D)]^2 \rightarrow \mathcal{L}^2\left(\mathbb{S}\right)$ maps $\boldsymbol{f} \in\left[H^{-1 / 2}(\partial D)\right]^2$ into the far-field pattern $\boldsymbol{v}^{\infty}=\boldsymbol{\mathcal{G}} \boldsymbol{f}$ of the exterior impedance boundary value problem with boundary data $\boldsymbol{f}$, and $
	\bm{\mathcal{T}}_{imp}: [H^{1 / 2}(\partial D)]^2 \rightarrow$ $[H^{-1 / 2}(\partial D)]^2$ is given by
	$$
	\bm{\mathcal{T}}_{i m p}=\bm{\mathcal{N}}+i\eta \bm{\mathcal{I}}-i\eta \bm{\mathcal{K}}^{\prime}+i\eta \bm{\mathcal{K}}+\eta^2\bm{\mathcal{S}} .
	$$
	Our objective is to establish that $\bm{\mathcal{T}}_{imp}$ is the Fredhom operator of index zero. Following the approach used in this paper, we need to decompose the operator $\bm{\mathcal{T}}_{imp}$ into a compact operator $\bm{\mathcal{C}}$  and a coercive operator $\bm{\mathcal{M}}$,
	where
	$\bm{\mathcal{C}}=\mathrm{i}\eta \bm{\mathcal{I}}+(\bm{\mathcal{N}}-\bm{\mathcal{N}}_0)
	-\mathrm{i}\eta (\bm{\mathcal{K}}^{\prime}-\bm{\mathcal{K}}_0^{\prime})
	+\mathrm{i}\eta (\bm{\mathcal{K}}-\bm{\mathcal{K}}_0)+\eta^2\bm{\mathcal{S}}$ and $\bm{\mathcal{M}}=\bm{\mathcal{N}}_0-\mathrm{i}\eta \bm{\mathcal{K}}_0^{\prime}+\mathrm{i}\eta \bm{\mathcal{K}}_0$. The compactness of $\bm{\mathcal{C}}$ holds since $\bm{\mathcal{S}}$, $\bm{\mathcal{N}}-\bm{\mathcal{N}}_0$, $\bm{\mathcal{K}}^{\prime}-\bm{\mathcal{K}}_0^{\prime}$, and the imbedding $\bm{\mathcal{I}}$ are compact.
	However, the coercivity of the operator $\bm{\mathcal{M}}$ is not easy to obtain.
	In future work we will try to address this issue, and we would like to apply this method to some other structures, e.g., cavities and layered cavities.

	\section{Appendix}
	\begin{proposition}\label{S_kernel}
		The kernel of the operator $\bm{\mathcal{S}}$ is weakly singular.
	\end{proposition}
	\begin{proof}
		Using the recursive formula for the Hankel function
		\begin{equation*}
			H^{(1)}_1(t)=-\frac{1}{t} H^{(1)}_1(z)+H^{(1)}_0(t),
		\end{equation*}
		a straightforward calculation yields
		\begin{equation}\label{Gamma_fac}
			\boldsymbol{\Gamma}(\boldsymbol{x}, \boldsymbol{y})=\gamma_1(|\bm{x}-\bm{y}|)\mathbf{I}
			+\gamma_2(|\bm{x}-\bm{y}|) \mathbf{J}(\bm{x}-\bm{y})
		\end{equation}
		where $\mathbf{J}(\bm{x}-\bm{y})=(\bm{x}-\bm{y})
		(\bm{x}-\bm{y})^\top/|\bm{x}-\bm{y}|^2$. Here,
		the functions $\gamma_1$ and $\gamma_2$ are given by
		$$
		\begin{aligned}
			& \gamma_1(v)=\frac{\mathrm{i}}{4 \mu} H_0^{(1)}\left(k_s v\right)-\frac{\mathrm{i}}{4 \omega^2 v}\left[k_s H_1^{(1)}\left(k_s v\right)-k_p H_1^{(1)}\left(k_p v\right)\right], \\
			& \gamma_2(v)=\frac{\mathrm{i}}{4 \omega^2}\left[\frac{2 k_s}{v} H_1^{(1)}\left(k_s v\right)-k_s^2 H_0^{(1)}\left(k_s v\right)-\frac{2 k_p}{v} H_1^{(1)}\left(k_p v\right)+k_p^2 H_0^{(1)}\left(k_p v\right)\right],
		\end{aligned}
		$$
		where $v=|\bm{x}-\bm{y}|$. We use the power series for the Bessel and Hankel functions to obtain the decomposition from
		$$
		\gamma_j(v)=\frac{1}{\pi} \ln v\xi_j(v)+\chi_j(v), \quad j=1,2,
		$$
		where
		$$
		\begin{aligned}
			& \xi_1(v):=-\frac{1}{2 \mu} J_0\left(k_s v\right)+\frac{1}{2 \omega^2 v}\left[k_s J_1\left(k_s v\right)-k_p J_1\left(k_p v\right)\right], \\
			& \xi_2(v):=\frac{1}{2 \omega^2}\left[k_s^2 J_0\left(k_s v\right)-\frac{2 k_s}{v} J_1\left(k_s v\right)-k_p^2 J_0\left(k_p v\right)+\frac{2 k_p}{v} J_1\left(k_p v\right)\right].
		\end{aligned}
		$$
		Then, apply the following asymptotic behavior of the bessel functions
		\begin{align}
			J_0(t)&=1-\frac{1}{4}t^2+\mathrm{O}\left(t^4\right),\quad t\rightarrow 0,\label{J0}\\
			J_1(t)&=\frac{1}{2}t-\frac{1}{16}t^3+\mathrm{O}
			\left(t^5\right),\quad t\rightarrow 0,\label{J1}\\
			Y_0(t)&=\frac{2}{\pi}\left\{\ln \frac{t}{2}+C\right\} J_0(t)+\frac{1}{2\pi}t^2+\mathrm{O}
			\left(t^4\right),\quad t\rightarrow 0,\label{Y0}\\
			Y_1(t)&=\frac{2}{\pi}\left\{\ln \frac{t}{2}+C\right\} J_1(t)-\frac{2}{\pi}\frac{1}{t}-\frac{1}{2\pi}t
			+\frac{5}{32}t^3+\mathrm{O}
			\left(t^5\right),\quad t\rightarrow 0,\label{Y1}
		\end{align}
		yields
		\begin{align*}
			\xi_1(v)&=-\frac{1}{4 \omega^2}\left(k_s^2+k_p^2\right)+
			\frac{1}{32 \omega^2}\left(3 k_s^4+k_p^4\right)
			v^2+\mathrm{O}\left(v^4\right), \\
			\xi_2(v)&=\frac{1}{16 \omega^2}\left(k_p^4-k_s^4\right)v^2
			+\mathrm{O}\left(v^4\right),\\
			\chi_1(v)&=-\frac{1}{4 \pi \omega^2}\left[k_s^2 \ln \frac{k_s}{2}+k_p^2 \ln \frac{k_p}{2}+\frac{1}{2}\left(k_s^2-k_p^2\right)
			+\left(C-\frac{\mathrm{i} \pi}{2}\right)\left(k_s^2+k_p^2\right)\right]
			+\mathrm{O}\left(v^2\right),\\
			\chi_2(v)&=\frac{1}{4 \pi \omega^2}\left(k_s^2-k_p^2\right)+\mathrm{O}\left(v^2\right),
		\end{align*}
		for $v\rightarrow 0$, where $C=0.57721 \ldots$ is Euler's constant.
		
		It can be seen that the functions $\xi_j$ and $\chi_j$ $(j=1,2)$ are analytic functions on $\mathbb{R}\times\mathbb{R}$. Thus, the kernel of the operator $\bm{\mathcal{S}}$ is logarithmic and weakly singular.
	\end{proof}
	
	\begin{proposition}\label{K'_kernel}
		The kernel of the operator $\bm{\mathcal{K}}'$ is  strongly singular, and the kernel of the operator $\bm{\mathcal{K}}'-\bm{\mathcal{K}}_0'$ is weakly singular.
	\end{proposition}
	\begin{proof}
		For a function $f$ and a matrix $\textbf{A}$, it follow from the product rule that
		\begin{equation}\label{T_PR}
			\bm{T}_{\bm{\nu}}(f\textbf{A})=\bm{T}_{\bm{\nu}}(f \textbf{I}) \textbf{A}+f \bm{T}_{\bm{\nu}}\textbf{A},
		\end{equation}
		Applying (\ref{T_PR}) to (\ref{Gamma_fac}) yields
		\begin{equation*}
			\bm{T}_{\bm{\nu(x)}}\bm{\Gamma}(\bm{x},\bm{y})
			=\sum^{2}_{j=1}\sum^{1}_{k=0}\gamma_j^{(k)}
			(|\bm{x}-\bm{y}|)\textbf{M}^{(k)}_j(\bm{x},\bm{y})
		\end{equation*}
		where
		\begin{equation*}
			\begin{aligned}
				\gamma_j^{(0)}(v):=\frac{1}{v^2}\gamma_j(v),\quad
				\gamma_j^{(1)}(v):=\frac{1}{v}\gamma'_j(v),
			\end{aligned}
		\end{equation*}
		and
		\begin{equation*}
			\begin{aligned}
				\textbf{M}^{(0)}_1(\bm{x}, \bm{y}):&=0,
				\\
				\textbf{M}^{(0)}_2(\bm{x}, \bm{y}):&=\left[(\lambda+2 \mu) \bm{\nu}(\bm{x})(\bm{x}-\bm{y})^{\top}+\mu(\bm{x}-\bm{y}) \bm{\nu}(\bm{x})^{\top} \right.\\
				&\left.\quad+\mu \bm{\nu}(\bm{x})^{\top}(\bm{x}-\bm{y})\left(\textbf{I}-4 \textbf{J}(\bm{x}-\bm{y})\right)\right]
				,\\
				\textbf{M}^{(1)}_1(\bm{x}, \bm{y}):&=\lambda \bm{\nu}(\bm{x})(\bm{x}-\bm{y})^{\top}+\mu(\bm{x}- \bm{y}) \bm{\nu}(\bm{x})^{\top}+\mu \bm{\nu}(\bm{x})^{\top}(\bm{x}- \bm{y})\mathbf{I},\\
				\textbf{M}^{(1)}_2(\bm{x}, \bm{y}):&=\textbf{M}^{(1)}_1(\bm{x}, \bm{y})\textbf{J}(|\bm{x}-\bm{y}|).
			\end{aligned}
		\end{equation*}
		Here, the functions $\gamma_j^{(1)}$ $(j=1,2)$ are given by
		\begin{equation*}
			\begin{aligned}
				\gamma_1^{(1)}(v)= &-\frac{\mathrm{i} k_s^3}{4 \omega^2v} H_1^{(1)}\left(k_s v\right)-\frac{\mathrm{i}}{4 \omega^2}\left\{\frac{1}{v^2}\left[k_s^2 H_0^{(1)}\left(k_s v\right)-k_p^2 H_0^{(1)}\left(k_p v\right)\right]\right. \\
				& \left.-\frac{2}{v^3}\left[k_s H_1^{(1)}\left(k_s v\right)-k_p H_1^{(1)}\left(k_p v\right)\right]\right\}, \\
				\gamma_2^{(1)}(v)= & \frac{\mathrm{i}}{4 \omega^2v}\left\{k_s^3 H_1^{(1)}\left(k_s v\right)-k_p^3 H_1^{(1)}\left(k_p v\right)+\frac{2}{v}\left[k_s^2 H_0^{(1)}\left(k_s v\right)-k_p^2 H_0^{(1)}\left(k_p v\right)\right]\right. \\
				& \left.-\frac{4}{v^2}\left[k_s H_1^{(1)}\left(k_s v\right)-k_p H_1^{(1)}\left(k_p v\right)\right]\right\}.
			\end{aligned}
		\end{equation*}
		We noted that these matrices $\textbf{M}^{(k)}_{j}$ $(j=1,2,k=0,1)$ are infinitely differentiable in $\mathbb{R}^2\times\mathbb{R}^2$. Utilizing the power series for both the Bessel and Hankel functions, we can derive the decomposition from
		\begin{align*}
			\gamma^{(k)}_{j}(v)&=\frac{1}{\pi} \ln v \xi^{(k)}_{j}(v)+\chi^{(k)}_{j}(v),
		\end{align*}
		where
		\begin{equation*}
			\begin{aligned}
				\xi^{(0)}_2(v):= &\frac{1}{2 \omega^2v^2}\left[k_s^2 J_0\left(k_s v\right)-\frac{2 k_s}{v} J_1\left(k_s v\right)-k_p^2 J_0\left(k_p v\right)+\frac{2 k_p}{v} J_1\left(k_p v\right)\right],\\
				\xi^{(1)}_1(v):= & \frac{k_s^3}{2 \omega^2v} J_1\left(k_s v\right)+\frac{1}{2 \omega^2}\left\{\frac{1}{v^2}\left[k_s^2 J_0\left(k_s v\right)-k_p^2 J_0\left(k_p v\right)\right] -\frac{2}{v^3}\left[k_s J_1\left(k_s v\right)-k_p J_1\left(k_p v\right)\right]\right\},\\
				\xi^{(1)}_2(v):= & -\frac{1}{2 \omega^2v}\left\{k_s^3 J_1\left(k_s v\right)-k_p^3 J_1\left(k_p v\right)+\frac{2}{v}\left[k_s^2 J_0\left(k_s v\right)-k_p^2 J_0\left(k_p v\right)\right]\right. \\
				& \left.-\frac{4}{v^2}\left[k_s J_1\left(k_s v\right)-k_p J_1\left(k_p v\right)\right]\right\}.
			\end{aligned}
		\end{equation*}
		Applying (\ref{J0})-(\ref{Y1}) again yields
		\begin{align}
			\xi^{(0)}_2(v)&=\frac{1}{16 \omega^2}\left(k_p^4-k_s^4\right) +\mathrm{O}\left(v^2\right),\notag\\
			\xi^{(1)}_1(v)&=\frac{1}{16 \omega^2}\left(3k_s^4+k_p^4\right) +\mathrm{O}\left(v^2\right),\notag\\
			\xi^{(1)}_2(v)&=\frac{1}{8 \omega^2}\left(k_p^4-k_s^4\right) +\mathrm{O}\left(v^2\right),\notag\\
			\chi^{(0)}_{2}(v)&=-\frac{1}{16 \pi \omega^2}\left\{k_s^4 \ln \frac{k_s}{2}-k_p^4 \ln \frac{k_p}{2}+\left(C-\frac{3}{4}-\frac{\mathrm{i} \pi}{2}\right)\left(k_s^4-k_p^4\right)\right\}
			+\frac{k_s^2-k_p^2}{4 \pi \omega^2v^2}+\mathrm{O}\left(v^2\right),\label{chi3}\\
			\chi^{(1)}_{1}(v)&=\frac{1}{16 \pi \omega^2}\left\{3 k_s^4 \ln \frac{k_s}{2}+k_p^4 \ln \frac{k_p}{2}
			-\frac{5}{4} k_s^4-\frac{3}{4} k_p^4+\left(C-\frac{\mathrm{i} \pi}{2}\right)\left(3 k_s^4+k_p^4\right)\right\} \notag\\
			&\quad-\frac{k_p^2+k_s^2}
			{4\pi\omega^2v^2}+\mathrm{O}\left(v^2\right),\label{chi1}\\
			\chi^{(1)}_{2}(v)&=-\frac{1}{8 \pi \omega^2}\left\{k_s^4 \ln \frac{k_s}{2}-k_p^4 \ln \frac{k_p}{2}+\left(C-\frac{1}{4}-\frac{\mathrm{i} \pi}{2}\right)\left(k_s^4-k_p^4\right)\right\}
			+\mathrm{O}\left(v^2\right)\notag,
		\end{align}
		for $v\rightarrow0$. It can be seen the functions $\chi^{(0)}_{2}$ and  $\chi^{(1)}_{1}$  are strongly singular. Thus, the kernel of $\bm{\mathcal{K}}'$ is strongly singular.

		However, the kernel of $\bm{\mathcal{K}}'-\bm{\mathcal{K}}_0'$ can be decomposed as follows:
		\begin{equation}\notag
			\begin{aligned}
				\boldsymbol{T}_{\boldsymbol{\nu}(\boldsymbol{x})} \left[\boldsymbol{\Gamma}(\boldsymbol{x}, \boldsymbol{y})-\boldsymbol{\Gamma}_0(\boldsymbol{x}, \boldsymbol{y})\right]&=\sum_{j=1}^2\sum_{k=0}^1
				\left\{\gamma^{(k)}_j
				(|\boldsymbol{x}-\boldsymbol{y}|)-\gamma^{(k)}_{0,j}
				(|\boldsymbol{x}-\boldsymbol{y}|)\right\}
				\mathbf{M}^{(k)}_{j}(\boldsymbol{x}, \boldsymbol{y}),
			\end{aligned}
		\end{equation}
		where
		\begin{equation}\notag
			\begin{aligned}
				&\gamma_{0,1}(v)=\frac{\lambda+3 \mu}{4 \pi \mu(\lambda+2 \mu)}\ln \frac{1}{v},\quad \gamma_{0,2}(v)=\frac{\lambda+\mu}{4 \pi \mu(\lambda+2 \mu)}, \\
				&\gamma'_{0,1}(v)=-\frac{\lambda+3 \mu}{4 \pi \mu(\lambda+2 \mu)}\frac{1}{v},\quad
				\gamma'_{0,2}(v)=0,\\
				&\gamma^{(0)}_{0,j}(v):=\frac{1}{v^2}\gamma_{0,j}(v),\quad \gamma^{(1)}_{0,j}(v):=\frac{1}{v}\gamma'_{0,j}(v),\quad j=1,2.
			\end{aligned}
		\end{equation}
		This exactly eliminates the singularities in (\ref{chi3}) and (\ref{chi1}). Thus, the kernel of $\bm{\mathcal{K}}'-\bm{\mathcal{K}}_0'$ is logarithmic and weakly singular.
	\end{proof}
	
	\begin{proposition}\label{K_kernel}
		The kernel of the operator $\bm{\mathcal{K}}$ is strongly singular, and the kernel of the operator $\bm{\mathcal{K}}-\bm{\mathcal{K}}_0$ is weakly singular.
	\end{proposition}
	\begin{proof}
		The kernel of the operator $\bm{\mathcal{K}}$ can be rewritten as
		\begin{equation}\notag
			\left[\boldsymbol{T}_{\boldsymbol{\nu}(\boldsymbol{y})} \boldsymbol{\Gamma}(\boldsymbol{x}, \boldsymbol{y})\right]^\top=\sum_{j=1}^2 \sum_{k=0}^1 \gamma_j^{(k)}(|\boldsymbol{x}-\boldsymbol{y}|) \left[\mathbf{M}_j^{(k)}(\boldsymbol{y}, \boldsymbol{x})\right]^\top,
		\end{equation}
		where the functions $\gamma_j^{k}$ is the same as above
		and the matrices $\left[\mathbf{M}_j^{(k)}\right]^\top$
		are still infinitely differentiable in $\mathbb{R}^2 \times \mathbb{R}^2$. The kernel of $\bm{\mathcal{K}}$ is strongly singular since the functions $\chi_2^{(0)}$ and $\chi_1^{(1)}$ are strongly singular.
		
		Similarly, the kernel of $\bm{\mathcal{K}}-\bm{\mathcal{K}_0}$ can be can be decomposed as follows:
		\begin{equation}\notag
			\begin{aligned}
				\left\{\boldsymbol{T}_{\boldsymbol{\nu}(\boldsymbol{y})}
				\left[\boldsymbol{\Gamma}(\boldsymbol{x}, \boldsymbol{y})-\boldsymbol{\Gamma}_0(\boldsymbol{x}, \boldsymbol{y})\right]\right\}^\top=
				\sum_{j=1}^2\sum_{k=0}^1
				\left\{\gamma^{(k)}_j
				(|\boldsymbol{x}-\boldsymbol{y}|)-\gamma^{(k)}_{0,j}
				(|\boldsymbol{x}-\boldsymbol{y}|)\right\} \left[\mathbf{M}_j^{(k)}(\boldsymbol{y}, \boldsymbol{x})\right]^\top.
			\end{aligned}
		\end{equation}
		Thus, the kernel of $\bm{\mathcal{K}}-\bm{\mathcal{K}}_0$ is weakly singular since $\gamma_{0,2}^{(0)}$ and $\gamma_{0,1}^{(1)}$ eliminates the singularities in (\ref{chi3}) and (\ref{chi1}).
	\end{proof}
	
	\begin{proposition}\label{N_Ker}
		The kernel of the operator $\bm{\mathcal{N}}$ is strongly singular, and the kernel of the operator $\bm{\mathcal{N}}'-\bm{\mathcal{N}}_0$ is weakly singular.
	\end{proposition}
	\begin{proof}
		The kernel of the operator $\bm{\mathcal{N}}$ can be rewritten as
		\begin{equation}\label{N_kernel}
			\boldsymbol{T}_{\nu(\boldsymbol{x})}
			\left[\boldsymbol{T}_{\nu(\boldsymbol{y})} \boldsymbol{\Gamma}(\boldsymbol{x}, \boldsymbol{y})\right]^{\top}=
			\sum_{j=1}^2 \sum_{k=0}^1 \boldsymbol{T}_{\nu(\boldsymbol{x})}\left\{ \gamma_j^{(k)}(|\boldsymbol{x}-\boldsymbol{y}|)
			\left[\mathbf{M}_j^{(k)}(\boldsymbol{y}, \boldsymbol{x})\right]^{\top}\right\},
		\end{equation}
		Applying (\ref{T_PR}) to the right hand of (\ref{N_kernel}) yields
		\begin{equation}\notag
			\boldsymbol{T}_{\nu(\boldsymbol{x})}\left[\boldsymbol{T}_{\nu(\boldsymbol{y})} \boldsymbol{\Gamma}(\boldsymbol{x}, \boldsymbol{y})\right]^{\top}=\sum_{j=1}^2 \sum_{k=0}^2 \gamma_j^{(k)}
			(|\boldsymbol{x}-\boldsymbol{y}|)\mathbf{N}_j^{(k)}(\boldsymbol{x}, \boldsymbol{y}),
		\end{equation}
		where $\gamma^{(2)}_{j}(v):=\gamma{''}_j(v)$, $j=1,2$, and
		$\mathbf{N}_1^{(0)}(\bm{x},\bm{y})=0$,
		\begin{equation}\notag
			\begin{aligned}
				\mathbf{N}_2^{(0)}(\bm{x},\bm{y})&:=\boldsymbol{T}
				_{\boldsymbol{\nu}(\boldsymbol{x})}
				\left[\mathbf{M}_2^{(0)}(\bm{y},\bm{x})^{\top}\right]
				-\frac{2}{|\bm{x}-\bm{y}|^2}\mathbf{M}_1^{(1)}
				(\bm{x},\bm{y})\mathbf{M}_2^{(0)}(\bm{y},\bm{x})^{\top},\\
				\mathbf{N}_1^{(1)}(\bm{x},\bm{y})&:=\boldsymbol{T}
				_{\boldsymbol{\nu}(\boldsymbol{x})}
				\left[\mathbf{M}_1^{(1)}(\bm{y},\bm{x})^{\top}\right]
				-\frac{1}{|\bm{x}-\bm{y}|^2}\mathbf{M}_1^{(1)}
				(\bm{x},\bm{y})\mathbf{M}_1^{(1)}(\bm{y},\bm{x})^{\top},\\
				\mathbf{N}_1^{(2)}(\bm{x},\bm{y})&:=\frac{1}{|\bm{x}-\bm{y}|^2}\mathbf{M}_1^{(1)}
				(\bm{x},\bm{y})\mathbf{M}_1^{(1)}(\bm{y},\bm{x})^{\top},\\
				\mathbf{N}_2^{(1)}(\bm{x},\bm{y})&:=\boldsymbol{T}
				_{\boldsymbol{\nu}(\boldsymbol{x})}
				\left[\mathbf{M}_2^{(1)}(\bm{y},\bm{x})^{\top}\right]
				-\frac{1}{|\bm{x}-\bm{y}|^2}\mathbf{M}_1^{(1)}
				(\bm{x},\bm{y})\mathbf{M}_2^{(1)}(\bm{y},\bm{x})^{\top}
				+\frac{1}{|\bm{x}-\bm{y}|^2}\mathbf{M}_1^{(1)}
				(\bm{x},\bm{y})\mathbf{M}_2^{(0)}(\bm{y},\bm{x})^{\top}
				\\
				\mathbf{N}_2^{(2)}(\bm{x},\bm{y})&:=\frac{1}{|\bm{x}-\bm{y}|
					^2}\mathbf{M}_1^{(1)}
				(\bm{x},\bm{y})\mathbf{M}_2^{(1)}(\bm{y},\bm{x})^{\top}.\\
			\end{aligned}
		\end{equation}
		Here, the functions $\gamma_j^{(2)}$ $(j=1,2)$ are given by
		\begin{equation}\notag
			\begin{aligned}
				\gamma_1^{(2)}(v)= & -\frac{\mathrm{i} k_s^4}{4 \omega^2} H_0^{(1)}\left(k_s v\right)+\frac{\mathrm{i}}{4 \omega^2}\left\{\frac{1}{v}\left[2 k_s^3 H_1^{(1)}\left(k_s v\right)-k_p^3 H_1^{(1)}\left(k_p v\right)\right]\right. \\
				& +\frac{3}{v^2}\left[k_s^2 H_0^{(1)}\left(k_s v\right)-k_p^2 H_0^{(1)}\left(k_p v\right)\right]
				\left.-\frac{6}{v^3}\left[k_s H_1^{(1)}\left(k_s v\right)-k_p H_1^{(1)}\left(k_p v\right)\right]\right\}, \\
				\gamma_2^{(2)}(v)= & \frac{\mathrm{i}}{4 \omega^2}\left\{k_s^4 H_0^{(1)}\left(k_s v\right)-k_p^4 H_0^{(1)}\left(k_p v\right)-\frac{3}{v}\left[k_s^3 H_1^{(1)}\left(k_s v\right)-k_p^3 H_1^{(1)}\left(k_p v\right)\right]\right. \\
				& \left.-\frac{6}{v^2}\left[k_s^2 H_0^{(1)}\left(k_s v\right)-k_p^2 H_0^{(1)}\left(k_p v\right)\right]+\frac{12}{v^3}\left[k_s H_1^{(1)}\left(k_s v\right)-k_p H_1^{(1)}\left(k_p v\right)\right]\right\} .
			\end{aligned}
		\end{equation}
		We noted that these matrices $\mathbf{N}_j^{(k)}(j=1,2, k=0,1,2)$ are infinitely differentiable in $\mathbb{R}^2 \times \mathbb{R}^2$; see \cite{CKM00}. Once again using the the power series for the Bessel and Hankel functions to obtain the decomposition from
		$$
		\gamma_j^{(k)}(v)=\frac{1}{\pi} \ln v \xi_j^{(k)}(v)+\chi_j^{(k)}(v), \quad j=1,2,k=0,1,2,
		$$
		where
		\begin{equation}\notag
			\begin{aligned}
				\xi_1^{(2)}(v):= & \frac{k_s^4}{2 \omega^2} J_0\left(k_s v\right)-\frac{1}{2 \omega^2}\left\{\frac{1}{v}\left[2 k_s^3 J_1\left(k_s v\right)-k_p^3 J_1\left(k_p v\right)\right]\right. \\
				& \left.+\frac{3}{v^2}\left[k_s^2 J_0\left(k_s v\right)-k_p^2 J_0\left(k_p v\right)\right]-\frac{6}{v^3}\left[k_s J_1\left(k_s v\right)-k_p J_1\left(k_p v\right)\right]\right\}, \\
				\xi_2^{(2)}(v):= & -\frac{1}{2 \omega^2}\left\{k_s^4 J_0\left(k_s v\right)-k_p^4 J_0\left(k_p v\right)-\frac{3}{v}\left[k_s^3 J_1\left(k_s v\right)-k_p^3 J_1\left(k_p v\right)\right]\right. \\
				& \left.-\frac{6}{v^2}\left[k_s^2 J_0\left(k_s v\right)-k_p^2 J_0\left(k_p v\right)\right]+\frac{12}{v^3}\left[k_s J_1\left(k_s v\right)-k_p J_1\left(k_p v\right)\right]\right\}.
			\end{aligned}
		\end{equation}
		Applying (\ref{J0})-(\ref{Y1}) again yields
		\begin{align}
			\xi_1^{(2)}(v)&=\frac{1}{16 \omega^2}\left(3 k_s^4+k_p^4\right)+\mathrm{O}\left(v^2\right),\notag\\
			\xi_2^{(2)}(v)&=\frac{1}{8 \omega^2}\left(k_p^4-k_s^4\right)
			+\mathrm{O}\left(v^2\right),\notag\\
			\chi_1^{(2)}(v)&= \frac{1}{16 \pi \omega^2}\left\{3 k_s^4 \ln \frac{k_s}{2}+k_p^4 \ln \frac{k_p}{2}+\frac{7}{4} k_s^4+\frac{1}{4} k_p^4+\left(C-\frac{\mathrm{i} \pi}{2}\right)\left(3 k_s^4+k_p^4\right)\right\}
			+\frac{k_p^2+k_s^2}{4\pi\omega^2v^2}
			+\mathrm{O}\left(v^2\right),\label{chi12}\\
			\chi_2^{(2)}(v)&=-\frac{1}{8 \pi \omega^2}\left\{k_s^4 \ln \frac{k_s}{2}-k_p^4 \ln \frac{k_p}{2}+\left(C+\frac{3}{4}-\frac{\mathrm{i} \pi}{2}\right)\left(k_s^4-k_p^4\right)\right\}
			+\mathrm{O}\left(v^2\right),\notag
		\end{align}
		for $v\rightarrow0$. It can be seen the functions $\chi^{(0)}_{2}$ and $\chi^{(k)}_{1}$ $(k=1,2)$ are strongly singular. Thus, the kernel of $\bm{\mathcal{N}}$ is strongly singular.
		
		Similarly, the kernel of $\bm{\mathcal{N}}-\bm{\mathcal{N}_0}$ can be decomposed as follows:
		\begin{equation}\notag
			\boldsymbol{T}_{\nu(\boldsymbol{x})}\left\{
			\boldsymbol{T}_{\nu(\boldsymbol{y})} \left[\boldsymbol{\Gamma}(\boldsymbol{x}, \boldsymbol{y})-\boldsymbol{\Gamma}_0(\boldsymbol{x}, \boldsymbol{y})\right]\right\}^{\top}
			=\sum_{j=1}^2 \sum_{k=0}^2 \left\{\gamma_j^{(k)}(|\boldsymbol{x}-\boldsymbol{y}|)
			-\gamma_{0,j}^{(k)}(|\boldsymbol{x}-\boldsymbol{y}|)\right\}
			\mathbf{N}_j^{(k)}(\boldsymbol{x}, \boldsymbol{y}),
		\end{equation}
		where
		\begin{equation}\notag
			\gamma''_{0,1}(v)=\frac{\lambda+3 \mu}{4 \pi \mu(\lambda+2 \mu)v^2}, \quad \gamma''_{0,2}(v)=0,\quad
			\gamma_{0, j}^{(2)}(v):=\gamma''_{0, j}(v), \quad j=1,2 .
		\end{equation}
		Thus, the kernel of $\bm{\mathcal{N}}-\bm{\mathcal{N}}_0$ is weakly singular since $\gamma_{0,2}^{(0)}$ and $\gamma_{0,1}^{(k)}$, $(k=1,2)$ eliminates the singularities in (\ref{chi3}), (\ref{chi1}) and
		(\ref{chi12}).

	\end{proof}
	
\end{document}